\documentclass[reqno]{amsart}
\renewcommand\today{January 3, 1999;\quad to appear in \emph{Handbook of
Geometrical Topology}}

\begin{document}
\author[K. B. Lee and F. Raymond]{Kyung Bai Lee and Frank Raymond}
\address{University of Oklahoma, Norman, OK 73019, U.S.A.\\
         University of Michigan, Ann Arbor, MI 48109, U.S.A.}
\email{kblee@math.ou.edu}
\email{Frank.Raymond@umich.edu}
%%%%%  57S30 discrete groups of transformations
%%%%%  53C30 homogeneous manifolds
\date{\today}

\title[Seifert Manifolds]{Seifert Manifolds}
%\begin{abstract}
%\end{abstract}
\maketitle
\tableofcontents

%       Theorem environments
\newtheorem{thm}{Theorem}[subsection]
\newtheorem{cor}[thm]{Corollary}
\newtheorem{lem}[thm]{Lemma}
\newtheorem{prop}[thm]{Proposition}
\newtheorem{prob}[thm]{Problem}
\newtheorem{conj}[thm]{Conjecture}
\newtheorem{exam}[thm]{Example}
\newtheorem{defn}[thm]{Definition}
\newtheorem{rem}[thm]{Remark}
\newtheorem{caution}[thm]{Caution}
\newtheorem{ax}{Axiom}
\newtheorem{step}[thm]{}

\newcommand\beginrem{\begin{rem}\rm}
\newtheorem{notation}[thm]{Notation}

%       Math definitions
\newcommand\on[1]{{\rm #1}}
\newcommand\col{\,:\,}
\newcommand\x{\times}
\newfont{\msbm}{msbm10}
\newcommand{\rx}{\mbox{\msbm \symbol{111}}}
\newcommand\wh{\widehat}
\newcommand\wt{\widetilde}
\newcommand\bs{{\backslash}}
\newcommand\ra{\rightarrow}
\newcommand\da{\downarrow}
\newcommand\ua{\uparrow}
\newcommand\bda{\Big\downarrow}
\newcommand\bua{\Big\uparrow}
\newcommand\bdal[1]{#1\Big\downarrow\phantom{#1}}
\newcommand\bdar[1]{\phantom{#1}\Big\downarrow#1}
\newcommand\lra{\longrightarrow}
\newcommand\hra{\hookrightarrow}

\newcommand\bbc{{\Bbb C}}
\newcommand\bbh{{\Bbb H}}
\newcommand\bbq{{\Bbb Q}}
\newcommand\bbr{{\Bbb R}}
\newcommand\bbe{{\Bbb E}}
\newcommand\bbz{{\Bbb Z}}
\newcommand{\Z}[1]{\mbox{${\Bbb Z}^{#1}$}}
\newcommand{\R}[1]{\mbox{${\Bbb R}^{#1}$}}
\newcommand{\N}[1]{\mbox{${\Bbb N}^{#1}$}}
\newcommand{\tf}{torsion-free }
\newcommand{\ds}{\displaystyle}
\newcommand{\kra}{\kern-7pt\rightarrow\kern-7pt}
\newcommand{\prend}{\rule{2.3mm}{2.3mm}}
\newcommand{\arst}{\mbox{A}\mbox{$(-2,s,t)$}}

\newcommand\cala{{\mathcal A}}
\newcommand\calb{{\mathcal B}}
\newcommand\calc{{\mathcal C}}
\newcommand\cale{{\mathcal E}}
\newcommand\calf{{\mathcal F}}
\newcommand\calh{{\mathcal H}}
\newcommand\calk{{\mathcal K}}
\newcommand\calm{{\mathcal M}}
\newcommand\calo{{\mathcal O}}
\newcommand\calp{{\mathcal P}}
\newcommand\calr{{\mathcal R}}
\newcommand\cals{{\mathcal S}}
\newcommand\calt{{\mathcal T}}
\newcommand\calu{{\mathcal U}}
\newcommand\calub{{\overline{\mathcal U}}}
\newcommand\calv{{\mathcal V}}
\newcommand\calw{{\mathcal W}}
\newcommand\calz{{\mathcal Z}}

\newcommand\bh{{\bold H}}

\newcommand\Gam{{\Gamma}}
\newcommand\gam{{\Gamma}}
\newcommand\alp{{\alpha}}
\newcommand\del{{\delta}}
\newcommand\Del{{\Delta}}
\newcommand\eps{{\epsilon}}
\newcommand\lam{{\lambda}}
\newcommand\ome{{\omega}}
\newcommand\sig{{\sigma}}
\newcommand\tth{{\theta}}
\newcommand\tht{{\widetilde{\theta}}}
\newcommand\tthb{{\ol{\theta}}}
\newcommand\scirc{\kern -2pt\circ\kern -2pt}

\newcommand\aff{\on{Aff}}
\newcommand\aut{\on{Aut}}
\newcommand\inn{\on{Inn}}
\newcommand\out{\on{Out}}
\newcommand\autg{\on{Aut}(G)}
\newcommand\inng{\on{Inn}(G)}
\newcommand\outg{\on{Out}(G)}
\newcommand\isom{\on{Isom}}
\newcommand\ttop{\on{TOP}}
\newcommand\diff{\on{Diff}}
\newcommand\hol{\on{Hol}}
\newcommand\uaep{\on{UAEP}}

\newcommand\topgp{\on {TOP}_G(G\x W)}
\newcommand\difgp{\on {Diff}_G(G\x W)}
\newcommand\holgp{\on {Hol}_G(G\x W)}
\newcommand\cs{\bbc^*}
\newcommand\mw[1]{\on{M}(W,#1)}
\newcommand\cw[1]{\on{\calc}(W,#1)}
\newcommand\hw[1]{\on{\calh}(W,#1)}
\newcommand\mwg{{\on{M}}(W,G)}
\newcommand\cwg{{\on{\calc}}(W,G)}
\newcommand\mwc{\on{M}(W,C_G(K))}
\newcommand\mwcz{\on{M}(W,C_G(K)_0)}
\newcommand\mwzg{\on{M}(W,\calz(G))}
\newcommand\mwn{\on{M}(W,N_G(K))}
\newcommand\mwnz{\on{M}(W,N_G(K)_0)}
\newcommand\topw{\on{TOP}(W)}

\newcommand\pinf{P_{\infty}}
\newcommand\pone{P_{1}}
\newcommand\psl{{\on{PSL}(2,\Bbb R)}}
\newcommand\pslt{\wt{\on{PSL}}(2,\Bbb R)}
\newcommand\tp{\widetilde{P}}
\newcommand\PSL{\on{PSL}}
\newcommand\gl{\on{GL}}
\newcommand\so{\on{SO}}
\newcommand\sltr{\on{SL}(2,\Bbb R)}
\newcommand\ol[1]{\overline{#1}}
\newcommand\hhom{\on{Hom}}
\newcommand\nil{\on{Nil}}
\newcommand\opext{\on{Opext}}
\newcommand\pic{\on{Pic}}
\newcommand\jac{\on{Jac}}
\newcommand\tor{\on{Tor}}
\newcommand\ev{\on{ev}}
\newcommand\uliep{$\on{ULIEP}$}
\newcommand\nb{\overline{N}}
\newcommand\ptm{\phantom{-}}
\newcommand\inv{^{-1}}
\newcommand\sdot{\cdot}
\newcommand\ecdot{\phantom{\cdot}}
\newcommand\vp{\varphi}
\newcommand\Ad{\on {Ad}}
\newcommand\ad{\on {ad}}
\newcommand\vapt{\widetilde{\varphi}}
\newcommand{\cstar}{{c^*}}
\newcommand\gb{{\ol{G}}}
\newcommand\kb{{\ol{K}}}
\newcommand\qb{{\ol{Q}}}
\newcommand\pib{\overline{\varPi}}
\newcommand\fg{\mathfrak g}
\newcommand\fk{\mathfrak k}
\newcommand\fm{\mathfrak m}
\newcommand\fb{\overline{f}}
\newcommand\fh{\widehat{f}}
\newcommand\tb{\overline{\theta}}
\newcommand\pd{properly discontinuous}
\newcommand\llra{\enspace\lra\enspace}
\newcommand\gq{G\x_{_{\kern-4pt(f,\widetilde{\varphi})}}Q}

\newcommand{\roster}{\renewcommand{\labelenumi}{\rm (\arabic{enumi})}
\begin{enumerate}}
\newcommand{\rosternb}{\renewcommand{\labelenumi}{\theenumi.}\begin{enumerate}}

\newcommand\mmm[1]{\par\noindent\hskip40pt\(#1\hfil\)\noindent}
\hfuzz2pc % Don't bother to report overfull boxes if overage is < 1pc
\vfuzz2pc % Don't bother to report overfull boxes if overage is < 2pc
%---end of macros    ----------------------------------------------------

\section{Introduction}
Seifert fibered spaces have played an interesting and important role in
topology and transformation groups.  The classical definition of Seifert
\cite{seif33-1} in 1933, concerning circle bundles with specified 
singularities over 2-dimensional manifolds, preceded much of the notions 
and development of bundle theory.  
In fact, Seifert's splendid paper had a profound influence on
both bundle theory and low dimensional topology.  Recognition of the importance
of Seifert fiber spaces to transformation groups was communicated a long time
ago to one of the authors by Deane Montgomery.

Many of the interesting spaces in topology and geometry admit the structure of
a bundle with singularities.  The geometric structure is usually displayed as
part of this bundle structure.  While Seifert fiberings, in their utmost
generality, could be said to encompass all of these phenomena, we have
concentrated here, in our exposition, on a particular type called {\it injective
Seifert fiberings.}  These seem to entail the simplest structure in that they
admit a global uniformization, have close connections with transformation
groups, and have many applications.  Moreover, because topologists have become
intrigued with the ``geometrization'' of manifolds, we hope that this paper
will also inform the reader of the significant role that Seifert fiber spaces
and Seifert constructions can play in these matters.

We have tried to make our exposition as concrete and accessible as
possible. 
This is especially true in the beginning and in the last two sections where
examples and illustrations of the theory are worked out in detail. 
Each section has its own introduction and it will be helpful to read 
section 2, 5 and 6 first.
To digest sections 3 and 4, it is advisable to read these through first for an
overview before checking the details.
\bigskip

%\subsection{Notation} 
\noindent{\bf Notation.} 
We shall use the following notation throughout the rest of the paper.
For a topological group $G$ with $K$ a closed subgroup and $X$ a nice
topological space, we put\par
\smallskip\noindent
\vbox{\halign{\hskip0.3in\hfil#&#\hfil&\quad#\hfil\cr
$\aut(G)$  &$={\rm Continuous\ automorphisms\ of\ } G$\cr
$\inn(G)$  &$={\rm Inner\ automorphisms\ of\ } G$\cr
$\out(G)$  &$=\aut(G)/\inn(G)$\cr
$N_G(K)$   &$=$ Normalizer of $K$ in $G$\cr
$C_G(K)$   &$=$ Centralizer of $K$ in $G$\cr
$\aut(G,K)$  &$=\{\alpha\in \aut(G) : \alpha|_K\in\aut(K)\}$\cr
$\inn(G,K)$  &$=\inn(G)\cap\aut(G,K)$\cr
$\out(G,K)$  &$=\aut(G,K)/\inn(G,K)$\cr
$\mu(a) $ &$=$ Conjugation by  $a$;\quad  so, 
$\mu(a)(x) = a x a^{-1}$ for $x\in G$.\cr
$\ttop(X)$  &$=$ The\ group of self-homeomorphisms of  $X$\cr
$\diff(X)$  &$=$ The\ group of self-diffeomorphisms of a smooth manifold $X$\cr
}}
%---end of section 1 ----------------------------------------------------

\section{Definitions, Motivation and Examples}

\subsection{Rough Definition of Seifert Fiber Spaces}  
In section 2, we explain what we mean by Seifert fiber spaces, give motivations,
precise definitions and some examples.

The Seifert fiber space construction was first defined and studied in
\cite{cr69-1} and \cite{cr71-3} and later reformulated in \cite{lee83-1},
\cite{lr84-1} and \cite{lr89-1}.
It is really a generalization of the classical Seifert spaces and homogeneous
spaces of a Lie group by its uniform lattices.  We shall start with a rough
definition, continue with some motivation, and finally, give an exact
definition.  

We then examine what one must do to effectively construct all possible
injective Seifert fiberings over a given base with a given typical fiber.
As a result of this investigation, we determine the precise nature of the
singular fibers.
\medskip

A Seifert fibering $F\ra E\stackrel{p}{\ra} B$ is a ``fibering'' with 
singularities.  For each $b\in B$, there is associated a finite subgroup 
$Q_b$ of $\ttop(F)$, the
group of self-homeomorphisms of $F$, so that $p^{-1}(b)=F/Q_b$.  If $Q_b$ is
trivial, $p^{-1}(b)=F$ is called a typical fiber.  Otherwise $F/Q_b$ is 
called a singular fiber.  
For example, if $F=T^k$ is a $k$-torus, and $Q_b$ acts freely on
$F$, then $F/Q_b$ will be a flat manifold which is finitely covered by 
the torus.
We shall require the groups $Q_b$ to be
controlled in a particular way as a replacement of the local triviality 
condition.  
Before giving a formal definition, we start with some examples.
\smallskip

\begin{exam}
{\rm
%\begin{figure}
%\centerline{\psfig{figure=mobius.eps,height=2in,width=2in}}
%\caption{1 M\"obius band as a Seifert fiber space}
%\end{figure}
As a simple example, consider the M\"obius band
$$
E=I\x [-1,1]/(0,t)\sim(1,-t).
$$
It has an obvious circle action so that its orbit space
$B=S^1\bs E$ is an interval.
Thus, we get $S^1\ra E\stackrel{p}\ra B$.
The orbits are circles parallel to the boundary
circle. The center circle orbit however is doubly covered by each of the
other orbits. $B$ is an arc $[0,1]$.
Here $p\inv(t)=S^1$, a free orbit for every $t\neq 0$, while
$p\inv(0)=\bbz_2\bs S^1$, the center circle, is the only singular fiber. 
}
\end{exam}

\begin{exam}  
{\rm
Let $\alpha: T^2\ra T^2$ be the map $\alpha(z_1,z_2)=(\ol{z}_1,\ol{z}_2)$ 
on a 2-torus.
Then $\alpha$ has period 2, and has 4 fixed points.
In fact, $B=T^2/{\bbz_2}$ is a sphere with 4 branch points.
Take the product $T^2\x I$, a $2$-torus and the unit interval. 
Identify the two ends by
\[
(w,0)\leftrightarrow (\alpha(w),1)
\]
to form the mapping torus $M$ of $\alpha$.
Denote an element of $M$ by $[w,t]$, and define an $S^1$-action on $M$ by
\[
e^{2\pi i\theta}([w,t])=[w,t+2\theta],
\]
where the second coordinate is taken modulo 1 and compatible with the identifications.   
The space $E=(T^2\x I)/\bbz_2$ is homeomorphic to an orientable 
$3$-dimensional flat Riemannian manifold with holonomy group $\bbz_2$.
The circle $S^1$ acts as isometries on $E$. If $p$ denotes the orbit
mapping, then the orbit space is $B=T^2/\bbz_2$.
The free orbits, away from the 4 branch points, are typical fibers while the
doubly covered orbits over the branch branch points are singular fibers.
These singular fibers are still circles, but are half of the length of the
typical fibers.
}
\end{exam}

\subsection{Torus Actions}
\label{torus-action}
To facilitate our understanding of Seifert fiber spaces we
shall examine first the situation arising from group actions.  
% New fro FAX
We shall assume familiarity with elementary properties of group actions.
Good references as we shall use these concepts are the book by Bredon
\cite{bred72-1}, the paper \cite{cr69-1}, 
and the chapter on {\it Topological Transformation Groups} of this book 
by Adem and Davis.
We recall first several definitions and facts.

Our spaces $X$ are path-connected, completely regular and Hausdorff 
so that the various slice theorems are valid. 
For covering space theory, we need and tacitly assume that our spaces are 
locally path-connected and semi 1-connected. 

A topological group $G$ acts {\it properly} on $X$ if the closure of 
$\{g\in G\mid g C\cap C\neq\emptyset\}$ is a compact subset of $G$, for all
compact subsets  $C$ of $X$.
See \cite{pala61-1} for properties of proper actions.
In particular, a proper action of a (not necessarily connected) 
locally compact Lie group on a completely regular space admits a slice, 
and the orbit space is Hausdorff.

If $Q$ is a discrete group acting properly on $X$, we shall say that $Q$
acts {\it \pd ly} on $X$. Note the isotropy subgroups are finite for a \pd\
action.

Let $G$ be a path-connected group acting on a path-connected Hausdorff 
space $X$.
Fix $x\in X$. The evaluation map $\ev^x\col(G,e)\ra(X,x)$ is defined by
\[\ev^x(g)=g\cdot x\,.\]
This induces the evaluation homomorphism
\[
\ev^x_*\col\pi_1(G,e)\ra\pi_1(X,x)\,.
\]
Put $H={\rm image}$ of $\ev^x_*\subset\pi_1(X,x)$.  Then $H$ is a central
subgroup of $\pi_1(X,x)$ which is independent of choices of $x$
\cite[Lemma 4.2]{cr69-1}.

Let $K$ be a normal subgroup of $\pi_1(X,x)$ and 
\par
\centerline{$X_K$ the covering space of $X$ associated with the subgroup $K$}
\par\noindent
so that $\pi_1(X_K,x')=K$.
Assume that $H\subset K$ and put
$$
Q=\pi_1(X,x)/\pi_1(X_K,x')=\pi_1(X,x)/K.
$$
Then we have the

\begin{prop}
[{\rm \cite[\S4]{cr69-1} and \cite[Chapter I, \S9]{bred72-1}}]
\label{proper}
\roster
\item
The $G$-action on $X$ lifts to an action of $X_K$ which commutes with the 
covering $Q$-action on $X_K$. The lifted action of $G$ on $X_K$ is 
proper if the action of $G$ on $X$ is proper.  
\item
The induced $Q$ action on the quotient space $W=G\bs X_K$ is properly
discontinuous $($i.e. behaves at least like branched coverings with Hausdorff
quotient$)$ provided $G$ acts properly on $X$.
\end{enumerate}
\end{prop}
\medskip

\begin{defn}
{\rm
A torus action $(T^k,X)$ is said to be {\it injective} if\linebreak 
$\ev^x_*\col \pi_1(T^k,e)\ra \pi_1(X,x)$ is injective.
}
\end{defn}

Suppose $T^k$ acts injectively on a path connected, locally path connected,
semi $1$-connected paracompact Hausdorff space $X$.  
Then, by the above proposition,
the $T^k$ action lifts to an action on $X_{\bbz^k}$, where $\bbz^k$ is the image of 
$\pi_1(T^k)\ra\pi_1(M)$.
We show that the lifted action $(T^k,X_{\bbz^k})$ is free and splits.
\medskip

\begin{prop}
[{\cite[\S3.1]{cr69-1}}]
\label{injective-torus}
If $T^k$ acts injectively on $X$, then
$X_{\bbz^k}$ splits into $T^k\times W$ so that $(T^k,X_{\bbz^k})=(T^k,T^k\times
W)$, where the $T^k$ action on $T^k\times W$ is via translation on the first
factor and trivial on the simply connected $W$ factor.
\end{prop}

\begin{proof}
Let $G=T^k$ and $S^1$ be a subgroup of $G$ isomorphic to the circle.
(Note that $S^1$ is a direct factor).
Lift the $S^1$ action to $X'=X_{\bbz^k}$. Let $y'\in X'$ and suppose
$S^1_{y'}\neq 1$ is the stabilizer of the lifted $S^1$ action at $y'$.
Choose paths $\wt{\gamma}$ in $X'$ from the base point $x'$ over $x$ to
$y'$ and $\alpha: (I,0,1)\ra (S^1,1,z)$ where $z$ is the ``first'' element
$\neq 1$ of $S^1$ for which $z\cdot y'=y'$.
Then $\alpha(t) y'$ defines a loop at $y'$; and $\gamma*\alpha*\ol{\gamma}$
is the associated loop based at $x'$.
Now $(\gamma*\alpha*\ol{\gamma})^n\sim \gamma*\alpha^n*\ol{\gamma}$
represents a generator of $\bbz^k=H$. Hence $n$ is the order of $z$ in
$S^1$.
This implies that $\gamma*\alpha*\ol{\gamma}$ represents an $n$th root of a
generator of $\bbz^k$ which is impossible unless $n=1$.
Thus $S^1_{y'}=1$, for all $y'\in X'$ and all circle subgroups of $S^1$ of
$T^k$. Hence the lifted toral action $(T^k,X')$ must be a free action and
$X\ra W=T^k\bs X'$ is a principal $T^k$ bundle.

We now show that this bundle is trivial.
Let $f: W\ra B_{T^k}$ be the classifying map from $W$ into the classifying
space for principal $T^k$-bundles.
The bundle $X'\ra W$ is represented by the homotopy class $[f]\in
[W,B_{T^k}]$.
Since $B_{T^k}$ is a $K(\bbz^k,2)$ space,  
$[f]$ is represented by the
image $f^*(u)$, where $u\in H^2(B_{T^k};\bbz^k)$ represents
$\on{id}: B_{T^k}\ra B_{T^k}$.

The cohomology sequence,
$$
0 \lra H^1(W;\bbz^k) \lra H^1(X';\bbz^k) \lra H^1(T^k;\bbz^k)  
  \lra H^2(W;\bbz^k) \lra H^2(X';\bbz^k),
$$
which arises from terms of low degree of the spectral sequence of the
fibering $X'\ra W$, is exact (e.g., \cite[p. 332 Ex 2]{mac:book}, or
\cite[p. 329 C]{cart-eil56-1}).

We observe that $H^1(X';\bbz^k) \lra H^1(T^k;\bbz^k)$ is a surjective
isomorphism because $\pi_1(T^k,1)\ra\pi_1(X',x')$ is an isomorphism.
Consequently, $f^*(u)$ injects into $H^2(X';\bbz^k)$.
On the other hand, $\pi^*(f^*(u))=\wt{f}^*({\pi'}^{*}(u))$, where
$\pi': E_{T^k}\ra B_{T^k}$, the universal $T^k$-bundle, and $\wt{f}$ is the
bundle map $X'\ra E_{T^k}$ induced by $f$.
But ${\pi'}^{*}(u)\in H^2(E_{T^k};\bbz^k)=0$. Therefore, $\pi^*(f^*(u))=0$,
which implies $f^*(u)=0$. Hence, $X'\ra W$ is a trivial fibering and
$X'=T^k\x W$.
\medskip

Finally, we remark that the group $T^k\x Q$ is acting properly on $T^k\x
W$, and $W$ is simply connected. The $T^k$ action does not lift to the
universal covering $\wt{X}$ of $X$ but the induced ineffective $\bbr^k$
action on $T^k\x W$ lifts to an effective $\bbr^k$ action on 
$\bbr^k\x W$ and commutes with the group $\pi_1(X,x)$ of covering
transformations on $\bbr^k\x W$.
\end{proof}

\subsection{Examples of Injective Actions}  
A manifold $M^n$ is called {\it hyper-aspherical} \cite{ds82-1} if there 
is a degree $1$-map from $M$ to a closed aspherical manifold of dimension $n$.

\begin{prop}
Any effective torus action on the following closed connected manifolds is 
injective:  
\roster
\item
All aspherical manifolds {\rm\cite{cr69-1}}, 
\item
All hyper-aspherical manifolds {\rm\cite{ds82-1}}, 
\item
More generally, all manifolds $M^n$ with 
$\xi^*: H^n(K(\pi_1(M),1);\bbq)\ra H^n(M;\bbq)$ surjective, 
where $\xi^*$ is induced by the classifying map $\xi: M\ra K(\pi_1(M),1)$,
with $\pi_1(M)$ torsion free. (The torsion freeness assumption of $\pi_1(M)$ can
be dropped if we assume the action is smooth, {\rm cf \cite{bh82-1} and 
\cite{lr87-1}.})
See {\rm\cite{sy79-1}, \cite{bh82-1}, \cite{ds82-1}, \cite{kk83-1}, 
\cite{ww83-1}}, {\rm\cite{glo85-1}} and {\rm\cite{lr87-1}} 
for such generalizations.
\item
Homologically K\"ahler manifolds all of whose isotropy subgroups are finite
(e.g., holomorphic actions) {\rm\cite[p. 170]{bore60-1}, 
\cite[p. 186]{cr71-3}}.
\end{enumerate}
\end{prop}

If every circle action on a closed manifold $M$ is injective, the only compact
connected Lie groups that can act on $M$ are tori.  See \cite{cr69-1} for 
a proof of this fact.  Therefore tori are the only compact connected Lie 
groups that act effectively on those manifolds listed in (1), (2) and (3).
Moreover, in each of these cases, $\pi_1(T^k)$ must inject into the center
of $\pi_1(M)$ if $T^k$ acts effectively.
Consequently, the only compact Lie groups that can act effectively on these
$M$ for which the center of $\pi_1(M)$ is finite are finite groups. 
\medskip

\subsection{Injective Seifert fibering}
\label{inj-seif-fib}
Perhaps the most important feature of an injective torus action
$(T^k,X)$ is the splitting $(T^k,T^k\times W)$ in Proposition
\ref{injective-torus}.
The universal covering group $\bbr^k$ acts ineffectively via $\bbr^k\ra T^k$
on $T^k\x W$ and lifts to the action of $\bbr^k$ as left translations on
$\bbr^k\x W$, the universal covering space of $T^k\x W$ and of $X$.
This $\bbr^k$ action contains the $\bbz^k$ covering transformations over
$T^k\x W$ and commutes with the group $\pi=\pi_1(X)$ of covering
transformations on $X$.
Together $\bbr^k$ and $\pi$ generate a subgroup $\bbr^k\cdot\pi$ of
$\ttop(\bbr^k\x W)$ isomorphic to $\bbr^k\x_{\bbz^k}\pi$, since
$\bbz^k=\bbr^k\cap\pi$.

We have been describing these features of injective actions to motivate our
definition of injective Seifert fiberings. It is apparent that we could
start with $\bbr^k\x W$ and using the reverse procedure, reconstruct the
injective toral actions on $X$. We will now describe this reverse procedure. 
We will also put it in a more general context where $\bbr^k$ is replaced 
by any Lie group $G$ and $W$ is any completely regular space admitting 
covering space theory. 

Let $\pi$ be a discrete subgroup of $\ttop(G\x W)$ and assume
\roster
\item[1.]
The left translational action $\ell(G)$ of $G$ on $G\x W$ is 
\emph{normalized} by $\pi$. 
\item[2.]
$\gam\subset \pi\cap\ell(G)$ is discrete and normal in $\pi$. 
\item[3.]
The induced action of $Q=\pi/\gam$ on $W$ is proper.
\end{enumerate}
These conditions imply that the group $\ell(G)\cdot\pi\subset\ttop(G\x W)$,
generated by $G$ and $\pi$, acts properly on $G\x W$.
\bigskip

With the above three conditions, we obtain the commutative diagram:
$$
\begin{array}{cccccl}
G              &\lra &G\x W         &\stackrel{G\bs}\lra  &W  
&\qquad\textrm{Product principal\ $G$-fibering}\\
\bdar{\gam\bs} &     &\bdar{\pi\bs} &      &\bdar{Q\bs}&\\
\gam\bs G      &\lra &X=\pi\bs(G\x W) &\stackrel{p}\lra  &Q\bs W
&\qquad\textrm{Seifert\ fibering}
\end{array}
$$
$X$ is called an {\it injective Seifert fiber space}, with {\it typical fiber}
$\Gam\bs G$, $B$ is called the {\it base} and the mapping $p$ is called
the {\it injective Seifert fibering}.  Since the actions of $\pi$ on $G\x W$ and
the action of $Q$ on $W$ are properly discontinuous, the quotient spaces $X$ and
$B$ are reasonable spaces.

If $\pi$ centralizes $\ell(G)$, then $\gam\subset\pi\cap\ell(G)$ is in the
center of $G$ and the $G$ action descends to a $G/\gam$ action on $X$.
In general, neither $\Gamma$ is central in $G$ nor $\pi$
centralizes the $G$ action but, as we shall see,
the obvious ``fibers'' (inverse images
$p\inv(b)$) still have a very nice description.

\subsection{$\topgp$}
\label{topgp}
To describe or construct all possible injective Seifert fiberings over $B$
with typical fiber $\Gam\bs G$, we begin with a proper action of a discrete
group $Q$ on $W$ so that $B=Q\bs W$. We form the product space $G\times W$
letting $G$ act on the first factor of $G\times W$ by left translations.
Denote this action by $\ell(G)\subset\ttop(G\times W)$.
We then search for discrete group $\pi$ contained in the normalizer of
$\ell(G)$ in $\ttop(G\times W)$ such that $\gam\subset\pi\cap G$ is discrete in
$G$ and the induced action of the quotient $\pi/\gam$ is equivalent to the
$Q$ action on $W$.
Our first task then is to describe the algebraic structure of the
normalizer of $\ell(G)$ in $\ttop(G\times W)$.
The topological structure of $\ttop(G\times W)$, 
while important and significant for some
applications, is not necessary for our problem now at hand.
\bigskip

A space is {\it admissible} if it is completely regular (Hausdorff)
locally path-connected and semi-locally simply-connected.
Let $G$ be a Lie group, $\aut(G)$ the Lie group of continuous
automorphisms of $G$. $\ttop(X)$ denotes the group of homeomorphisms of $X$.
See section 1 for a list of notation.
\bigskip

Recall that a homeomorphism $f$ of $G\x W$ onto itself is
{\it weakly $G$-equivariant} if and only if there exists a continuous 
automorphism  $\alpha_f$ of $G$ so that
\[
f(a\sdot x,w) = \alpha_f(a)f(x,w)
\]
for all $a\in G$ and $(x,w)\in G\x W.$

\begin{defn}{\rm
The group of all weakly $G$-equivariant self-homeomorphisms of $G\x W$ 
is denoted by $\topgp$.
}\end{defn}

\begin{lem}
\label{pre-universal-group}
$\topgp$ is the  normalizer of $\ell(G)$ in $\ttop(G\x W)$.
\end{lem}

\begin{proof}
Let $f\in\ttop(G\x W)$. Then the following are all equivalent:
\roster
\item
$f\in\topgp$
\item
$f$ is weakly $G$-equivariant.
\item
There exists $\alpha_f\in\aut(G)$ such that 
$f(a\sdot x,w) = \alpha_f(a)f(x,w)$ for all $a\in G$ and $(x,w)\in G\x W$.
\item
There exists $\alpha_f\in\aut(G)$ such that 
$f\circ\ell(a)\circ f\inv=\ell(\alpha_f(a))$ for all $a\in G$.
\item
$f$ normalizes $\ell(G)$.
\end{enumerate}
\end{proof}

Each element $f\in\topgp$ sends fibers of $G\x W\ra W$ to fibers.
That is, each $f$ induces a map $\ol{f}\in\ttop(W)$ so that
$$
\begin{array}{ccccccccc}
G\times W &\stackrel{f}\lra &G\times W\\
\Big\da &&\Big\da\\
W &\stackrel{\overline f}\lra &W
\end{array}
$$
commutes. Therefore, to each $f\in\topgp$, we can assign a pair
$(\alpha_f,\ol{f})\in\aut(G)\x\ttop(W)$.
This assignment
$$
\topgp\lra \aut(G)\x\ttop(W)
$$
is a surjective
homomorphism and it splits.
For, if we define $(\alpha,h)(x,w)=(\alpha(x),h(w))$ for $(\alpha,h)\in
\aut(G)\x\ttop(W)$, then
\[
\begin{array}{ll}
(\alpha,h)(a\cdot x,w)
&=(\alpha(a\cdot x),h(w))\\
&=(\alpha(a)(\alpha(x),h(w))\\
&=\ell(\alpha(a))((\alpha(x),h(w))\\
&=\ell(\alpha(a))(\alpha,h)(x,w).
\end{array}
\]
Let $K$ be the kernel of $\topgp\ra \aut(G)\x\ttop(W)$.
If $f\in K$, since $f$ moves only along the fibers,
there exists a unique continuous function
$\lambda : G\x W\ra G$ such that $f(x,w)=(x\cdot(\lambda(x,w))\inv,w)$.
We show that $\lambda$ only depends upon $W$.
For any $a\in G$,
\[
\begin{array}{ll}
(ax\cdot(\lambda(ax,w)\inv,w)
&=f(a\cdot x,w)\\
&=\ell(a) f(x,w)\quad{\rm since\ } \alpha_f={\rm id}\\
&=\ell(a)(x\cdot(\lambda(x,w)\inv,w)\\
&=(ax\cdot(\lambda(x,w)\inv,w).
\end{array}
\]
So, $\lambda(ax,w)=\lambda(x,w)$ which means $\lambda$ only depends upon
$W$.
Let
$$
\mwg=\{{\rm
the\ continuous\ maps\ of\ } W\ {\rm into\ } G\}.
$$
For each $\lambda\in\mwg$, define $f_{\lambda}\in K$ by
$f_{\lambda}(x,w)=(x\cdot(\lambda(w))\inv,w)$.
The assignment $\mwg\ni\lambda\mapsto f_{\lambda}\in K$ is easily checked
to be an isomorphism, where the group law in $\mwg$ is   
\[
(\lambda_1\cdot\lambda_2)(w)=\lambda_1(w)\cdot\lambda_2(w).
\]
Conjugation in $\topgp$ induces the action of $\aut(G)\x\ttop(W)$ on
$\mwg$.
It is given by
\[
^{(\alpha,h)}\lambda=\alpha\circ\lambda\circ h\inv,
\]
that is,
\[
(\alpha,h)\circ\lambda\circ(\alpha,h)\inv(x,w)
=(x\cdot\alpha(\lambda(h\inv(w)))\inv,w).
\]
For,
\[
\begin{array}{ll}
(\alpha,h)\circ\lambda\circ(\alpha,h)\inv(x,w)
&=(\alpha,h)\circ\lambda(\alpha\inv(x),h\inv(w))\\
&=(\alpha,h)(\alpha\inv(x)\cdot(\lambda(h\inv(w)))\inv,h\inv(w))\\
&=(x\cdot\alpha(\lambda(h\inv(w)))\inv,w)
\end{array}
\]
and
\[
\begin{array}{ll}
(^{(\alpha,h)}\lambda)(x,w)
&=\alpha\circ\lambda\circ h\inv(x,w)\\
&=(x\cdot\alpha(\lambda(h\inv(w)))\inv,w).
\end{array}
\]
Therefore we have shown

\begin{lem} %(\cite{{lr81-1}})
\label{universal-group}
$\displaystyle{
\topgp = \mwg \rx (\autg \x \topw),
}$ and it is the  normalizer of $\ell(G)$ in $\ttop(G\x W)$.
\end{lem}

The group operation in $\topgp$ is
\[
\begin{array}{ll}
(\lambda_1,\alpha_1,h_1)\cdot(\lambda_2,\alpha_2,h_2)
&=(\lambda_1\cdot ^{(\alpha_1,h_1)}\lambda_2,
\alpha_1\scirc\alpha_2,h_1\scirc h_2)\\
&=(\lambda_1\cdot (\alpha_1\scirc\lambda_2\scirc h_1\inv),
\alpha_1\scirc \alpha_2,h_1\scirc h_2)
\end{array}
\]
Specifically, $(\lambda, \alpha, h)$ acts on $(x,w)$ by
\[
\begin{array}{ll}
(\lambda, \alpha, h) \cdot (x,w) 
&= (\lambda,1,1)\circ(1,\alpha,h))(x,w)\\
&= (\lambda,1,1)(\alpha(x),h(w))\\
&= (\alpha(x)\cdot(\lambda(h(w)))\inv,h(w)).
\end{array}
\]
Then $\mwg\rx\aut(G)$ is the group of all weakly $G$-equivariant
homeomorphisms of $G\x W$ which move only along the fibers.
The group $\mwg$ is the gauge group for the trivial fiber bundle 
$G\x W\ra W$.

For $a\in G$, the constant map $W\ra G$ sending $W$ to $a$ is denoted by
$r(a)$.
Clearly,
\[
r(a)=(a,1,1)\in  \mwg \rx (\autg \x \topw).
\]
This is a right translation by $a^{-1}$ on the first factor of $G\x W$
so that $r(a)(x,w)=(x\sdot a^{-1},w)$, and the subgroup of all such right
translations is denoted by $r(G)\subset \mwg$.
Let $\ell(G)$ denote the group of left translations on the first
factor so that $\ell(a)(x,w)=(a\sdot x,w)$.
Then elements of $\ell(G)$ are of the form
\[
\ell(a)=(a^{-1}, \mu(a),1)\in  \mwg\rx(\autg \x \topw),
\]
where $\mu(a)\in\inng$ is conjugation by $a$.

\begin{rem}{\rm
If $W$ is a smooth manifold and we take $\topgp\cap\diff(G\x W)$, then this
coincides with the weakly $G$-equivariant diffeomorphisms of $G\x W$ and
$$
\diff_G(G\x W)=\calc(W,G)\rx(\aut(G)\x\diff(W)),
$$
where $\calc(W,G)$ is the group of smooth maps of $W$ into $G$.
}\end{rem}

\subsection{Example: $3$-dimensional Seifert manifolds with base the
2-torus}
\label{3-dim-example}
We take $W=\bbr^2$. A group $Q=\bbz^2$ acts on $W$ as translations.
Let 
$$1\ra \bbz\ra\pi\ra Q\ra 1$$
be a central extension of $\bbz$ by $Q$. 
Then $\pi$ has a presentation
$$
\pi=\langle \alpha,\beta,\gamma\;\;|\;\; [\alpha,\beta]=\gamma^p, 
[\alpha,\gamma]=[\beta,\gamma]=1\rangle,
$$
where $\gamma$ is a generator of the center $\bbz$, and the images of
$\alpha,\beta$ in $Q$ are generators of $Q$.
Suppose $p\neq 0$.
Using $\bbz\subset\bbr$, one can obtain an effective action of $\pi$ on
the product $\bbr\x W$ as follows:
For $(z,x,y)\in\bbr\x W$,
\newcommand\eee{\;=\;}
\begin{equation}
\begin{array}{cclll}  
\alpha(z,x,y)&\eee&(z+y,&x+1,&y)\\
\beta (z,x,y)&\eee&(z,&x,&y+1)\\
\gamma(z,x,y)&=&(z+{1\over p},&x,&y).\\
\end{array}
\end{equation}
Notice that these maps are of the form
\[
(z,x,y)\mapsto (\phi(z)-\lambda(h(x,y)),h(x,y)),
\]
where $\phi$ is an automorphism of $\bbr$ and $h$ is an action of $Q$ on
$W$, and $\lambda$ is a map $W\ra\bbr$.
In our case $\phi$ is always the identity automorphism.
Consequently, the group $\pi$ lies in $\ttop_{\bbr}(\bbr\x W)$ as
$$
(\lambda,\phi,h)\; \in\; \on{M}(W,\bbr)\rx(\gl(1,\bbr)\x\ttop(W)).
$$
The action of $\pi$ on $\bbr\x W$ described above can be explained
differently as follows. Consider the Heisenberg group
$$
N=\left\{
\left[
\begin{matrix}
1 &\;\;\; x &\;\;\; z \\
0 &\;\;\; 1 &\;\;\; y \\
0 &\;\;\; 0 &\;\;\; 1
\end{matrix}
\right]
\; :\; x,y,z\in\bbr\right\}
$$
which is connected, simply connected and two-step nilpotent.
We denote such a matrix by $(z,x,y)$ so that
$$
\left[
\begin{matrix}
1 &\;\;\; x &\;\;\; z \\
0 &\;\;\; 1 &\;\;\; y \\
0 &\;\;\; 0 &\;\;\; 1
\end{matrix}
\right]
\quad\longleftrightarrow\quad (z,x,y).
$$
Then the group operation is
\[
(z',x',y')\cdot(z,x,y)=(z'+z+x'y,\ x'+x,\ y'+y),
\]
and the center of $N$ is 1-dimensional $\calz=\bbr$, consisting of all 
matrices with $x=y=0$. The quotient $W=N/\calz$ is isomorphic to $\bbr^2$
so that 
\[
1\ra\bbr\ra N\ra\bbr^2=W\ra 1
\]
is an exact sequence of groups. As spaces, this is a smooth fibration which
is also a product $N=\bbr\x W$.
Suppose $p\neq 0$ in the presentation of $\pi$.
Let
$$
\alpha=(0,1,0),\; \beta=(0,0,1),\;{\;\rm and\;} \gamma=({1/p},0,0)\in N.
$$
Then these satisfy the relations of the group $\pi$ so that $\pi$ sits in
$N$ as a discrete subgroup,
and furthermore, the action of $\pi$ defined above is nothing but the left
multiplication on $N$ by elements of $\pi$.

Actually this means that $\pi\bs N$ is a nilmanifold which is a matter that
we shall explore more later. The free action of $\pi$ commutes with the
$\bbr$ action on $\bbr\x\bbr^2$ and $\bbr\cap\pi$ is the central $\bbz$
subgroup generated by $\gamma=({1\over p},0,0)$.
The $\bbr$ action on $\bbr\x\bbr^2$ descends to an effective
$S^1=\langle\gamma\rangle\bs\bbr$ action on $\pi\bs(\bbr\x\bbr^2)$.
The map $\pi\bs(\bbr\x\bbr^2)\ra \bbz^2\bs\bbr^2=T^2$ is the orbit mapping 
of a free $S^1$ action. 
It is the principal $S^1$-bundle over $T^2$ with Euler class $-p$.

\subsection{Injective Seifert Constructions}
\label{inj-seif-const}
In section \ref{topgp}, we described the structure of the group of weakly
$G$-equivariant   self homeomorphisms of $G\x W$ as
$$
\topgp=\mwg\rx(\aut(G)\x\ttop(W))
$$
and showed that it was the normalizer of $\ell(G)$ in $\ttop(G\x W)$.
The action of $\mwg\rx(\aut(G)\x\ttop(W))$ on $G\x W$ is given by
\[
(\lambda, \alpha, h) \cdot (x,w) = (\alpha(x)\cdot(\lambda(hw))\inv,hw).
\]
We observed that $\ell(G)$ embeds in $\topgp$ as
\[
\ell_a\mapsto (a\inv,\mu(a),1)\in\mwg\rx\inn(G)\subset\topgp,
\]
where $\inn(G)$ denotes the inner automorphisms of $G$.
Both $\inn(G)$ and
\linebreak  
$\mwg\rx\inn(G)$ are normal in $\mwg\rx\aut(G)$ and in $\topgp$, respectively.
We may therefore rewrite $\topgp$ as the group extension
\[
1\ra \mwg\rx\inn(G) \ra \topgp\ra\out(G)\x\ttop(W)\ra 1
\]
where $\out(G)=\aut(G)/\inn(G)$.

In section \ref{inj-seif-fib}, we described the condition that must be 
satisfied if a discrete subgroup $\pi$ of $\topgp$ is to yield an 
injective Seifert fibering.

Suppose
\roster
\item[1.]
$G$ is a Lie group.
\item[2.]
$\gam$ is isomorphic to a discrete subgroup of $G$.
\item[3.]
$W$ is an admissible space (see subsection \ref{topgp}).
\item[4.]
$Q$ is a discrete group  acting properly on $W$ via $\rho: Q\ra\ttop(W)$.
\item[5.]
$1\ra\gam\ra\pi\ra Q\ra 1$ is an extension.
\end{enumerate}
\bigskip

\noindent
{\bf Definition}.
A homomorphism of $\pi$ into $\topgp$ so that the diagram of extensions
$$
\begin{array}{ccccccccc}
1 &\lra &\gam          &\lra  &\pi          &\lra  &Q   &\lra  &1\\
  &     &\bdar{\ell}   &      &\bdar{\theta}&      &\bdar{\varphi\x\rho}&  &\\
1 &\lra &\mwg\rx\inn(G)&\lra  &\topgp       &\lra  &\out(G)\x\ttop(W)&\lra &1\\
\end{array}
$$
commutes is called an \emph{injective Seifert Construction}.
The group $\topgp$ is called the \emph{universal group} for this Seifert
Construction. 
It follows that the group $\ell(G)\cdot\theta(\pi)$, generated by $\ell(G)$
and $\theta(\pi)$ in $\topgp$, acts properly on $G\x W$. The construction
then yields the injective Seifert fibering
\[
\ell(\gam)\bs\ell(G)\lra \theta(\pi)\bs(G\x W)\stackrel{p}\lra Q\bs W
\]
with the typical fiber $\gam\bs G$ and base $Q\bs W$. 

Note that the total space of our Seifert fibering is
$\pi\bs(G\x W)=Q\bs((\Gam\bs G)\x W)$.
We say it is {\it modelled on} $G\times W$. 
The base space is denoted by $B=Q\bs W$.

\subsection{The singular fibers} 
\label{sing-fiber}
\begin{prop}
Each singular fiber is the quotient of the homogeneous space $\gam\bs G$ by
a finite group of affine diffeomorphisms.
\end{prop}

\begin{proof}
There is an intermediate fibering
$\Gam\bs G \ra (\Gam\bs G)\times W \ra W$
so that the following diagram is commutative:
$$
\begin{array}{cccccl}
G &\lra &G\times W           &\lra &W        &\\
  &     &\bdal{\Gam\bs}      &     &\bdar{=} &\\
\Gam\bs G 
  &\lra &(\Gam\bs G)\times W &\lra &W        &\qquad\textrm{fibering}\\
  &     &\bdal{Q\bs}        &     &\bdar{Q\bs} &\\
\Gam\bs G
  &\lra &\pi\bs(G\times W)  &\lra &Q\bs W=B &\qquad\textrm{Seifert
fibering}\\
\end{array}
$$

The typical fiber of the Seifert fibering
$$
\Gam\bs G \lra E \stackrel{p}\lra B
$$
is the homogeneous space $\Gam\bs G$.  Then, what are the singular fibers?
\roster
\item[(1)]  Pick $w\in W$ over $b\in B$.
\item[(2)]  Over $w$, find $(\Gam\bs G)\times w$ in $(\Gam\bs G)\times W$.
\item[(3)]  Let $Q_w$ be the stabilizer of $w$.  
            Then $p^{-1}(b)=Q_w\bs(\Gam\bs G)\x w)$.
\end{enumerate}
\par\noindent
Another way:\par
\roster
\item[(4)]  Over $w$, find $G\times w$ in $G\times W$.
\item[(5)]  Let $\pi_w$ be the subgroup of $\pi$ leaving $G\times w$ invariant.
That is,
$1 \ra \Gam \ra \pi_w \ra Q_w \ra 1$
is the pullback of
$1 \ra \Gam \ra \pi \ra Q \ra 1$
via $Q_w\hra Q$.
\end{enumerate}
Then $p^{-1}(b)=Q_w\bs((\Gam\bs G)\x w)=\pi_w\bs(G\x w)$.  

The diagram (of the Seifert fibering) restricted to $G\times w$ is
$$
\begin{array}{ccccccccc}
G         &\lra  &G\times w &\lra  &w\\
          &      &\bda      &      &\bdar{=}  &     & \\
\Gam\bs G &\lra  &(\Gam\bs G)\x w 
                            &\lra  &w\\
          &      &\bda      &     &\bda      &      & \\
\Gam\bs G &\lra  &\ds{\frac{(\Gam\bs G)\x w}{Q_w}}
                           &\lra  &b\\
          &      &||  &     & \\
          &      &\pi_w\bs(G\x w)&&&
\end{array}
$$
We need to understand the action of $Q_w$ on $(\Gam\bs G)\times w$, or
equivalently, the action of $\pi_w$ on $G\times w$.

The typical fibers, which are homogeneous spaces $\gam\bs G$, occur over
$b\in B$ for which $Q_w=1$; that is, where $Q$ acts freely.
For example, if $Q$ acts effectively on $W$ and $W$ is a connected
manifold, then the set of typical fibers will be open and dense in $E$.
Each of singular fibers is a quotient of $\gam\bs G$ by a finite group of
affine diffeomorphisms of $G$. This is because when one restricts 
$G\times W$ to $G\times w$,
$\on{M}(w,G)\rx\left(\aut(G)\times\ttop(w)\right)$ becomes
$r(G)\rx\aut(G)$. 

Our group $\gam$ injects into $\ell(G)$ and $\pi_w$ goes into
$\on{M}(w,G)\rx\aut(G)$ which is the same as $r(G)\rx\aut(G)$.
But $r(G)\rx\aut(G)=\ell(G)\rx\aut(G)=\aff(G)$,
because
$(a,\alpha) = (a,\mu(a)\inv)\circ(1,\mu(a)\circ\alpha)$ in
$\mwg\rx\aut(G)$.
Therefore, we have
$$
\begin{array}{ccccccccc}
1 &\lra &\Gam     &\lra   &\pi_w       &\lra &Q_w     &\lra &1\\
  &    &\bda     &      &\bda        &    &\bda    &    &\\
1 &\lra &\ell(G) &\lra   &\aff(G) &\lra &\aut(G) &\lra &1
\end{array}
$$
Thus, $\pi_w$ acts on $G\x w\subset G\x W$ as a group of affine
diffeomorphisms of $G$.
The singular fiber $p^{-1}(b)=\pi_w\bs(G\x w)$ is an infra-homogeneous space 
(see section \ref{infra-homo}) if the finite group $Q_w=\pi_w/\gam$ 
acts on $\gam\bs G$ freely, and is an orbifold covered by $\gam\bs G$ otherwise.
\end{proof}

\subsection{Infra-homogeneous Spaces}\label{infra-homo}  
For a discrete $\pi\subset\aff(G)$ for which $\pi$ acts properly on $G$ 
and $\pi\cap G=\gam$ is of finite index in $\pi$,
we obtain an example of an injective Seifert fibering by choosing 
$W=\{w\}$, a point. We get
\[
\pi\bs(G\x\{w\}) \stackrel{p}\lra B=\{w\}.
\]
This is an injective Seifert fibering with only one fiber, which is
singular, and with typical fiber the homogeneous $\gam\bs(G\x\{w\})$.

We like to make more precise the geometry carried by the singular fibers of
our Seifert fiberings.
Endow $G$ with the linear connection defined by the left-invariant
vector fields. Since the parallel transport is the effect of the left
translations on the tangent vectors of $G$, and hence clearly independent
of paths; the connection is flat. A geodesic through the identity element
$e\in G$ are the $1$-parameter subgroups of $G$ and thus defined for any
real value of the affine parameter. All geodesics are translates of
geodesics through $e$ and thus the connection is complete.
One easily checks that the torsion tensor has vanishing covariant
derivative. According to \cite[Proposition 2.1]{kt68-1}, 
$$
\aff(G)=\ell(G)\rx\aut(G)
$$
is the group of affine diffeomorphisms
(= connection-preserving diffeomorphisms) of $G$, and
$(a,\alpha)\in G\rx\aut(G)$ acts on $G$ by $(a,\alpha)(x)=a\cdot\alpha(x)$
for all $x\in G$.
For example, if $G=\bbr^n$,
$\aff(\bbr^n)=\bbr^n\rx \gl(n,\bbr)$, the ordinary affine group of $\bbr^n$.
\bigskip

Suppose $\pi$ is a discrete subgroup of $\aff(G)$ and acts properly on $G$.
Suppose $\Gam\subset\pi\cap G$ and $Q=\pi/\Gam$ is finite.  
If $\pi$ acts freely on $G$, we call
$\gam\bs G$ a \emph{homogeneous space} and the quotient 
$\pi\bs G=Q\bs(\gam\bs G)$ an \emph{infra-homogeneous space}.
However in general, $\pi$ will not
necessarily act freely and $\pi\bs G$ is then an orbifold (i.e. a
$V$-manifold).  One should notice that the finite regular (possibly
branched) covering $\Gam\bs
G\ra\pi\bs G$ is extremely nice since the action of $Q$ came from $\aff(G)$.
\bigskip

The case where $G$ is a simply connected abelian or nilpotent Lie group is
especially important for us.  For $G=\bbr^n$, $\on{O}(n)$ is a maximal compact
subgroup of $\gl(n,\bbr)$.  A uniform (= cocompact) discrete subgroup $\pi$ of 
$\bbr^n\rx\on{O}(n)$ is called a {\it crystallographic group}.  
By a theorem of Bieberbach,
$\pi\cap\bbr^n$ is isomorphic to $\bbz^n$, and is a lattice 
(=uniform discrete subgroup) of $\bbr^n$.  
[In general, a lattice $\gam$ of a Lie group $G$ is a discrete subgroup such
that $\gam\bs G$ has finite volume. For a solvable Lie group $G$, this is
equivalent to $\gam\bs G$ being compact.]
If $\pi$ is torsion-free, we call $\pi$ a {\it Bieberbach group}
(=torsion free crystallographic group).  
Flat manifolds are the orbit spaces $\pi\bs\bbr^n$, where $\pi$ is a 
Bieberbach group.  Note that each flat manifold is finitely covered by a 
flat torus $\bbz^n\bs\bbr^n$.

Let $G$ be a connected, simply connected nilpotent Lie group.
Choose a maximal compact subgroup $C$ of $\aut(G)$.  
A uniform discrete subgroup $\pi$ of $G\rx C$ is called an {\it almost
crystallographic group}.  
If it is torsion-free, it is called an {\it almost Bieberbach group}.
An almost Bieberbach group $\pi$ yields an {\it infra-nilmanifold} $\pi\bs G$.
Note here again that any infra-nilmanifold is finitely covered by the 
nilmanifold $\Gam\bs G$, where $\Gam=\pi\cap G$.  
See \cite{lr85-1} and \cite{ausl60-1} for more details.

Thus, nilmanifolds are a generalization of tori, and infra-nilmanifolds
correspond to flat manifolds.  It is also known that a manifold is
diffeomorphic to an infra-nilmanifold if and only if it is {\it almost flat}.
This term is due to Gromov.  See \cite{fh83-1} for a proof of the above fact.
%---end of section 2 ----------------------------------------------------

\section{Group Cohomology}

\subsection{Introduction}
The injective Seifert Construction entails embedding a discrete group $\pi$
into the universal group $\topgp$ in such a way that the following diagram
of short exact sequences of groups commutes:
$$
\begin{array}{ccccccccc}
1 &\lra &\gam          &\lra  &\pi          &\lra  &Q   &\lra  &1\\
  &     &\bdar{\ell}   &      &\bdar{\theta}&      &\bdar{\varphi\x\rho}&  &\\
1 &\lra &\mwg\rx\inn(G)&\lra  &\topgp       &\lra  &\out(G)\x\ttop(W)&\lra &1\\
\end{array}
$$
Therefore, in order to understand when such constructions are possible and to
classify them when they do exist, we need to investigate the general
procedure for mapping one short exact sequence of groups into another.
These results will be formulated in terms of the first and second
cohomology of groups with not necessarily abelian coefficients. We will
need the full generality of these results later as we build in more
geometry to our Seifert Constructions by ``reducing'' our universal groups
$\topgp$.

First, we shall recall some definitions and elementary properties of group
extensions. Then we explain the correspondence between congruence classes
of extensions of abelian group $A$ by a group $Q$ in terms of the second
cohomology of $Q$ with coefficients in $A$, $H^2(Q;A)$.
If $A$ is replaced by a non-abelian group $G$, we show that the congruence
classes of extensions of $G$ by $Q$, $\opext(G,Q,\vp)$ is in one--one
correspondence with $H^2(Q;\calz(G))$, where $\calz(G)$ is the center of
$G$.

Next, we study the  problem of mapping homomorphically one short exact
sequence of groups into another. Given one such homomorphism $\theta_0$, we
show that all other possible homomorphisms are measured by the first
cohomology group.
Finally, if we change the extension of the initial exact sequence (measured
by a second cohomology group) we determine when it too maps into the target
exact sequence.

\subsection{Group Extensions}
\label{grp-ext}
A group extension is a short exact sequence
$$
1\ra G\ra E\ra Q\ra 1
$$ 
of not necessarily abelian groups. There is a naturally associated
homomorphism $\vp:Q\ra\out(G)$, called the {\it abstract kernel} of
the extension. This comes about as follows.
Pick a ``section'' $s: Q\ra E$ so that the composite $Q\ra E\ra Q$
is the identity map ($s$ is not necessarily a homomorphism).
We pick $s(1)=1$.
This defines a map $\vapt:Q\ra\aut(G)$ 
(which is a lift of the abstract kernel $\vp$)
 given by
$$
\vapt(\alpha)=\mu(s(\alpha)),
$$
where $\mu$ is the conjugation map.
Even if $\vapt$ is not a homomorphism, it induces a homomorphism 
$\vp: Q\ra\out(G)$, which is our abstract kernel. Of course, $\vp$ does
not depend on the choice of the section $s$.
We say {\it the group $E$ is an extension associated with the abstract
kernel  $\vp: Q\ra\out(G)$}.

Define $f: Q\x Q\ra G$ by
$$
s(\alpha)\cdot s(\beta) = f(\alpha,\beta)\cdot s(\alpha\beta).
$$
Then one can easily verify that
\begin{equation}
\label{co-one}
\vapt(\alpha)\circ\vapt(\beta)
= \mu(f(\alpha,\beta))\circ\vapt(\alpha\beta)
\end{equation}
\begin{equation}
\label{co-two}
f(\alpha,1) = 1 = f(1,\beta)
\end{equation}
\begin{equation}
\label{co-three}
f(\alpha,\beta)\cdot f(\alpha\beta,\gamma) 
=\vapt(\alpha)(f(\beta,\gamma))\cdot f(\alpha,\beta\gamma)
\end{equation}
for every $\alpha, \beta, \gamma\in Q$. 
In fact, (\ref{co-one}) follows from the definition of the map $f$ above,
(\ref{co-two}) follows from $s(1)=1$, and (\ref{co-three})
is a result of the associative law of $G$. 

Conversely, suppose we have maps $\vapt: Q\ra\aut(G)$
(not necessarily a homomorphism) lifting a homomorphism 
$\vp :Q\ra\out(G)$  and $f: Q\x Q\ra G$ satisfying
(\ref{co-one}), (\ref{co-two}) and (\ref{co-three}). 
Then there exists an extension
$1\ra G\ra E\ra Q\ra 1$ with the given abstract kernel $\vp$.
In fact, $E=G\x Q$ as a set, and has group operation
$$
(a,\alpha)\cdot (b,\beta)
=(a\cdot\vapt(\alpha)(b)\cdot f(\alpha,\beta),\;\alpha\beta).
$$
Thus we have

\begin{prop}
For a given abstract kernel $\vp: Q\ra\out(G)$, pick a lift
$\vapt: Q\ra\aut(G)$. Then the set of all extensions with abstract 
kernel $\vp$ is in one--one correspondence with the set of all maps
$f: Q\x Q\ra G$ satisfying 
{\rm (\ref{co-one}), (\ref{co-two}) and (\ref{co-three}).} 
\end{prop}

\begin{defn}
\label{gen-ext}
{\rm
Let $\vp: Q\ra\out(G)$ be a homomorphism, $\vapt: Q\ra\aut(G)$ be a lift, and
$f: Q\x Q\ra G$ a map satisfying 
{\rm (\ref{co-one}), (\ref{co-two}) and (\ref{co-three}).} 
Then the group $E=G\x Q$ with  group operation
$$
(a,\alpha)\cdot (b,\beta)
=(a\cdot\vapt(\alpha)(b)\cdot f(\alpha,\beta),\;\alpha\beta).
$$
is denoted by $G\x_{(f,\vapt)}Q$.
}\end{defn}

Note that some abstract kernel
$\vp: Q\ra\out(G)$ may not have any extension at all.
When there is a lift $\vapt$ which is a homomorphism (for example, when $G$ is
abelian), there exists at least one extension (by taking $f=1$, the
constant map),
the {\it semi-direct product} $G\x_{(1,\vapt)}Q$. 
As a set, this is $G\x Q$, and the group operation is 
$$
(a,\alpha)\cdot(b,\beta)=(a\cdot\vapt(\alpha)(b),\alpha\beta).
$$

\begin{defn}{\rm
Two extensions
$$
1\ra G\ra E\ra Q\ra 1\quad\on{and}\quad 1\ra G\ra E'\ra Q\ra 1
$$ 
with the same abstract kernels are {\it congruent} if there is an
isomorphism $\theta: E\ra E'$ which restricts to
the identity map on $G$, and induces the identity on $Q$.  
That is, the diagram
$$
\begin{array}{cccccccccc}
1 &\lra  &G    &\lra     &E  &\lra   &Q   &\lra  &1\\
  &      &\bdar{=} &     &\bdar{\theta}    &      &\bdar{=}&&\\
1 &\lra  &G    &\lra     &E'  &\lra   &Q   &\lra  &1\\
\end{array}
$$
is commutative.
}\end{defn}

\begin{defn}{\rm
Let $\vp: Q\ra\out(G)$ be a homomorphism, and 
$\vapt: Q\ra\aut(G)$ be a lift of $\vp$.
Then $\opext(Q,G,\vapt)$ denotes the set of all congruence classes of
extensions  of the  group $G$ by $Q$ with operator $\vapt$.

Therefore, an element $[f]\in\opext(Q,G,\vapt)$ is represented by an extension
$1\ra G\ra E\ra Q\ra 1$ with $E=G\x_{(f,\vapt)}Q$.
}\end{defn}

\subsection{Pullback and Pushout}
\begin{step}
\label{pushout-ext}
{\rm (Pushout)
Let $\gam$ be a subgroup of $G$. We say $(G,\gam)$ has \uaep (unique
automorphism extension property) if every
automorphism of $\gam$ extends to an automorphism of $G$ uniquely.
For example, if $G$ is a connected, simply connected nilpotent Lie group,
and $\gam$ is any lattice, then $(G,\gam)$ has \uaep.
In particular, $(\bbr^n,\bbz^n)$, where $\bbz^n$ spans $\bbr^n$, has \uaep.
}\end{step}

Let $\epsilon:\gam\hra G$ be an injective homomorphism.
Suppose $(G,\epsilon(\gam))$ has \uaep, and 
$$
1\lra\gam\lra\pi\lra Q\lra 1
$$
is exact. Then there is a unique group $E$ fitting the commuting diagram
$$
\begin{array}{cccccccccc}
1 &\lra  &\Gam        &\lra     &\pi         &\lra   &Q   &\lra  &1\\
  &      &\bdar{\epsilon} &         &\bda    &      &\bdar{=}&&\\
1 &\lra &G             &\lra  &E            &\lra  &Q   &\lra  &1
&\qquad{\it pushout\ \rm by\ }\epsilon:\gam\ra G \\
\end{array}
$$
Such a group $E$ can be constructed as follows.
Suppose $[f]\in\opext(Q,\gam,\vapt)$ is represented by the extension
$\pi=\gam\x_{(f,\vapt)} Q$.
That is, $\pi=\gam\x Q$ as a set, and has group operation
$$
(a,\alpha)\cdot (b,\beta)
=(a\cdot\vapt(\alpha)(b)\cdot f(\alpha,\beta),\;\alpha\beta)
$$
for $a,b\in \gam$ and $\alpha,\beta\in Q$.
Recall that $\vapt$ and $f$ satisfy the equalities (3.1), (3.2) and (3.3)
in section \ref{grp-ext}.
Since $(G,\gam)$ has the \uaep, $\vapt:Q\ra\aut(\gam)$ can be viewed as
$\vapt:Q\ra\aut(\gam)\ra\aut(G)$ and $f:Q\x Q\ra \gam$ as
$f:Q\x Q\ra \gam\stackrel{\epsilon}\ra G$.
With these new interpretations of $\vapt$ and $f$,
the equalities (3.1), (3.2) and (3.3) still hold.
Now we can define $E$ by $E=G\x_{(f,\vapt)} Q$. That is,
$E=G\x Q$ as a set, and has group operation
$$
(a,\alpha)\cdot (b,\beta)
=(a\cdot\vapt(\alpha)(b)\cdot f(\alpha,\beta),\;\alpha\beta)
$$
for $a,b\in G$ and $\alpha,\beta\in Q$.

When $\epsilon$ is an obvious inclusion $i:\gam\hra G$, we denote such 
an extension $E$ simply by $G\pi$ so that
$$
\begin{array}{cccccccccc}
1 &\lra  &\Gam        &\lra     &\pi         &\lra   &Q   &\lra  &1\\
  &      &\bdar{i} &         &\bdar{\cap}&      &\bdar{=}&&\\
1 &\lra &G             &\lra  &G\pi            &\lra  &Q   &\lra  &1
&\qquad{\it pushout\ \rm by\ }i:\gam\ra G \\
\end{array}
$$
\bigskip

\begin{step}
\label{pullback-ext}
{\rm (Pullback)
Let
$$
1\ra G\ra \tp\ra P\ra 1
$$
be exact, and $\rho:Q\ra P$ be a homomorphism. 
Then there is a unique group $\wt{Q}$ fitting the commuting diagram
}\end{step}
$$
\begin{array}{cccccccccc}
1 &\lra  &G        &\lra     &\wt{Q}         &\lra   &Q   &\lra  &1
&\qquad{\it pullback\ \rm by\ }\rho:Q\ra P \\
  &      &\bdar{=} &         &\bda    &      &\bdar{\rho}&&\\
1 &\lra &G             &\lra  &\tp            &\lra  &P   &\lra  &1\\
\end{array}
$$
Such a group $\wt{Q}$ can be constructed as follows.
Suppose $[f]\in\opext(P,G,\vapt)$ is represented by the extension
$\tp=G\x_{(f,\vapt)} P$. That is, $\tp=G\x P$ as a set, 
and has group operation
$$
(a,\alpha)\cdot (b,\beta)
=(a\cdot\vapt(\alpha)(b)\cdot f(\alpha,\beta),\;\alpha\beta),
$$
where $\vapt:P\ra\aut(G)$ and $f:P\x P\ra G$.
Now, $\vapt:P\ra\aut(G)$ gives rise to
$\vapt:Q\stackrel{\rho}\ra P\ra\aut(G)$, 
and $f:P\x P\ra G$ gives rise to
$f:Q\x Q\stackrel{\rho\x\rho}\ra P\x P\ra G$.
With these new interpretations of $\vapt$ and $f$, 
the equalities (3.1), (3.2) and (3.3) still hold. 
Now we can define $\wt{Q}$ by, $\wt{Q}=G\x_{(f,\vapt)} Q$.
That is, $\wt{Q}=G\x Q$ as a set, and has group operation
$$
(a,\alpha)\cdot (b,\beta)
=(a\cdot\vapt(\alpha)(b)\cdot f(\alpha,\beta),\;\alpha\beta)
$$
for $a,b\in G$ and $\alpha,\beta\in Q$.

\subsection{Extensions with abelian kernel $A$ and $H^2(Q;A)$}
Consider the case where $G=A$ is abelian in the previous section.
In other words, we try to classify all group extensions of an abelian group
$A$ by $Q$ with abstract kernel $\vp: Q\ra\out(A)$.
Since $A$ is abelian, $\aut(A)=\out(A)$, and there is no need of having 
a lift $\vapt$. Moreover (\ref{co-one}) holds automatically.
A map $f :Q\x Q\ra A$ satisfying the two equalities 
(\ref{co-two}) and (\ref{co-three}) is called a {\it factor set} 
(or {\it 2--cocycle}) for the abstract kernel $\vp$.
The set of all factor sets is denoted by $Z_{\vp}^2(Q;A)$ or simply, 
by $Z^2(Q;A)$.

Let $\lambda:Q\ra A$ be any map such that $\lambda(1)=0$. 
Define $\del\lambda: Q\x Q\ra A$ by
\begin{equation}
\label{delta-lambda}
\del\lambda(\alpha,\beta)
=\lambda(\alpha)+\vp(\alpha)(\lambda(\beta))-\lambda(\alpha\beta).
\end{equation}
It turns out that such a $\del\lambda$ is a 2-cocycle and is called  a
{\it 2-coboundary}. The set of all 2-coboundaries is denoted by
$B^2_\vp(Q;A)$. Clearly, $B^2_\vp(Q;A)$ is a subgroup of $Z^2_\vp(Q;A)$.
Let $f_1, f_2: Q\x Q\ra A$ be 2-cocycles. We say $f_1$ is {\it
cohomologous} to $f_2$ if $f_1-f_2 \in B^2_\vp(Q;A)$; that is,
if there is a map $\lambda: Q\ra A$ such that
$f_2-f_1=\del\lambda$.

We define the second cohomology group as the quotient group
$$
H^2_\vp(Q;A) = Z^2_\vp(Q;A)/B^2_\vp(Q;A).
$$

Since the abstract kernel $\vp:Q\ra\out(A)$ lifts to a homomorphism
$\vp:Q\ra\out(A)=\aut(A)$, there is always an extension; namely, the
semi-direct product $A\rx Q$.
Then every 2-cocycle $f\in Z^2_{\vp}(Q;A)$ gives rise to an extension of $A$ by
$Q$, which can be denoted by $A\x_{(f,\vp)} Q$. Furthermore, it is easy to see
that two cocycles $f_1$ and $f_2$ yield congruent extensions 
$A\x_{(f_1,\vp)} Q$ and $A\x_{(f_2,\vp)} Q$
if and only if $f_1$ and $f_2$ are cohomologous.
Consequently, we have 
$$
\opext(A,Q,\vp)\equiv H^2_\vp(Q;A)
$$ 
as sets, and they classify the extensions of $A$ by $Q$ with abstract 
kernel $\vp$, up to congruence. It is customary to identify the zero 
element of $H^2_\vp(Q;A)$ with the extension class of the semi-direct product
$A\rx Q$.
Then addition of the 2-cocycles in $Z^2_{\vp}(Q;A)$ induces the group
operations in $H^2_\vp(Q;A)$.

\subsection{Extensions with non-abelian kernel $G$ and $H^2(Q;\calz(G))$}
\label{ext-non-abel}
For any homomorphism $\vp: Q\ra\out(G)$, take any lift
$\vapt: Q\ra\aut(G)$. This, in turn, induces a unique homomorphism
$\vapt: Q\ra\aut(\calz(G))$ by restriction (regardless which lift of $\vp$
is chosen, and even if $\vapt$ was not a homomorphism).
We shall define an action of $H_{\vapt}^2(Q;\calz(G))$ on the set
$\opext(G,Q;\vapt)$ which will turn out to be simply transitive.

Let $g\in Z^2_{\vapt}(Q;\calz(G))$; that is, $g:Q\x Q\ra\calz(G)$
satisfying the equalities
(\ref{co-two}) and (\ref{co-three}) of section \ref{grp-ext}.
For any $[f_0]\in\opext(G,Q;\vapt)$, define
$$
(g\cdot f_0)(\alpha,\beta)=g(\alpha,\beta)\cdot f_0(\alpha,\beta).
$$
Then, clearly, $g\cdot f_0$ satisfies the equalities
(\ref{co-one}), (\ref{co-two}) and (\ref{co-three}).
Moreover, if two extensions 
$G\x_{(f_0,\vapt)} Q$ and $G\x_{(f_1,\vapt)} Q$ are congruent, so are 
$G\x_{(g\cdot f_0,\vapt)} Q$ and $G\x_{(g\cdot f_1,\vapt)} Q$.
Therefore
$$
(g,f_0)\mapsto g\cdot f_0
$$
is an action of $Z^2(Q;\calz(G))$ on $\opext(G,Q;\vapt)$.

Suppose $[g\cdot f]=[f]$ for some $[f]\in\opext(G,Q;\vapt)$. 
Then one easily sees
that $g\in B^2_{\vapt}(Q;\calz(G))$. This gives us a simple action of 
$H^2_{\vapt}(Q;\calz(G))$ on $\opext(G,Q;\vapt)$.

Now we show that this action is transitive. Suppose
$f,f'\in \opext(G,Q;\vapt)$. Let
$$
g(\alpha,\beta)=f(\alpha,\beta)\cdot f'(\alpha,\beta)\inv
$$
for all $\alpha,\beta\in Q$.
Then $g(\alpha,\beta)$ lies in $\calz(G)$, and it is easy to check
$$
g:Q\x Q\ra\calz(G)
$$
satisfies the equalities (\ref{co-two}) and (\ref{co-three}).
Thus, $g\in Z^2_{\vapt}(Q;\calz(G))$.
Consequently, we have shown the action of $H^2_{\vapt}(Q;\calz(G))$ 
on $\opext(G,Q;\vapt)$ is simply transitive. Thus,
$$
\opext(G,Q,\vapt)\approx H^2_{\vapt}(Q;\calz(G))
$$
as sets, and they classify the extensions of $G$ by $Q$ with abstract 
kernel $\vp$, up to congruence.
In general, however, there may not exist any extension of $G$ by $Q$ with 
the given abstract kernel.
Moreover, even when there exists one, there is no {\it a priori}
special element, like the semi-direct product $G\rx Q$, because the 
semi-direct product can be formed when and only when the abstract kernel
lifts to an homomorphism into $\aut(G)$.

\subsection{$H^1(Q;C)$ with a non-abelian $C$}
\label{non-abel-first}
The first cohomology is useful in classifying homomorphisms from one group
to another.
Let $C$ be a group (not necessarily abelian) and 
$\vp: Q\ra\aut(C)$ be a homomorphism.
We shall define the first cohomology set $H^1_{\vp}(Q;C)$.
A map $\eta:Q\ra C$ is called a {\it $1$--cocycle} if it satisfies
\[
\eta(\alpha\beta)=\eta(\alpha)\cdot\vp(\alpha)(\beta)
\]
for all $\alpha,\beta\in Q$.
The set of all $1$--cocycles is denoted by $Z^1_{\vp}(Q;C)$.
Two $1$-cocycles $\eta, \eta':Q\ra C$ are {\it cohomologous} if there exists
$c\in C$ such that
\[
\eta'(\alpha)=c\cdot\eta(\alpha)\cdot\vp(\alpha)(c\inv)
\]
for all $\alpha\in Q$.
The set of cohomologous classes of $Z^1_{\vp}(Q;C)$ is denoted by 
$H^1_{\vp}(Q;C)$.
When $C$ is abelian, $H^1_{\vp}(Q;C)$ is the ordinary first cohomology group.

\subsection{Mapping one group extension into another}
\label{two-cond}
\begin{defn}{\rm
A diagram of short exact sequences of groups
\begin{equation}
\label{half-comm}
\begin{array}{ccccccccc}
1 &\llra     &G        &\llra  &E  &\llra &Q &\llra  &1\\
  &         &\bdar{i}  &     &\bdar{\psi} &     &\bdar{\tthb}&\\
1 &\llra     &G'       &\llra  &E' &\llra &Q'&\llra  &1\\
\end{array}
\end{equation}
is called \emph{half-commutative} if $i$ is an inclusion and
\roster
\item[(A)] The right hand side square is commutative, and
\item[(B)] For every $x\in E$,
$i(x\cdot a\cdot x\inv)=\psi(x)\cdot i(a)\cdot\psi(x)\inv$
for all $a\in G$.
\end{enumerate}
}\end{defn}
\bigskip

Condition (B) is weaker than the left hand side square being
commutative (i.e., $\psi|_G=i$).
It simply means that
$$
i\circ\mu(x)=\mu'(\psi(x))\circ i: G\ra G'
$$
for all $x\in E$. It is not required that $\psi|_G=i$. 

Let $[f]\in\opext(Q,G,\vapt)$ represent the congruence class of the 
extension $E$. That is, $E=\gq$.
The homomorphism $\psi:E\ra E'$ gives rise to a map
$$
\vapt':Q \lra\aut(G')
$$
by $\vapt'(\alpha)=\mu'(\psi((1,\alpha))$, where $\mu'$ is conjugation in $E'$.
Since $\vapt(\alpha)=\mu(1,\alpha)$, 
the two maps $\vapt:Q\ra\aut(G)$ and $\vapt':Q\ra\aut(G')$ are related by
the commuting diagram
$$
\begin{array}{ccc}
G &\stackrel{\vapt(\alpha)}\lra  &G\\
\bdar{i} &                       &\bdar{i}\\
G' &\stackrel{\vapt'(\alpha)}\lra  &G'
\end{array}
$$
This shows that $\vapt'(\alpha)$ maps $i(G)$ onto itself, and hence
$C=C_{G'}(i(G))$, the centralizer of $i(G)$ in $G'$, onto itself also.

\begin{thm}
\label{main-alg}
Let $\psi:E\ra E'$ be a homomorphism which makes the diagram
{\rm (\ref{half-comm})} half commutative. Let $E=\gq$, and define
$\vapt'(\alpha)=\mu'(\psi(1,\alpha))$ for $\alpha\in Q$ 
{\rm ($\mu'$ is conjugation in $E'$)}.
Then,
\roster
\item
There exists a homomorphism $\theta: E\ra E'$ making the diagram
\begin{equation}
\label{const-diag}
\begin{array}{ccccccccc}
1 &\llra     &G        &\llra  &E  &\llra &Q &\llra  &1\\
  &         &\Big\da{i}  &      &\bdar{\theta}   &     &\Big\da{\tthb}&&\\
1 &\llra     &G'       &\llra  &E' &\llra &Q'&\llra  &1\\
\end{array}
\end{equation}
commutative if and only if there exists a map $\lambda: Q\ra C=C_{G'}(i(G))$
satisfying
\begin{equation}
\label{cocycle-condition}
i(f(\alpha,\beta),1)\cdot \lambda(\alpha\beta)\inv
=\lambda(\alpha)\inv\cdot\vapt'(\alpha)(\lambda(\beta)\inv)
\cdot\psi(f(\alpha,\beta),1).
\end{equation}
for all $\alpha,\beta\in Q$.
\item
Suppose that there exists a homomorphism $\theta_0: E\ra E'$ as above.
Then $\mu'\circ\theta_0$ induces a 
homomorphism  $\vp'': Q\ra\aut(C)$.

The set of all $\theta$'s fitting the diagram, up to conjugation by
elements of $G'$ (in fact, of $C$), is in one--one
correspondence with $H^1_{\vp''}(Q;C)$.
\item
Suppose that there exists a homomorphism $\theta_0: E\ra E'$ as above,
and that $i(\calz(G))\subset\calz(G')$.
Let $[E'']\in \opext(Q,G,\vapt)$ be an extension. Then there exists a
homomorphism $E''\ra E'$ completing the diagram if and only if
$$
[E'']-[E]\ \in\ \on{Ker}\left\{i_*:H^2(Q;\calz(G))\ra H^2(Q;\calz(G'))\right\}.
$$
\end{enumerate}
\end{thm}

\begin{proof}
(1) The extension $E$ is $G\x Q$ with group law 
$$
(a,\alpha)\cdot (b,\beta)
=(a\cdot\vapt(\alpha)(b)\cdot f(\alpha,\beta),\;\alpha\beta).
$$
Then 
$$
(1,\alpha)(a,1)(1,\alpha)\inv=(\vapt(\alpha)(a),1).
$$

\noindent
Suppose there is a homomorphism $\theta : E\ra E'$ making both squares
commutative. Then 
\[
\begin{array}{llll}
\theta(1,\alpha)\cdot i(a,1)\cdot\theta(1,\alpha)\inv
&=\theta((1,\alpha)(a,1)(1,\alpha)\inv)\cr
&=i((1,\alpha)(a,1)(1,\alpha)\inv)\cr
&=\psi(1,\alpha)\cdot i(a,1)\cdot\psi(1,\alpha)\inv.
\end{array}
\]
Since conjugations of $i(a,1)$ 
by $\theta(1,\alpha)$ and $\psi(1,\alpha)$ are the same, 
their difference must lie in the centralizer of $i(G)$ in $G'$.
In other words, there exists a map
$$\lambda: Q\ra C=C_{G'}(i(G))$$
so that 
$\theta(1,\alpha)=\lambda(\alpha)\inv\cdot\psi(1,\alpha)$.
Thus, in general,
\[
\begin{array}{ll}
\theta(a,\alpha)
&=\theta(a,1)\cdot\theta(1,\alpha)\cr
&=i(a,1)\cdot\lambda(\alpha)\inv\cdot\psi(1,\alpha).
\end{array}
\]
Now 
\begin{equation}
\label{prf-1}
\begin{array}{ll}
\theta((1,\alpha)(1,\beta))
&=\theta(f(\alpha,\beta),\alpha\beta))\cr
&=i(f(\alpha,\beta),1)\cdot \lambda(\alpha\beta)\inv\cdot\psi(1,\alpha\beta).\cr
\end{array}
\end{equation}
\noindent 
On the other hand,
\begin{equation}
\label{prf-2}
\begin{array}{ll}
\theta(1,\alpha)\cdot \theta(1,\beta)
&=(\lambda(\alpha)\inv\cdot\psi(1,\alpha))\cdot(\lambda(\beta)\inv\cdot\psi(1,\beta))\cr
&=\lambda(\alpha)\inv\cdot\{\psi(1,\alpha)\cdot\lambda(\beta)\inv
\cdot\psi(1,\alpha)\inv\}\cdot\{\psi(1,\alpha)\cdot\psi(1,\beta)\}\cr
&=\lambda(\alpha)\inv\cdot\vapt'(\alpha)(\lambda(\beta)\inv)
\cdot\psi(f(\alpha,\beta),\alpha\beta)\cr
&=\lambda(\alpha)\inv\cdot\vapt'(\alpha)(\lambda(\beta)\inv)
\cdot\psi(f(\alpha,\beta),1)\cdot\psi(1,\alpha\beta)\cr
\end{array}
\end{equation}
From the equalities (\ref{prf-1}) and (\ref{prf-2}), we get
\[
i(f(\alpha,\beta),1)\cdot \lambda(\alpha\beta)\inv
=\lambda(\alpha)\inv\cdot\vapt'(\alpha)(\lambda(\beta)\inv)
\cdot\psi(f(\alpha,\beta),1).
\]

Conversely, suppose there exists $\lambda: Q\ra C$ satisfying this
equality. One simply defines $\theta :E\ra E'$ by
\[
\theta(a,\alpha) =i(a,1)\cdot\lambda(\alpha)\inv\cdot\psi(1,\alpha).
\]
We claim that $\theta$ is a desired homomorphism.
This is shown by a series of calculations:
\bigskip

\(\theta(a,\alpha)\cdot\theta(b,\beta)\)
\mmm{=\{i(a,1)\cdot\lambda(\alpha)\inv\cdot\psi(1,\alpha)\}\cdot
  \{i(b,1)\cdot\lambda(\beta)\inv\cdot\psi(1,\beta)\}\hfill
  ({\rm definition\ of\ } \theta)}
\mmm{=i(a,1)\cdot\lambda(\alpha)\inv\cdot
  \{\psi(1,\alpha)\cdot i(b,1)\cdot\psi(1,\alpha)\inv\} \cdot}
\mmm{\phantom{=\, }
  \{\psi(1,\alpha) \cdot\lambda(\beta)\inv\cdot\psi(1,\alpha)\inv\}\cdot
  \{\psi(1,\alpha)\cdot\psi(1,\beta)\}}
\mmm{=i(a,1)\cdot\lambda(\alpha)\inv\cdot
  \vapt'(\alpha)(i(b,1))\cdot\vapt'(\alpha)(\lambda(\beta)\inv)
  \cdot\psi(f(\alpha,\beta),\alpha\beta)}
\rightline{(since $\psi$ is a homomorphism and
$(1,\alpha)(1,\beta)=(f(\alpha,\beta),\alpha\beta)$)}
\mmm{=i(a,1)\cdot\vapt'(\alpha)(i(b,1))\cdot\lambda(\alpha)\inv\cdot
  \vapt'(\alpha)(\lambda(\beta)\inv)
  \cdot\psi(f(\alpha,\beta),1)\cdot\psi(1,\alpha\beta)}
\rightline{($\lambda(\alpha)\inv\in C$)}
\mmm{=i(a,1)\cdot i(\vapt(\alpha)(b,1)\cdot
  \{\lambda(\alpha)\inv\cdot \vapt'(\alpha)(\lambda(\beta)\inv)\cdot
  \psi(f(\alpha,\beta),1)\}
  \cdot\psi(1,\alpha\beta)}
\rightline{($\vapt'(\alpha)\circ i=i\circ\vapt(\alpha)$)}
\mmm{=i(a,1)\cdot i(\vapt(\alpha)(b,1)\cdot
  \{i(f(\alpha,\beta),1)\cdot\lambda(\alpha\beta)\inv\}
  \cdot\psi(1,\alpha\beta)
  \hfill{\rm (equality\ (\ref{cocycle-condition}))}}
\mmm{=i(a\cdot \vapt(\alpha)(b)\cdot f(\alpha,\beta))
  \cdot\lambda(\alpha\beta)\inv \cdot\psi(1,\alpha\beta)}
\mmm{=\theta(a\cdot \vapt(\alpha)(b)\cdot f(\alpha,\beta),\alpha\beta)
\hfill({\rm definition\ of\ } \theta)}
\mmm{=\theta((a,\alpha)\cdot(b,\beta)).
\hfill({\rm group\ operation\ in\ } E)}

\bigskip

(2)
For $x\in E$, $c\in C=C_{G'}(i(G))$ and $a\in G$, it is easy to see that
$\mu'(\theta_0(x))(c)$ and $i(a)$ commute with each other.
This shows that the homomorphism
$E \stackrel{\theta_0}\lra E' \stackrel{\mu'}\lra \aut(G')$, where $\mu'$
is conjugation in $E'$, leaves $C$ invariant. Therefore we have a 
homomorphism
$$
E \stackrel{\theta_0}\lra E' \stackrel{\mu'}\lra \aut(C).
$$
Furthermore, for $a\in G$, 
$$
(\mu'\circ\theta_0)(a)=\mu'(i(a))
$$ 
is trivial on $C$. Thus $E\ra\aut(C)$ factors through $Q$. 
We have obtained a homomorphism 
$$
\vp''=\mu'\circ\theta_0 : Q\ra\aut(C).
$$
(Note that $\vp''$ is actually a homomorphism, even though
$\vapt:Q\ra\aut(G')$ is not).
Let $\theta : E\ra E'$ be a homomorphism fitting the diagram.
Then $\theta$ must be of the form
$$
\theta(\alpha)=\eta(\alpha)\theta_0(\alpha)
$$
for some map $\eta : E\ra G'$.
It is easy to verify that $\eta$ satisfies
$$
\eta(\alpha\beta)
=\eta(\alpha)\cdot\theta_0(\alpha)\eta(\beta)\theta_0(\alpha)^{-1}
=\eta(\alpha)\cdot \vp''(\alpha)(\eta(\beta)).
$$
Since $\theta|_{G}=i=\theta_0|_{G}$, we have $\eta(a)=1$ 
for all $a\in G$. For any $a\in G$ and $\alpha\in E$,
$\eta(\alpha a)=\eta(\alpha)$ from the above equality.
Moreover,
$\eta(\alpha)=\eta(\alpha\cdot(\alpha^{-1} a \alpha))
=\eta(a \alpha)
=\eta(a)\cdot\theta_0(a)\eta(\alpha)\theta_0(a)^{-1}
=1\cdot i(a)\eta(\alpha)i(a)^{-1}
=i(a)\eta(\alpha)i(a)^{-1}$.
This shows $\eta$ has values in the centralizer $C=C_{G'}(i(G))$.
Thus
$$
\eta : Q\ra C=C_{G'}(i(G)),
$$
and hence $\eta\in Z^1_{\vp''}(Q;C)$. 
See section \ref{non-abel-first}.
\bigskip

Conversely, a map  $\eta : Q\ra C$ satisfying the cocycle
condition gives rise to a homomorphism $\eta\cdot\theta_0:E\ra E'$.
Consequently,
$$
\eta\longleftrightarrow\eta\cdot\theta_0
$$
is a one--one correspondence between $Z^1_{\vp''}(Q;C)$ and the
set of all homomorphisms $E\ra E'$ inducing $i$ and $\tthb$ on $G$ and $Q$,
respectively.   

Let $\eta,\eta'\in Z^1_{\vp''}(Q;C)$, and suppose there exists 
$c\in C$ such that
$$
(\eta'\cdot\theta_0)(\alpha)=
c\cdot(\eta\cdot\theta_0)(\alpha)\cdot c\inv
$$
for all $\alpha\in E$. This is equal to
$$
\eta'(\alpha)=c\cdot\eta(\alpha)\cdot\vp''(\alpha)(c\inv);
$$
i.e.,
$$
[\eta']=[\eta]\in H^1_{\vp''}(Q;C).
$$
Consequently, $H^1_{\vp''}(Q;C)$ classifies all homomorphisms $E\ra E'$
inducing $i$ and $\tthb$ on $G$ and $Q$, respectively, up to conjugation by
elements of $C$.

\bigskip
(3) 
Let
$$
1\ra G'\ra \cale\ra Q\ra 1
$$
be the pullback (see subsection \ref{pullback-ext}) of the exact sequence 
$1\ra G'\ra E'\ra Q'\ra 1$
via $\tthb:Q\ra Q'$.
For any $1\ra G\ra E''\ra Q\ra 1$, there exists a homomorphism $E''\ra E'$
making the diagram 
$$
\begin{array}{ccccccccc}
1 &\llra     &G        &\llra  &E''  &\llra &Q &\llra  &1\\
  &         &\Big\da{i}  &      &\bda   &     &\Big\da{\tthb}&&\\
1 &\llra     &G'       &\llra  &E' &\llra &Q'&\llra  &1\\
\end{array}
$$
commutative if and only if there exists a homomorphism $E''\ra \cale$
making the diagram 
$$
\begin{array}{ccccccccc}
1 &\llra     &G        &\llra  &E''  &\llra &Q &\llra  &1\\
  &         &\Big\da{i}  &      &\bda   &     &\Big\da{\tthb}&&\\
1 &\llra     &G'       &\llra  &\cale &\llra &Q'&\llra  &1\\
\end{array}
$$
commutative.
Suppose $\theta_0:E\ra E'$ exists as in the diagram (\ref{const-diag}).
This yields $E\ra \cale$.
We identify
$$
\begin{array}{ccccccccc}
\opext(G,Q,\vapt) &\approx &H^2(Q;\calz(G))\\
\opext(G,Q,\vapt') &\approx &H^2(Q;\calz(G')).
\end{array}
$$
Since $i(\calz(G))\subset\calz(G')$, $i:\calz(G)\ra\calz(G')$ induces a
homomorphism 
$$
i_*: H^2(Q;\calz(G))\lra H^2(Q;\calz(G')).
$$
Then, clearly, $[E'']\in\opext(G,Q,\vapt)$ can map into $[E']$ (or,
equivalently into $\cale$) if and only if 
$[E'']-[E]\in H^2(Q;\calz(G))$ maps to $0\in H^2(Q;\calz(G'))$.
\end{proof}

Unless there is one definite extension $E$ and a homomorphism 
$\theta_0$, one cannot conclude that there is a 
homomorphism $E''\ra E'$ using this cohomology homomorphism.
On the other hand, suppose that there is an extension $E$ and a homomorphism
$E\ra E'$, and  that $H^2(Q;\calz(G'))=0$. 
Then every extension is in the kernel of 
$H^2(Q;\calz(G)) \ra H^2(Q;\calz(G'))$, and hence
every extension with the abstract kernel $\vp$ can be mapped into $E'$.
%---end of section 3 ----------------------------------------------------

\section{Seifert Fiber Space Construction with $\topgp$}
\subsection{Introduction}

Recall, from Lemma \ref{universal-group}, that the group of all weakly 
$G$-equivariant self-homeomorphisms of $G\x W$ is 
$$
\topgp = \mwg \rx (\autg \x \topw),
$$
the universal group for a Seifert Construction.

Let $\gam\subset G$ be a discrete subgroup,  
$\rho: Q\ra\ttop(W)$ a \pd\ action, and
$1\ra\gam\ra \pi\ra Q\ra 1$ a group extension (with abstract kernel 
$\vp : Q\ra\out(G)$).

Our goal is to find a homomorphism $\theta: \pi\ra\topgp$ which makes the 
diagram
$$
\begin{array}{ccccccccc}
1 &\lra  &\Gam        &\lra     &\pi         &\lra   &Q   &\lra  &1\\
  &      &\bdar{\ell} &         &\bdar{\theta}&      &\bdar{\vp\times\rho}
                                                          &&\\
1 &\lra  &\mwg\rx\inn(G)&\lra   &\topgp      &\lra   &\out(G)\rx\ttop(W) 
                                                          &\lra  &1\\
\end{array}
$$
commutative.

The construction of a homomorphism $\theta:\pi\ra\topgp$ will be achieved in
two steps. First, we map $\pi$ into a group $G\pi$ (see subsection 
\ref{pushout-ext}), and second $G\pi$ into $\topgp$. 

For the first step, we require, in addition, that $\gam$ be a lattice in a
simply connected completely solvable Lie group or a semi-simple Lie group 
for which Mostow's rigidity theorem holds. Then we stack the two steps 
together to get our desired homomorphisms.
A complete proof for the abelian case is given.
The other cases rely heavily on the abelian case, and we sketch and
reference the proofs for them.

As it turns out, the proofs of the existence also leads to uniqueness and 
rigidity theorems for these homomorphisms.
The meaning of existence, uniqueness and rigidity are then explored,
for it is these properties on which many of the applications of the Seifert
Construction are based.
Finally, we observe that the main theorem \ref{const-thm} of this section
remains valid in the smooth category.
We treat the second step first.

\subsection{Mapping an extension of $G$ into $\topgp$}
Interpreting the algebraic criterion of Theorem \ref{main-alg} 
with $E=G\pi$ and $E'=\topgp$, we get 

\begin{prop}
\label{g-into-topgp}
Let $\rho: Q\ra\ttop(W)$ be a \pd\ action, and
$1\ra G\ra E\ra Q\ra 1$ an extension with abstract kernel $\vp$.
Suppose $[f]\in\opext(Q,G,\vapt)$ represents the extension $E$.
That is, $E=G\x_{(f,\vapt)} Q$ (see subsection \ref{gen-ext}).
(Then $\vapt: Q\ra\aut(G)$ and $f: Q\x Q\ra G$ satisfy the equalities
{\rm (\ref{co-one}), (\ref{co-two})} and {\rm (\ref{co-three})} of section 
{\rm \ref{grp-ext}}).
\roster
\item
There exists a homomorphism $\theta: E\ra\topgp$ making
$$
\begin{array}{ccccccccc}
1 &\lra &G             &\lra  &E            &\lra  &Q   &\lra  &1\\
  &     &\bdal{\ell}   &      &\bdar{\theta}& &\bdar{\vp\times\rho}&      &\\
1 &\lra &\mwg\rx\inn(G)&\lra  &\topgp       &\lra  &\out(G)\x\ttop(W)&\lra &1\\
\end{array}
$$
commutative if and only if there exists a map $\lambda: Q\ra\mwg$
satisfying
\begin{equation}
\label{cond-1}
f(\alpha,\beta)=
\{\vapt(\alpha)\circ\lambda(\beta)\circ\rho(\alpha)\inv\}
\cdot\lambda(\alpha)\cdot\lambda(\alpha\beta)\inv
\in r(G)\subset\mwg
\end{equation}
for all $\alpha,\beta\in Q$. {\rm [Compare this with the equality 
(\ref{delta-lambda}) defining $\delta\lambda$].}
\item
Suppose there exists a homomorphism $\theta_0: E\ra\topgp$ as above.
Then $\mu'\circ\theta_0$ {\rm ($\mu'$ is conjugation in $\topgp$)} induces a 
homomorphism  $\vp'': Q\ra\aut(\mwg)$.

The set of all $\theta$'s fitting the diagram, up to conjugation by
elements of $\mwg$, is in one--one
correspondence with $H^1_{\vp''}(Q;\mwg)$.
\item
Suppose that there exists a homomorphism $\theta_0: E\ra\topgp$ as above.
Then $\theta_0(\calz(G))\subset\mwzg$.
Let $[E'']\in \opext(Q,G,\vapt)$ be an extension. Then there exist a
homomorphism $E''\ra\topgp$ completing the diagram if and only if
$[E'']-[E]$ lies in the kernel of 
$\ell_*: H^2(Q;\calz(G))\ra H^2(Q;\mwzg)$.
\end{enumerate}
\end{prop}

\begin{proof}
\newcommand\wa{\widetilde{\alpha}}
(1) We apply Theorem \ref{main-alg} (1) to the diagram
$$
\begin{array}{ccccccccc}
1 &\lra &G             &\lra  &E            &\lra  &Q   &\lra  &1\\
  &     &\bdal{\ell}   &      &             &      &\bdar{\vp\times\rho}& &\\
1 &\lra &\mwg\rx\inn(G)&\lra  &\topgp       &\lra  &\out(G)\x\ttop(W)&\lra &1\\
\end{array}
$$
Let $\psi:E\ra\topgp$ be given by
\[
\psi(\wa)=(1,\vapt(\wa),\rho(\wa))\in\mwg\rx(\aut(G)\x\ttop(W)),
\]
where $\vapt(\wa)\in\aut(G)$ is conjugation in $E$, and 
$\rho: E\ra Q\ra\ttop(W)$ (abuse of notation) is the given 
\pd\ action of $Q$ on $W$. 
In particular,
$$
\psi(1,\alpha)=(1,\vapt(\alpha),\rho(\alpha)).
$$
Conditions (A) in section \ref{two-cond} is
easily checked. For (B), Let $a\in G$ and $\wa\in E$.
Then
\[
\begin{array}{lll}
\ell(\wa\cdot a\cdot\wa\inv)
&=\ell(\vapt(\wa)(a))\\
&=(\vapt(\wa)(a)\inv,\mu(\vapt(\wa)(a)),1).
\end{array}
\]
On the other hand,
\[
\begin{array}{lll}
\psi(\wa)\cdot\ell(a)\cdot\psi(\wa)\inv
&=(1,\vapt(\wa),\rho(\wa))\cdot(a\inv,\mu(a),1)\cdot
(1,\vapt(\wa)\inv,\rho(\wa)\inv)\\
&=(\vapt(\wa)(a)\inv,\mu(\vapt(\wa)(a)),1).\\
\end{array}
\]
Thus, 
\[
\ell(\wa\cdot a\cdot\wa\inv)=\psi(\wa)\cdot\ell(a)\cdot\psi(\wa)\inv
\]
as we wanted.

Now we examine the equality (\ref{cocycle-condition}).
\[
\begin{array}{lll}
\ell(f(\alpha,\beta))\cdot\lambda(\alpha\beta)\inv
&=(f(\alpha,\beta)\inv,\mu(f(\alpha,\beta)),1)\cdot
(\lambda(\alpha\beta)\inv,1,1)\\
&=(\lambda(\alpha\beta)\inv\cdot f(\alpha,\beta)\inv,\mu(f(\alpha,\beta)),1),
\end{array}
\]
while
\[
\begin{array}{lll}
\lambda(\alpha)\inv
&\cdot\vapt'(\alpha)(\lambda(\beta)\inv)\cdot\psi(f(\alpha,\beta))\\
&=(\lambda(\alpha)\inv,1,1)\cdot
\{\psi(1,\alpha)\cdot\lambda(\beta)\inv\cdot\psi(1,\alpha)\inv\}
\cdot(1,\mu(f(\alpha,\beta)),1)\\
&=(\lambda(\alpha)\inv,1,1)\cdot
\{(1,\vapt(\alpha),\rho(\alpha))\cdot(\lambda(\beta)\inv,1,1)\cdot
(1,\vapt(\alpha),\rho(\alpha))\inv\}\cdot\\
&\phantom{=\ } (1,\mu(f(\alpha,\beta)),1)\\
&=(\lambda(\alpha)\inv,1,1)\cdot
(\vapt(\alpha)\circ\lambda(\beta)\inv\circ\rho(\alpha)\inv,1,1)
\cdot(1,\mu(f(\alpha,\beta)),1)\\
&=(\lambda(\alpha)\inv\cdot\{\vapt(\alpha)\circ\lambda(\beta)\inv\circ\rho(\alpha)\inv\},
\mu(f(\alpha,\beta)),1)\\
\end{array}
\]
Therefore,
\[
\ell(f(\alpha,\beta))\cdot\lambda(\alpha\beta)\inv
=\lambda(\alpha)\inv\cdot\vapt(\alpha)
(\lambda(\beta)\inv)\cdot\psi(f(\alpha,\beta))
\]
if and only if
$$
\lambda(\alpha\beta)\inv\cdot f(\alpha,\beta)\inv
=\lambda(\alpha)\inv\cdot
\{\vapt(\alpha)\circ\lambda(\beta)\inv\circ\rho(\alpha)\inv\}
$$
or equivalently,
$$
f(\alpha,\beta)=
\{\vapt(\alpha)\circ\lambda(\beta)\circ\rho(\alpha)\inv\}
\cdot\lambda(\alpha)\cdot\lambda(\alpha\beta)\inv
$$
as elements of $G$. 
This is the identity (\ref{cond-1}).

(2) Suppose there exists a homomorphism $\theta_0: E\ra\topgp$ as above.
Since $\ell(G)$ is normal in $\topgp$, and
the centralizer of $\ell(G)$ in $\mwg\rx\inn(G)$ is $\mwg$,
$\mwg$ is normal in $\topgp$.
Consequently, conjugation by elements of $\topgp$ leaves $\mwg$ invariant,
and $\mu'\circ\theta_0$ induces a homomorphism  $\vp'': Q\ra\aut(\mwg)$.
The rest is clear.

(3) Observe that the center of $\mwg\rx\inn(G)$ is $\mwzg$ so that
$\ell(\calz(G))\subset\mwzg$. Now apply Theorem \ref{main-alg} (3).
\end{proof}
\bigskip

\beginrem
When the condition (\ref{cond-1}) is satisfied, vanishing of $H^1(Q;\mwg)$ 
and $H^2(Q;\mwzg)$ (e.g., if $\calz(G)$ is contractible)
yields strong consequences by the statements (2) and (3).
Namely,
{\sl
\roster
\item[{\rm(2)}]
The homomorphism $\theta$ fitting the diagram is unique up to conjugation by
elements of $\mwg$.
\item[{\rm(3)}]
Every extension with the abstract kernel $\vp$ can be mapped into $\topgp$.
\end{enumerate}
}
\end{rem}
\bigskip

\noindent
The vanishing of $H^2(Q;\mwzg)$ only says that there is at most one
extension of  $\mwg\rx\inn(G)$  by $Q$ for a given abstract kernel.
This alone is not necessarily enough to guarantee that 
the condition (\ref{cond-1}) is automatically satisfied
because the coefficient groups $G$ and $\mwg\rx\inn(G)$ have distinct
``operator groups''.
On the other hand, we have the following important observation.
\bigskip

\begin{thm}
\label{new-remark}
If $\vapt_C:Q\stackrel{\vapt}\lra\aut(\mwg\rx\inn(G))\ra\aut(\mwg)$
is a homomorphism, and $H^2(Q;\mwzg)=0$, then $\theta$ exists as in
Proposition \ref{g-into-topgp}, and consequently condition (\ref{cond-1})
is satisfied.
\end{thm}

\noindent{\bf Proof.}
Let
$$
1\lra \mwg\rx\inn(G)\lra \cale\lra Q\lra 1
$$
be the pullback of 
$1\lra \mwg\rx\inn(G)\lra \mwg\rx(\aut(G)\x\ttop(W))\lra 
\out(G)\x\ttop(W) \lra 1$.
Then 
$$
\cale=(\mwg\rx\inn(G))\x_{(h,\vapt)} Q
$$ 
for some maps
$\vapt:Q\ra\aut(\mwg\rx\inn(G))$ and $h:Q\x Q\ra \mwg\rx\inn(G)$
satisfying
\[
\vapt(\alpha)\circ\vapt(\beta)
= \mu(h(\alpha,\beta))\circ\vapt(\alpha\beta)
\]
\[
h(\alpha,1) = 1 = h(1,\beta)
\]
\[
h(\alpha,\beta)\cdot h(\alpha\beta,\gamma) 
=\vapt(\alpha)(h(\beta,\gamma))\cdot h(\alpha,\beta\gamma)
\]
for every $\alpha, \beta, \gamma\in Q$. 
Note that the pair $(h,\vapt)$ is not unique. (In fact, each
$\vapt(\alpha)$ can change by conjugation by an element of
$\mwg\rx\inn(G)$,
and $h$ changes accordingly). 
The map $\vapt$ induces two maps $\vapt_G:Q\ra\aut(G)$ and
$\vapt_{C}:Q\ra\aut(\mwg)$.
The second map in the definition of $\vapt_C$ is obtained as follows.
Since $\mwg\rx\inn(G)=\ell(G)\x_{\calz(G)}\mwg$, and $\ell(G)$ is normal
in $\topgp$, an automorphism of $\mwg\rx\inn(G)$ leaving $\ell(G)$ invariant
induces an automorphism of $\mwg$. 
Also $\vapt_G$ induces a homomorphism $\vp:Q\ra\out(G)$.

Now let 
$$
1\ra G\ra E\ra Q\ra 1
$$
be any extension with abstract kernel $\vp:Q\ra\out(G)$.
Then there exists a map $f:Q\x Q\ra G$ such that $E=G\x_{(f,\vapt_G)} Q$
satisfies
\[
\vapt_G(\alpha)\circ\vapt_G(\beta)
= \mu(f(\alpha,\beta))\circ\vapt_G(\alpha\beta)
\]
\[
f(\alpha,1) = 1 = f(1,\beta)
\]
\[
f(\alpha,\beta)\cdot f(\alpha\beta,\gamma) 
=\vapt_G(\alpha)(f(\beta,\gamma))\cdot f(\alpha,\beta\gamma)
\]
for every $\alpha, \beta, \gamma\in Q$. 
Let us examine the first equality. This equality holds in the group $\aut(G)$.
Now $f(\alpha,\beta)\in\ell(G)$ and $\mwg=C$ centralizes $\ell(G)$.
Hence $\mu(f(\alpha,\beta))$ on $C$ is trivial.
If we replace $\vapt_G$ in the first equation by $\vapt_C$, then the
equation is satisfied on $\mwg$ in $\mwg\rx\inn(G)$ because 
$\mu(f(\alpha,\beta))$ is trivial and $\vapt_C$ is a homomorphism.
Since $\ell(G)$ and $\mwg$ generate $\mwg\rx\inn(G)$,
the first equality still holds in $\mwg\rx\inn(G)$
if we replace $\vapt_G$ with $\vapt$.
Because
$\vapt|_{\ell(G)}(\alpha)(\ell(g))=\ell(\vapt_G(\alpha)(g))$,
the remaining equalities hold inside $\mwg\rx\inn(G)$ if we interpret 
$f:Q\x Q\ra G$ as $f:Q\x Q\ra \ell(G)\subset \mwg\rx\inn(G)$.
Consequently, $E=G\x_{(f,\vapt_G)}Q$ is sitting in
$\cale'=(\mwg\rx\inn(G))\x_{(f,\vapt)}Q$.

Now, if $H^2(Q;M(W,\calz(G))=0$, the extensions $\cale$ and $\cale'$ are
congruent, and we have a desired homomorphism $E\ra \cale$.
This theorem will be used to prove Theorem \ref{const-thm}.

\subsection{Construction with $\topgp$}
\begin{step}
\label{special}
{\rm
We shall say a discrete group $\gam$ is \emph{special} if $\gam$ is
isomorphic to a lattice in any one of the following Lie groups $G$:
}\end{step}
\roster
\item[(S1)]
$\bbr^k$ for some $k>0$,
\item[(S2)]
a simply connected nilpotent Lie group,
\item[(S3)]
a simply connected completely solvable Lie group; that is,
for each $X\in\mathcal G$, the Lie algebra of $G$,
$\on{ad}(X):\mathcal G\ra \mathcal G$ has only real eigenvalues.
\item[(S4)]
a semi-simple centerless Lie group without any normal compact factors and 
if $G$ contains any 3-dimensional factors (i.e., $\psl$), then the
projection of the lattice to each of these factors is dense.
\end{enumerate}

We shall also call the Lie group $G$ \emph{special}. Such groups possess
\uliep (Unique Lattice Isomorphism Extension Property). That is,
any isomorphism between such lattices 
extends uniquely to an isomorphism of $G$.
Even more generally, if $\gam_i$ is a lattice in a special $G_i$ ($i=1,2$),
then any isomorphism $\vp:\gam_1\ra\gam_2$  extends uniquely to an
isomorphism $\Psi: G_1\ra G_2$.

Notice that:\ type (S1)$\Longrightarrow$  type (S2)$\Longrightarrow$
type (S3).
\bigskip

The following is the main construction. Its proof is deferred until section 
\ref{proof-exist-uniq-rigid}.

\begin{thm}[\rm \cite{cr69-1},\cite{klr83-1},\cite{rw77-1}, \cite{llr96-1}]  
\label{const-thm}
Let $\gam\subset G$ be a \emph{special} lattice.
Let $\rho:Q\ra\ttop(W)$ be a \pd\ action.
Then for any extension
$1\ra\Gam\ra\pi\ra Q\ra 1$ (with abstract kernel $\vp : Q\ra\out(G)$)
the following are true:
\roster
\item 
\emph{Existence:}  
There exists a proper action $\theta\col\pi\ra\topgp$ making the diagram 
\begin{equation}
\label{seif-const-diag}
\begin{array}{ccccccccc}
1 &\lra  &\Gam        &\lra     &\pi         &\lra   &Q   &\lra  &1\\
  &      &\bdar{\ell} &         &\bdar{\theta} &       &\bdar{\vp\times\rho}
                                                          &&\\
1 &\lra  &\mwg\rx\inn(G)&\lra   &\topgp      &\lra   &\out(G)\rx\ttop(W) 
                                                          &\lra  &1\\
\end{array}
\end{equation}
commutative.
For completely solvable $G$ we have in addition:
\item  
\emph{Uniqueness:}  
Congruent extensions are conjugate in
$\mwg\subset\topgp$.  More precisely, suppose
$\theta_1,\theta_2\col\pi\ra\topgp$ are two homomorphisms which fit in
diagram {\rm (\ref{seif-const-diag})}
 with fixed $\ell$ and $\vp\times\rho$, then there exists
$\lam\in\mwg$ such that $\theta_2=\mu(\lam)\circ\theta_1$.
\item
\emph{Rigidity:}  
Suppose $\theta_1,\theta_2\col\pi\ra\topgp$) are two homomorphisms which fit
in the diagram {\rm (\ref{seif-const-diag})} 
{\rm (possibly with distinct $\ell$ and $\rho$)}.  
Then there exists
$(\lam,a,h)\in\topgp$ such that $\theta_2=\mu(\lam,a,h)\circ\theta_1$,
provided that $\rho_2=\mu(h)\circ\rho_1$.
\end{enumerate}
\end{thm}

Uniqueness and rigidity for the semisimple case (S4)  will be
discussed in sections \ref{4.5.6} and \ref{4.5.7}

\subsection{The Meaning of Existence, Uniqueness and Rigidity}  
\label{meaning}
\begin{step}{\rm (Existence)
Suppose there exists $\theta\col\pi\ra\topgp$ making the diagram 
(\ref{seif-const-diag}) commutative.
Then $\theta(\pi)$ acts properly on $G\x W$. This yields the injective 
Seifert fibering
}\end{step}
\[
\ell(\gam)\bs\ell(G)\lra \theta(\pi)\bs(G\x W)\stackrel{p}\lra Q\bs W
\]
with the typical fiber $\gam\bs G$ and base $Q\bs W$. 
We have explained the properties of such a fibering, and its singular
fibers in subsections \ref{inj-seif-const}, \ref{sing-fiber} and 
\ref{infra-homo}.

\begin{step}{\rm (Uniqueness)
\label{uniqueness}
Let $\theta_0,\theta_1$ be
homomorphisms of $\pi$ into $\topgp$ such that
}\end{step}
$$
\begin{array}{ccccccccc}
1&\lra &\Gam          &\lra &\pi           &\lra &Q &\lra &1\\
 &     &\bdal{\ell}   &     &\bdar{\theta_i} &     &\bdar{\phi\times\rho}&&\\
1&\lra &\mwg\rx\inn(G)&\lra &\topgp        &\lra &\out(G)\times\ttop(W)&\lra &1
\end{array}
$$
is commutative.  If there exists $\lam\in \mwg$ such that
$\theta_1=\mu(\lam)\circ\theta_0$, the map $\lam\col W\ra G$ 
induces a homeomorphism of $G\x W$ by
$$
\lambda(x,w)=(x\cdot(\lambda(w))\inv,w),
$$
see section \ref{topgp}. This, in turn, yields a homeomorphism
$[\lam]\col {{G\x W}\over{\theta_0(\pi)}}\lra{{G\x W}\over{\theta_1(\pi)}}$.
As the commuting diagram shows
$$
\begin{array}{ccccccccc}
M_0&=\frac{G\x W}{\theta_0(\pi)} &\stackrel{[\lam]}\lra 
&\frac{G\times W}{\theta_1(\pi)}&=M_1\\
&\bda &&\bda&\\
&\rho(Q)\bs W &\stackrel{=}\lra   &\rho(Q)\bs W&
\end{array}
$$
the map $[\lambda]$ induces the identity map on the base space.
In fact, $\lambda$ is $G$-equivariant (that is, it commutes with the left
principal $G$-action).
Such spaces $M_0$ and $M_1$ are said to be \emph{strictly equivalent}.

When $G$ is special of type (S1), (S2) or (S3); or $W$ is contractible, of
type (S4) (see subsection \ref{special}),
$\lambda$ is homotopic to the constant map $e:W\ra G$ ($e$ is the identity
element of $G$).
Then the path $\{\lambda_t: 0\leq t\leq 1\}$ gives rise to a continuous 
family of homomorphisms $\theta_t: \pi\ra\topgp$ by
$$
\theta_t(\wt{\alpha})=\lambda_t\cdot \theta_0(\wt{\alpha})\cdot \lambda_t\inv
$$
for $\wt{\alpha}\in\pi$, and consequently a continuous family of homeomorphisms 
$$
[\theta_t]\col {{G\x W}\over{\theta_0(\pi)}}\lra{{G\x W}\over{\theta_t(\pi)}}.
$$
Thus ${{G\x W}\over{\theta_0(\pi)}}$ can be deformed to 
${{G\x W}\over{\theta_1(\pi)}}$ by moving just along the fibers.
In fact, the family $\lambda_t: G\x W\ra G\x W$ is $G$-equivariant.
If $\theta_0(\pi)$ commutes with $\ell(G)$, then the deformation
$[\theta_t]$ is $G$-equivariant.

\begin{step}{\rm (Rigidity)
Let $\theta_0,\theta_1$ be
homomorphisms of $\pi$ into $\topgp$ such that
}\end{step}
$$
\begin{array}{ccccccccc}
1&\lra &\Gam          &\lra &\pi           &\lra &Q &\lra &1\\
 &     &\bdal{\ell_i} &     &\bdar{\theta_i} &     &\bdar{\phi_i\times\rho_i}&&\\
1&\lra &\mwg\rx\inn(G)&\lra &\topgp        &\lra &\out(G)\times\ttop(W)&\lra &1
\end{array}
$$
is commutative.  If there exists $(\lam,a,h)\in \topgp$ such that
$\theta_1=\mu(\lam,a,h)\circ\theta_0$, the map $(\lam,a,h)\col G\times W\ra G\times
W$ induces a homeomorphism $[\lam,a,h]$
$$
\begin{array}{ccccccccc}
M_0&=\frac{G\x W}{\theta_0(\pi)} &\stackrel{[\lam,a,h]}\lra 
&\frac{G\times W}{\theta_1(\pi)}&=M_1\\
&\bda &&\bda&\\
&\rho_0(Q)\bs W &\stackrel{\overline h}\lra   &\rho_1(Q)\bs W&
\end{array}
$$
which preserves the fibers.  More precisely, $[\lam,a,h]$ is defined by
$$[\lam,a,h]\left([x,w]\right)=\left[(\lam,a,h)(x,w)\right]\,.$$

Since $(\lam,a,h)\circ\theta_0(\alp)=\theta_1(\alp)\circ(\lam,a,h)$ for all $\alp\in\pi$,
the map $[\lam,a,h]$ is well-defined.  Further,
$(\lam,a,h)(x,w)=\left(a(x)\cdot\lam h(w),h(w)\right)$ shows that $[\lam,a,h]$ is
a map from $M_0$ to $M_1$.  That is, $[\lam,a,h]$ is the descent of the
weakly equivariant fiber preserving map $(\lam,a,h)$ sending $G$ fibers to 
$G$ fibers.
Such spaces $M_0$ and $M_1$ are said to be \emph{equivalent}.

The conjugation of $\theta_0$ by $(\lambda,a,h)$ is called a \emph{Seifert
automorphism} of $G\x W$. The induced homeomorphism 
$[\lambda,a,h]: M_0\ra M_1$ is called a \emph{Seifert equivalence} or 
\emph{Seifert automorphism}.

If $\Gam$ is characteristic in $\pi$, any automorphism of $\pi$ induces an
automorphism of $\Gam$ and $Q$.  Consequently, if $M_0$ and $M_1$ are Seifert
fiber spaces modelled on $G\times W$ which have the same fundamental group (or
``orbifold fundamental group'') they are equivalent as Seifert fiber spaces,
provided that the base spaces are rigidly related.  (i.e., there exists
$h\in\ttop(W)$ for which $\rho_1=\mu(h)\circ\rho_0)$.

Similarly, if $\tau:\pi\ra\pi'$ is an isomorphism carrying $\gam$
isomorphically onto $\gam'$, and inducing an isomorphism $\ol{\tau}:Q\ra
Q'$, then the analogous rigidity statement for 
$M=\theta(\pi)\bs (G\x W)$ and $M'=\theta'(\pi')\bs (G\x W)$ holds, cf 
\cite[8.5]{cr69-1} and \cite[2.4]{klr83-1}.
More precisely, if there exists $h\in\ttop(W)$ such that
$\mu(h)\circ\rho=\rho'\circ\ol{\tau}$, then there exists
$\wh{h}=(\lambda,a,h)\in\topgp$ such that 
$\theta'\circ\tau=\mu(\wh{h})\circ\theta$.
\medskip

\begin{exam}
[When $W$ is a point]
\label{w-point}
{\rm
Suppose $W=\{p\}$, a point.  Then 
$$
\begin{array}{lll}
\topgp&=\on{M}(p,G)\rx(\aut(G)\x\ttop(p))\\
      &=r(G)\rx\aut(G)\\
      &=\ell(G)\rx\aut(G).\\
\end{array}
$$
Also, note that $\on{M}(p,G)\rx\inn(G)=r(G)\rx\inn(G)=\ell(G)\rx\inn(G)$.
Since $Q$ acts on $W$ properly, $Q$ must be a finite group.

We assume that $\gam\subset G$ is a \emph{special} lattice, and apply
Theorem \ref{const-thm}.
Then diagram (\ref{seif-const-diag}) becomes
$$
\begin{array}{ccccccccc}
1 &\lra  &\Gam        &\lra     &\pi         &\lra   &Q   &\lra  &1\\
  &      &\bdar{\ell} &         &\bdar{\theta} & &\bdar{\vp} &&\\
1 &\lra &\ell(G)\rx\inn(G)&\lra&\ell(G)\rx\aut(G) &\lra   &\out(G) &\lra  &1\\
\end{array}
$$
Furthermore, since $\ell(\gam)\subset\ell(G)$, the above diagram induces
$$
\begin{array}{ccccccccc}
1 &\lra  &\Gam        &\lra     &\pi         &\lra   &Q   &\lra  &1\\
  &      &\bdar{\ell} &         &\bdar{\theta} & &\bdar{\wt{\vp}} &&\\
1 &\lra &\ell(G)&\lra&\ell(G)\rx\aut(G) &\lra   &\aut(G) &\lra  &1\\
\end{array}
$$
where $\wt{\vp}$ is a homomorphism lifting $\vp$.
[In particular, this shows that, for special $G$, the abstract kernel 
$\vp:Q\ra\out(G)$ (with $Q$ finite) lifts to a homomorphism 
$\wt{\vp}:Q\ra\aut(G)$].

Consequently, $\theta(\pi)\subset\ell(G)\rx C$ for a finite subgroup
$C\subset\aut(G)$.
This means that there exists a Riemannian metric on $G$ for which 
$\theta(\pi)\subset\isom(G)$.

Assume that $\theta$ is injective. (Recall that the kernel of $\theta$ is
finite. Also, note that if $\pi$ is torsion free, $\theta$ is injective.)
For $\bbr^n = G$,  $\vapt(Q)$ lies in a compact subgroup of $\gl(n,\bbr)$
and can be conjugated within $\gl(n,\bbr)$ to $O(n)$ and consequently,
by rigidity, $\theta(\pi)$ is conjugated by an affine diffeomorphism  into  
the group of Euclidean motions $E(n) = \bbr^n \rx O(n)$.
An injective Seifert fiber space modelled on 
$\bbr^n\times\textrm{point}$ is a flat manifold if $\pi$ is torsion free,
and a flat orbifold otherwise.  
Similarly, if $G$ is a connected, simply connected nilpotent Lie group, 
an injective Seifert fiber 
space modelled on $G\times\textrm{point}$ is an infra-nilmanifold
if $\pi$ is torsion free, and an infranil-orbifold otherwise.  
(Moreover, in both cases, each such manifold (or orbifold) must arise in 
this fashion.)  Since $W$ is a point, the rigidity on $W$ always holds, 
and by the rigidity theorem, any two infra-$G$-manifolds modelled on special
$G$ with isomorphic fundamental groups are ``affinely'' diffeomorphic.
(Recall that $\topgp$ reduces to $G\rx\aut(G)$ in this case.)
This generalizes the classical result of Bieberbach's from crystallographic
groups to almost crystallographic groups in the nilpotent case.
}
\end{exam}
\medskip
\hfuzz=1.4pt

\subsection{Proofs of Theorem \ref{const-thm}}
\label{proof-exist-uniq-rigid}
There are two steps. First we map $\pi$ into an extension $E$ of $G$ by $Q$, 
and second we map $E$ into $\topgp$.

\begin{step}{\rm
Every special Lie group $G$ has the \uliep.
For every $G$ with \uliep, 
let $1\ra \gam\ra \pi\ra Q\ra 1$ be an extension which represents 
$f\in\opext(\gam,Q,\vapt)$ satisfying {\rm (\ref{co-one}), (\ref{co-two})} 
and {\rm (\ref{co-three})} of section {\rm \ref{grp-ext}}.
One can form a pushout by the inclusion $i:\gam\hra G$ to obtain $E=G\pi$.
See subsection \ref{pushout-ext}.
Then $f\in\opext(G,Q,\vapt)$ and the diagram
$$
\begin{array}{ccccccccc}
1 &\lra &\Gam     &\lra   &\pi     &\lra  &Q     &\lra &1\\
  &     &\bda     &       &\bda    &      &\bda  &     &\\
1 &\lra &G        &\lra   &E=G\pi       &\lra  &Q     &\lra &1\\
\end{array}
$$
is commutative.
}\end{step}

\begin{step}{\rm
Now in order to map $E$ into $\topgp$, we apply Proposition \ref{g-into-topgp}
for the cases (S1), (S2), and (S3). 
These cases depend upon the abelian case whose proof we give in complete
detail.
When $G=\bbr^k$ is abelian, the existence of $\lambda$ in the equality 
(\ref{cond-1})
\[
f(\alpha,\beta)=(\vapt(\alpha)\circ\lambda(\beta)\circ\rho(\alpha)\inv)
\cdot\lambda(\alpha)\cdot\lambda(\alpha\beta)\inv
\]
means that $[f]$ will have to be $0$ in $H^2(Q;\on{M}(W,\bbr^k))$.
The following Lemma takes care of this problem:
}\end{step}

\begin{lem}
[{\cite[Lemma 8.4]{cr69-1} and \cite[Theorem 7.1]{cr71-3}}]
Let $\rho: Q\ra W$ be a properly discontinuous action with $Q\bs W$ 
paracompact.
For any homomorphism $\vp: Q\ra\gl(k,\bbr)$, and an action of $Q$ on 
$\mw{\bbr^k}$ by
$\alpha\cdot\lambda = \vp(\alpha)\circ\lambda\circ\rho(\alpha)\inv$,
$H^i(Q;\mw{\bbr^k}) = 0$ for all $i>0$.
\end{lem}

\begin{proof}
We shall prove the lemma, as was done in \cite{cr69-1}, for $Q\bs W$
compact. For the general case we refer the reader to \cite{cr71-3}.
Let $U_x$ be a neighborhood of $x\in W$ with the property that
$Q_xU_x=U_x$ and such that if $U_x\cap\alpha U_x\neq\emptyset$, then
$\alpha\in Q_x$. Let $V_x=Q U_x$.
Take $M(V_x)\subset\mw{\bbr^k}$ to be the submodule of all maps with
support in $V_x$.
Similarly, $M(U_x)\subset M(V_x)\subset\mw{\bbr^k}$ is the subspace of all
maps with support in $U_x$, and $M(U_x)$ is a $Q_x$--module.
Then
$$\hhom_{\bbz Q_x}(\bbz Q;M(U_x))\approx M(V_x)
$$
because $\bbz Q$ is a free $\bbz Q_x$ module with a basis given by choosing
coset representatives. Using Shapiro's lemma, it follows that
$$
H^*(Q_x;M(U_x))\approx H^*(Q;M(V_x)).
$$
(See, e.g., \cite[Chapter X, Proposition 7.4]{cart-eil56-1} or
\cite[Chapter III 5.8 and 6.2]{brow82-1}).
Now $H^i(Q_x;M(U_x))=0$ for all $i\geq 1$, because $Q_x$ is a finite group
and $M(U_x)$ is an $\bbr$-vector space.

If $y\not\in Q(x)$, then neighborhoods  $U_x$ and $U'_y$ can be chosen so 
that  $V_x\cap U'_y = \emptyset$ because $Q$ acts properly
discontinously on $W$. 
This insures that $Q\bs W$ is Hausdorff (and completely
regular). 
Because $Q\bs W$ is compact, there exists a finite set of points
$x_1, x_2,\cdots, x_n$ with neighborhoods
$U_{x_1}, U_{x_2},\cdots, U_{x_n}$ as above whose images,
$V^*_{x_1}, V^*_{x_2},\cdots, V^*_{x_n}$ in $Q\bs W$ cover $Q\bs W$.
Take a partition of unity subordinate to this open covering.
By composing each partition function with the quotient map
$W\stackrel{p}\lra Q\bs W$, we get a partition of unity
$\epsilon_1,\epsilon_2\cdots,\epsilon_n$ on $W$ subordinate to 
$V_1,V_2,\cdots, V_n$, where $V_i=p\inv(V^*_{x_i})=Q_{x_i}(U_{x_i})$.
Each partition function satisfies $\epsilon_j(\alpha w)=\epsilon_j(w)$,
for all $\alpha\in Q$, $w\in W$.
Now multiplication by $\epsilon_j$ is a $Q$-module homomorphism of
$\mw{\bbr^k}$ into $M(V_j)$.
Consequently, there is induced homomorphisms ${\epsilon_j}_*$ and the 
commutative diagram
\begin{equation}
\begin{array}{ccccccccc}
H^*(Q;\mw{\bbr^k}   &\stackrel{{\epsilon_j}_*}\lra &H^*(Q;\mw{\bbr^k}\\
\Big\da{\epsilon_j}_* &                          &i_* \Big\ua\\
H^*(Q;M(V_j))       &\stackrel{=}\lra          &H^*(Q;M(V_j))
\end{array}
\end{equation}
with ${\epsilon_1}_*+{\epsilon_2}_*+\cdots+{\epsilon_n}_*=\textrm{identity}$.
This yields the lemma.
\end{proof}

\begin{step}{\rm
For a nilpotent Lie group $G$, one shows that the equation 
\[
f(\alpha,\beta)=(\vapt(\alpha)\circ\lambda(\beta)\circ\rho(\alpha)\inv)
\cdot\lambda(\alpha)\cdot\lambda(\alpha\beta)\inv
\]
has a solution for $\lambda$. 
One reduces this  problem to iterated applications of 
$H^2(Q; \mw{\bbr^k} = 0$ to guarantee the existence
while $H^1(Q; \mw{\bbr^k} = 0$ for the uniqueness. 
The readers are referred to \cite{klr83-1} for the details in 
the nilpotent case.  
}\end{step}

\begin{step}{\rm
Theorem \ref{new-remark} is used to prove the existence in the completely
solvable and type (S4) cases as well as providing a different procedure to
that used in \cite{klr83-1} for the nilpotent case. Details can be found in
\cite[Theorem 3]{lr89-1}.
The uniqueness for the completely solvable case is not given in 
\cite{lr89-1} but it can be obtained inductively from the nilpotent case.

Rigidity in the abelian case was first proved by a cohomological argument,
\cite[\S 8.5]{cr69-1}. In the nilpotent case, rigidity follows
from uniqueness and \uliep, \cite[\S 2.4]{klr83-1}.
In fact, the same argument works for completely solvable $G$.
}\end{step}

\begin{step}
\label{4.5.6}
{\rm
Suppose $\aut(G)$ splits as $\inn(G)\rx\out(G)$ , and $\calz(G) = 0$ 
(for example, $G$ of type (S4)).
If we replace $\mwg$ by $r(G)$ and $\topgp$ by
$r(G)\rx(\aut(G) \x \ttop(W))\subset\topgp$,
then a homomorphism $\theta: E\ra r(G)\rx(\aut(G) \x \ttop(W))$ exists 
(with $\theta|_G=\ell$ on $G$) by Theorem \ref{new-remark}.
% The reason is that the equality (\ref{cond-1}) is now satisfied by virtue 
% of the splitting and 
% $\psi(\wt{\alpha}) = (1, \mu_E (\wt{\alpha}), \rho(\alpha))$, as an element
% of $r(G)\rx(\aut(G) \x \ttop(W))$ satisfies  conditions A and B of
% subsection \ref{two-cond} (hence the hypothesis of Theorem \ref{main-alg}).
% So the condition (\ref{cond-1}) is satisfied as in Theorem \ref{new-remark}.
% Since $\calz(r(G)) =0$, $\theta$ then exists.
Thus $E$ maps into 
$$
r(G)\rx(\aut(G)\x\ttop(W))=(\ell(G)\rx\aut(G))\x\ttop(W).
$$
This induces two homomorphisms $E\ra \ell(G)\rx\aut(G)$ and $E\ra\ttop(W)$.
Consequently, the action of $E$ on $G\x W$ is diagonal.

Note that uniqueness does not hold in $(\ell(G)\rx\aut(G))\x\ttop(W)$ for
$\gam\x Q$. The different strict equivalence classes are in 1--1
correspondence with $H^1(Q;r(G))=$ conjugacy classes of representations of
$Q$ into $r(G)$.
}\end{step}

\begin{step}
\label{4.5.7}
{\rm
In \cite{rw77-1}, yet another  construction is given in the (S4) case.
Instead of $G\times W$ being the
uniformizing space, $G/K\times W$ is chosen, where $K$ is a maximal compact
subgroup of $G$.  Consequently, $G/K$ is a simply connected symmetric space.
The typical fiber becomes the locally symmetric space $F=\gam\bs G/K$ and 
singular fiber is the quotient of $F$ by a {\it finite group of isometries.}  
The space $E$ has a finite covering $E'$ (perhaps branched) where $E'$ is 
$F\times(Q'\bs W)$ and $Q'$ is the kernel of $Q\ra\out(\Gam)$.  
Since $\out(\gam)$ is finite, $Q'$ has finite index in $Q$ and the finite 
group $Q/Q'$ acts diagonally on $F\times(Q'\bs W)$.

In \cite{lr96-1}, this is explained in more detail by showing that the
uniformizing group $\ttop_{(G,K)}(G/K\x W)$ is isomorphic to
$\isom(G/K)\x\ttop(W)$. In this context, existence, uniqueness and rigidity
are all shown to be valid.
Typical fibers which are locally symmetric spaces have a great deal of
geometric interest and the uniqueness and rigidity of this modified
construction has significant applications.
When modeling on $G/K\x W$ instead of $G\x W$ we shall call this the
modified injective Seifert Construction of type (S4).
}\end{step}
%\bigskip\hrule\bigskip

\begin{step}{\rm
Now we stack the two constructions together to obtain a homomorphism
$\gam\ra\topgp$:
$$
\begin{array}{ccccccccc}
1 &\lra &\Gam     &\lra   &\pi     &\lra  &Q     &\lra &1\\
  &     &\bda     &       &\bda    &      &\bda  &     &\\
1 &\lra &G             &\lra  &E            &\lra  &Q   &\lra  &1\\
  &     &\bdal{\ell}   &      &\bdar{\theta}&      &\bda&      &\\
1 &\lra &\mwg\rx\inn(G)&\lra  &\topgp       &\lra  &\out(G)\x\ttop(W)&\lra &1\\
\end{array}
$$
This completes the proof of Theorem \ref{const-thm}. \hfill Q.E.D.
}\end{step}

\subsection{Smooth Case}
If $W$ is a smooth manifold and $Q$ acts on it smoothly so
that $\rho\col Q\ra\diff(W)$, then the construction can be done smoothly.  The
universal group is now $\calc(W,G)\rx\left(\aut(G)\times\diff(W)\right)$
which is certainly the subgroup of weakly $G$-equivariant diffeomorphisms 
of $G\times W$.  
The same proof as for $\mwg$ works in proving $H^i(Q;\calc(W,\calz(G)))=0$,
$i=1,2$ in the case of $G$ completely solvable.  
We use this fact by simply referring to the ``smooth Seifert
fiber space construction''.

If $W$ is a complex manifold and $Q$ acts holomorphically, then one can also
ask whether the construction can be done holomorphically on $\bbc^k\times W$.
There are two types of obstructions.  It could happen that the necessary
$\bbc^k$-bundle over $W$, while trivial as a smooth bundle, may not be 
trivial as a holomorphic bundle.  A more serious matter is that
$H^2\left(Q;\calh(W,\bbc^k)\right)$, where $\calh(W,\bbc^k)$
denotes holomorphic maps, does not necessarily vanish
and so not every smooth realization has a holomorphic one.  
Even when the first difficulty does not arise,
%If $W$ is a Stein manifold, the first difficulty does not arise but 
the latter one may still persist.  
The solution to these problems and the corresponding theory is
carefully worked out in \cite{cr71-3}.  One particular feature of the extended
theory is that the uniformizing space need no longer be a product $G\times W$
but may be any principal $G$-bundle over $W$.
See \cite{lr89-1} for Seifert fibered spaces modelled on principal 
bundles and \ref{holomorphic} for an illustration of the holomorphic theory.
\medskip

Another general approach to Seifert fiber spaces is due to A.~Holmann 
\cite{holm64-1}.  His procedure extends the classical fiber bundle methods.
%---end of section 4 ----------------------------------------------------

\section{Applications}
\subsection{Introduction}
The Seifert Construction, which is a special embedding, $\theta:\pi\ra\topgp$,
of the group $\pi$ into $\topgp$ such
that $\pi$ acts properly on $G\x W$, preserves some of  the properties
of both $G$ and $W$ on $\theta(\pi)\bs(G\x W)$.
Furthermore, the action of $\pi$ on $G\x W$ ``twists'' the topology and
geometry of $G$ and $W$ to create the orbit space  $\theta(\pi)\bs (G\x W)$  
in the same way that the group structures of $\Gamma$ and $Q$ ``twist''
to create the group $\pi$. In other words, this algebraic twisting of
$\pi$ makes the geometric twisting of the ``bundle with singularities''
$$
\Gamma\bs G\ra \theta(\pi)\bs (G\x W)\ra Q\bs W,
$$
where the homogeneous space $\Gamma\bs G$ is a typical fiber. 
In the several applications, we have selected to include here this features
seems especially prominent.

One of the important geometric problems that has motivated the development
of Seifert fiberings is the construction of closed aspherical manifolds
realizing Poincar\'e duality groups $\pi$ of the form
$1\ra\gam\ra\pi\ra Q\ra 1$.
The Seifert Construction enables one to find explicit aspherical manifolds
$M(\pi)$ when $Q$ acts on a contractible manifold $W$ and $\gam$ is a
torsion free lattice in a Lie group.

For a topological manifold, the homotopy classes of self-homotopy
equivalences can be regarded as algebraic data.
We shall show how the Seifert Construction can be used to lift finite
subgroups of homotopy classes to an action on the manifold.

The final illustration is a description of the classical 3-dimensional
Seifert manifolds and how they fit into our scheme of general Seifert
fiberings. This description will also be useful for the final section.

\subsection{Existence of Closed Smooth $K(\pi,1)$-manifolds}
There are two difficult problems related to the title. They are:

\noindent$\bullet$
\emph{Which groups can be the fundamental group of a closed 
aspherical manifold?} and\par
\noindent$\bullet$
\emph{If $\pi$ is the fundamental group of an aspherical manifold,
can we give an actual explicit construction of an aspherical manifold
for the group $\pi$?}

There are some general criteria for the first problem such as $\pi$ must be
finitely presented, have finite cohomological dimension and satisfy
Poincar\'e duality in that dimension.
The Seifert Construction gives answers to both questions for large classes of
groups $\pi$.
The idea is that if $1\ra\gam\ra\pi\ra Q\ra 1$ is a torsion free extension
where $\gam$ is the fundamental group of a closed aspherical manifold and
$Q$ is a group acting properly on a contractible manifold $W$ with compact
quotient, then $\pi$ should be the fundamental group of a closed aspherical
manifold. We have the following

\begin{thm}
\label{6-1-1}
Let $\gam$ be a cocompact special lattice in $G$ (see subsection \ref{special})
and $\rho:Q\ra\ttop(W)$ be
a \pd\ action on a contractible manifold $W$ with compact quotient.
If $1\ra\gam\ra\pi\ra Q\ra 1$ is a torsion free extension of $\gam$ by $Q$,
then for any Seifert Construction $\theta:\pi\ra\topgp$,
\roster
\item
$M(\theta(\pi))=\theta(\pi)\bs(G\x W)$ is a closed aspherical manifold if 
$\gam$ is of type {\rm (S1), (S2)} or {\rm (S3)},
\item
$M(\theta(\pi))=\theta(\pi)\bs((G/K)\x W)$,
where $K$ is a maximal compact subgroup of $G$,
is a closed aspherical manifold
if $\gam$ is of type {\rm (S4)}.
\end{enumerate}
\end{thm}

\begin{proof}
We know from Theorem \ref{const-thm} that for each extension
$1\ra\gam\ra\pi\ra Q\ra 1$, there exists a homomorphism of $\pi$ into
$\topgp$ (resp., $\ttop_{(G,K)}(G/K\x W)$ in the case of type (S4), see
\cite{lr96-1} for this notation).
We need only to check that this homomorphism is injective.
Suppose $Q_0$ is the kernel of $\phi\x\rho: Q\ra\out(G)\x\ttop(W)$.
Then $Q_0$ is finite since the $Q$ action on $W$ is \pd.
Let $1\ra\gam\ra E\ra Q_0\ra 1$ be the pullback via $Q_0\subset Q$.
The group $E$ is torsion free since $\pi$ is assumed to be torsion free.
The restriction $\theta|_E$ defines an action of $E$ on $G\x W$ (resp.,
$G/K\x W$).
Since $\gam$ is of finite index in $E$, no non-trivial element of $E$ can
fix $G\x W$ (resp., $G/K\x W$).
For if it did, then some power would be a non-trivial element of $\gam$
which does not fix $G$ (resp., $G/K$).
Therefore, $\theta$ is injective;
$\theta(\pi)$ acts properly and freely since it is torsion free.
\end{proof}

\begin{rem}{\rm
1. If $W$ is a smooth contractible manifold and $\rho: Q\ra\diff(W)$,
then the construction can be done smoothly and $M(\theta(\pi))$ is smooth.\par
2. If $\rho_1$ and $\rho_2$ are rigidly related
(i.e., there exists
$h\in\ttop(W)$ for which $\rho_2=\mu(h)\circ\rho_1)$)
and $\gam$ is characteristic in
$\pi$, then $M(\theta_1(\pi))$ and $M(\theta_2(\pi))$ are homeomorphic via
a Seifert automorphism.
Moreover, if we fix $\ell$ and $\rho$, then the constructed $M(\theta(\pi))$
are all strictly equivalent.\par
3. When $W=\{p\}$ is a point (a $0$-dimensional contractible manifold),
then $Q$ must be finite for $Q$ to act \pd ly, 
and every $\rho:Q\ra\ttop(\{p\})$ is rigidly related.
The closed aspherical manifolds constructed are infra-$G$-manifolds.
cf. Example \ref{w-point}.\par
4. One important application of these constructions is that they provide
model  aspherical manifolds with often strong geometric properties.
If one wants to study the famous conjecture that two closed aspherical
manifolds with isomorphic fundamental groups are homeomorphic via the
methods of controlled surgery, then the constructed aspherical Seifert
manifolds are excellent model manifolds.
}\end{rem}

The above procedure can be extended for even more general extensions. As an
example,

\begin{thm}
\label{6-1-2}
Let $\pi$ be a torsion-free extension of a virtually poly-$\bbz$ group $\gam$
by $Q$, where $Q$ acts on a contractible manifold $W$ 
properly discontinuously with compact quotient.  
Then there exists a closed $K(\pi,1)$-manifold.
\end{thm}

\begin{proof}
A torsion-free virtually poly-$\bbz$ group $\gam$ has a unique maximal normal
nilpotent subgroup $\Delta$, which is called the \emph{discrete nilradical}
of $\gam$. See \cite{ausl63-1}.
Then the quotient $\gam/\Delta$ is virtually free abelian of finite rank.
Furthermore, since $\Delta$ is a characteristic subgroup of $\gam$,
it is normal in $\pi$. Consider the commuting diagram with exact rows and
columns:
$$
\begin{array}{ccccccccc}
  &     &1             &      &1            &      &      &      & \\
  &     &\bda          &      &\bda         &      &      &  &\\
  &     &\Delta        &\stackrel{=}\lra  &\Delta & &      &      & \\
  &     &\bda          &      &\bda         &      &      &      & \\
1 &\lra &\gam          &\lra  &\pi          &\lra  &Q     &\lra  &1\\
  &     &\bda          &      &\bda         &      &\bdar{=}&  &\\
1 &\lra &\gam/\Delta   &\lra  &\pi/\Delta   &\lra  &Q     &\lra  &1\\
  &     &\bda          &      &\bda         &      &      &  &\\
  &     &1             &      &1            &      &      &      & \\
\end{array}
$$
Since $\gam/\Delta$ is virtually free abelian of finite rank (say of $s$),
it contains a characteristic subgroup $\bbz^s$.
Let $Q'=(\pi/\Delta)/\bbz^s$. Then the natural projection $Q'\lra Q$
has a finite kernel. Therefore, if we let $Q'$ act on $W$ via $Q$, the
action will still be properly discontinuous.

One can do a Seifert fiber space construction with the exact sequence
$$
1\lra\bbz^s\lra\pi/\Delta\lra Q'\lra 1
$$
which yields a \pd\ action of $\pi/\Delta$ on $\bbr^s\x W$ with
compact quotient.
Using this action of $\pi/\Delta$ on $\bbr^s\x W$, one does a Seifert 
fiber space construction with the exact sequence
$$
1\lra\Delta\lra\pi\lra \pi/\Delta\lra 1.
$$
This gives rise to a \pd\ action of $\pi$ on $N\x(\bbr^s\x W)$,
where $N$ is the unique simply connected nilpotent Lie group containing 
$\Delta$ as a lattice,  with compact quotient.

If the space $W$ is smooth, and the action of $Q$ on $W$ is smooth, both 
constructions can be done smoothly so that the \pd\ action of $\pi$ on 
$N\x(\bbr^s\x W)$ is smooth.

In any case, since the group $\pi$ is torsion free, the resulting action
of $\pi$ on $N\x(\bbr^s\x W)$ is free. Consequently, we get a closed
$K(\pi,1)$-manifold 
$$M=\pi\bs (N\x\bbr^s\x W).$$
It has a Seifert fiber structure
$$
F\lra M \lra Q\bs W
$$
where the typical fiber $F$ itself has a Seifert fiber structure
$$
\Delta\bs N\lra F\lra T^s=\bbz^s\bs\bbr^s.
$$
In fact, since the action of the characteristic subgroup $\bbz^s$ on
$\bbr^s$ is free, $F$ is a genuine fiber bundle, with fiber a nilmanifold
$\Delta\bs N$ over the base torus $T^s$.
\end{proof}

The space $W$ does  not have to be aspherical. As far as the action of $Q$
is \pd, the construction works. The resulting action of $\pi$ is free if
and only 
if the pre-image of $Q_w$ (the stabilizer of the $Q$ action at each $w\in W$) in
$\pi$ is torsion free. In this case, the space $\pi\bs(G\x W)$ will not 
be aspherical. See \ref{const-thm}.

\subsection{When is $\theta$ injective?}
Let
$$
\begin{array}{ccccccccc}
1 &\lra  &\Gam        &\lra     &\pi         &\lra   &Q   &\lra  &1\\
  &      &\bdar{\ell} &         &\bdar{\theta} &       &\bdar{\vp\times\rho}
                                                          &&\\
1 &\lra  &\mwg\rx\inn(G)&\lra   &\topgp      &\lra   &\out(G)\rx\ttop(W) 
                                                          &\lra  &1\\
\end{array}
$$
be a Seifert construction.
The argument as mentioned just above also shows that the action of $\pi$ on
a general $G\x W$ will be free (and hence $\theta$ is injective) 
if the pre-image of $Q_w$ (the stabilizer of the $Q$ action at each $w\in W$) in
$\pi$ is torsion free.
Of course $\theta$ may be injective, that is, $\pi$ acts effectively on
$G\x W$, without necessarily acting freely.

We determine when the homomorphism $\theta:\pi\ra\topgp$ has
trivial kernel for special lattices $\gam$ in $G$ of types (S1), (S2) or
(S3) in subsection \ref{special}.

\begin{prop}
\label{inj-prime}
$\theta$ is injective and $\theta(\pi)\cap\ell(G)=\theta(\gam)$
if and only if $\vp\x\rho$ is injective.
\end{prop}

\begin{proof}
If $\vp\x\rho$ is injective, then clearly $\theta$ is injective and 
$\theta(\pi)\cap\ell(G)=\theta(\gam)$. Now assume that $\theta$ is
injective.
Let $K$ be the kernel of $\vp\x\rho$, and $1\ra\gam\ra E\ra K\ra 1$ 
the pullback via $K\subset Q$. 
Since $K$ acts trivially on $W$, it is finite and $E$ must map injectively
into $\mwg\rx\inn(G)=\ell(G)\x_{\calz(G)}\mwg$.
The group $K$ maps injectively into $\mwg/\calz(G)$ by the natural
projection $\ell(G)\x_{\calz(G)}\mwg\ra\mwg/\calz(G)$, because we
hypothesized that $\theta(\pi)\cap\ell(G)=\theta(\gam)$.
We complete the argument by showing $\mwg/\calz(G)$ has no non-trivial
elements of finite order.
This follows from the fact that the exponential map ${\mathcal G}\ra G$ 
is bijective for our special $G$.
Therefore the equation $x^n=a$ in $G$ has a unique solution for any 
$n\in\bbz$ and every $a\in G$. 
So if $f\in\mwg$ represents an element of $K$ in $\mwg/\calz(G)$, then
$f^n\in\calz(G)$, where $n$ is the order of $f$.
Therefore, $f^n(w)=a$, for some $a\in\calz(G)$, and all $w\in W$.
But there exists a unique $x\in G$ such that $f(w)=x$. Now, $\calz(G)$ is
$\bbr^k$ for some $k$, and each $a\in\calz(G)$ has a unique $n$th root in
$\calz(G)$, $x\in\calz(G)$.
Since $f$ has order $n$, $n$ must be 1 and $\gam=E$.
\end{proof}

\subsection{Rigidity of infra-homogeneous spaces}
Recall from \ref{infra-homo} the following definitions.
Let $G$ be a connected Lie group. A quotient of $G$ by a lattice $\gam$
is called a \emph{homogeneous space}.
More generally, let $\aff(G)=G\rx\aut(G)$ act on $G$ by
$$
(a,\alpha)\cdot x=a\cdot \alpha(x).
$$
An \emph{infra-homogeneous space} is a quotient of $G$ by a group
$\pi\subset\aff(G)$, acting \pd ly and freely, 
such that $\gam=\pi\cap G$ is a uniform lattice of $G$
and $\gam$ has a finite index in $\pi$.
Therefore, an infra-homogeneous space is finitely covered by a homogeneous
space. 
For example, a flat Riemannian manifold may be called an ``infra-torus''.
To emphasize $G$, sometimes we use the term 
\emph{infra-$G$-manifold} for an infra-homogeneous space.

Now we explain a little more detail for the claims made in Example 
\ref{w-point}.
Let $\Delta_i$ be a special lattice of type (S1), (S2) or (S3).
That is,
$$
1\lra\Delta_i\lra\pi_i\lra Q_i\lra 1
$$
is an extension, $\theta_i:\pi_i\ra\aff(G_i)$, $Q_i$ finite,
$\rho_i:Q_i\ra\out(G_i)$ is injective, and hence by Proposition
\ref{inj-prime}, 
$\Delta_i=\theta_i(\pi_i)\cap\ell(G_i)$, for $i=1,2$.
Suppose $h:\pi_1\ra\pi_2$ is an isomorphism.
Then, we claim, $\Delta_1$ is mapped isomorphically onto $\Delta_2$ and hence 
$Q_1$ onto $Q_2$.
This isomorphism extends to an isomorphism of $G_1$ onto $G_2$.
Thus we may think of
$$
\theta_i:\pi\ra\aff(G)=\ttop_G(G\x\{p\})
$$
such that $\theta_i(\pi)\cap G=\Delta_i$, $i=1,2$.
Now because $\rho_i:Q_i\ra\{p\}$ is trivial, the embeddings are conjugate
by an element of $\aff(G)$.

It remains to verify that $h$ maps $\Delta_1$ isomorphically onto $\Delta_2$. 
It is well known that $\Delta_i$ is characteristic in the abelian and
nilpotent cases.
For a proof of $\Delta_i$ being characteristic  in the completely solvable
case, we refer the reader to \cite[Proposition 4.5]{lee95-1}.
Therefore $h$ maps $\Delta_1$ isomorphically onto  $\Delta_2$.
\bigskip

Instead of appealing to the rigidity theorems for the Seifert Construction
as we did above, we shall instead obtain these same extensions of the
classical Bieberbach theorems directly.
In the next theorem, the calculations involved resemble the proof  of 
uniqueness of the Seifert Construction when $W$ is a point, and show how
closely linked the Seifert Construction is to the Bieberbach theorems.

\begin{thm}
\label{gen-rigid}
Suppose $G$ has {\rm ULIEP} (see subsection \ref{special}), 
and $H^1(\Psi;G)$ {\rm(non-abelian cohomology)} 
vanishes for every finite subgroup $\Psi\subset\aut(G)$.
Let $\pi, \pi'\in \aff(G)$ be finite extensions of lattices in $G$.
Then every isomorphism $\theta:\pi\ra\pi'$ is a conjugation by an element
of $\aff(G)$.
\end{thm}

\begin{proof}
Let $\gam=G\cap\pi$, $\gam'=G\cap\pi'$ be the pure translations.
Let $\Lambda=\gam\cap\theta^{-1}(\gam')$.
Then $\theta|_{\Lambda}: \Lambda\ra \theta(\Lambda)$ is an isomorphism of
lattices of $G$.
By ULIEP, $\theta|_{\Lambda}$ extends uniquely to an automorphism
$C: G\ra G$.
Thus, $\theta|_{\Lambda}=C|_{\Lambda}$, and hence,
$\theta(z,I)=(Cz,I)$ for all $(z,I)\in\Lambda$.
Let us denote the composite homomorphism
$$
\pi\stackrel{\theta}\lra\pi'\hra \aff(G)\ra\aut(G)
$$
by $\tb$, where $\aff(G) = G\rx \autg\ra\autg$ is the projection.
Let $\Psi=\pi/\Lambda$.
Then clearly, $\tb$ factors through $\Psi$.
Define a map $\lambda: \pi\ra G$ by
$$
\theta(w,K)=(Cw\cdot \lambda(w,K),\tb(w,K))
\eqno (1)
$$
It is easy to see that $\lambda(zw,K)=\lambda(w,K)$ for all $z\in\Lambda$.
Therefore, $\lambda$ is really a map
$$
\lambda: \Psi\ra G.
$$
For any $(z,I)\in\Lambda$ and $(w,K)\in\pi$, apply $\theta$ to both sides of
$(w,K)(z,I)(w,K)^{-1}=(w\cdot Kz\cdot w^{-1},I)$ to get
$$
Cw\cdot \lambda(w,K)\cdot\tb(w,K)(Cz)\cdot \lambda(w,K)^{-1}\cdot(Cw)^{-1}
=\theta(w\cdot Kz\cdot w^{-1}).
$$
However, $w\cdot Kz\cdot w^{-1}\in \Lambda$ since $\Lambda$ is normal in
$\pi$, and the right hand side equals to $C(w\cdot Kz\cdot w^{-1})
=Cw\cdot CKz\cdot (Cw)^{-1}$ since $C: G\ra G$ is an automorphism.
From this, we have
$$
\tb(w,K) (Cz)= \lambda(w,K)^{-1}\cdot CKz\cdot \lambda(w,K)
\eqno (2)
$$
This is true for all $z\in\gam$. Note that $\tb(w,K)$, $C$ and $K$ are
automorphisms of the Lie group $G$.
By ULIEP of $G$, {\sl the equality ($2$) holds true for all $z\in G$.}

We claim that, with the $\Psi$-structure on $G$ via $\tb: \Psi\ra \autg$,
$\lambda\in Z^1(\Psi;G)$; i.e., $\lambda:\Psi\ra G$ is a crossed homomorphism.
We shall show
$$
\lambda((w,K)\cdot(w',K'))=\lambda(w,K)\cdot\tb(w,K)(\lambda(w',K'))
$$
for all $(w,K), (w',K')\in \pi$.
(Note that we are using the elements of $\pi$ to denote the elements
of $\Psi$).
Apply $\theta$ to both sides of
$(w,K)(w',K')=(w\cdot Kw',KK')$ to get
$Cw\cdot \lambda(w,K)\cdot\tb(w,K)[Cw'\cdot \lambda(w',K')]
=C(w\cdot Kw')\cdot \lambda((w,K)(w',K'))$.
From this, it follows that
$$
\lambda((w,K)(w',K'))
=(CKw')^{-1}\cdot \lambda(w,K)\cdot\tb(w,K)(Cw')\cdot\tb(w,K)\lambda(w',K').
$$
By ($2$), we have
$\tb(w,K)(Cw')= \lambda(w,K)^{-1}\cdot CKw'\cdot \lambda(w,K)$ so that
$\lambda((w,K)\cdot(w',K'))=\lambda(w,K)\cdot\tb(w,K) \lambda(w',K')$.
This shows that $\lambda$ is a crossed homomorphism.

By the assumption, $H^1(\Psi;G)$ vanishes.
This means that any crossed homomorphism is ``principal''. 
In other words, there exists $d\in G$ such that
$$
\lambda(w,K)=d\cdot\tb(w,K)(d^{-1})
\eqno (3)
$$
Let $D=\mu(d^{-1})\circ C$.
We check that $\theta$ is the conjugation
by $(d,D)=(d,\, \mu(d^{-1})\circ C)\in \aff(G)$.
Using (1), (2) and (3), one can show
$\tb(w,K)\circ\mu(d^{-1})\circ C= \mu(d^{-1})\circ C\circ K$.
Thus, for any $(w,K)\in\pi$,
$$
\begin{array}{lll}
\theta(w,K)\cdot (d,D)
        &=(Cw\cdot \lambda(w,K),\;\tb(w,K))\cdot(d,\;\mu(d^{-1})\circ C)\cr
        &=(Cw\cdot \lambda(w,K)\cdot\tb(w,K)(d),\;
          \tb(w,K)\circ\mu(d^{-1})\circ C)\cr
        &=(Cw\cdot d\cdot\tb(w,K)(d^{-1})\cdot\tb(w,K)(d),\;
                        \tb(w,K)\circ\mu(d^{-1})\circ C)\cr
        &=(Cw\cdot d,\;\mu(d^{-1})\circ C\circ K)\cr
        &=(d,D)\cdot(w,K).\cr
\end{array}
$$
This finishes the proof of our claim.
\end{proof}
\bigskip

\begin{cor}
Suppose $G$ has {\rm ULIEP}, and $H^1(\Psi;G)$ {\rm(non-abelian cohomology)} 
vanishes for every finite subgroup $\Psi\subset\aut(G)$.
Then homotopy equivalent infra-G-manifolds are affinely diffeomorphic.
\end{cor}

\begin{proof}
Let $M=\pi\bs G$, $M'=\pi'\bs G'$ be infra-$G$-manifolds.
A homotopy equivalence $M\ra M'$  induces an isomorphism $\theta: \pi\ra\pi'$.
Since $\theta$ maps some subgroup of $G\cap\pi$ into a lattice of $G'$,
there is an isomorphism of $G$ onto $G'$.
Using this isomorphism, we identify $G$ with $G'$.
Now apply Theorem \ref{gen-rigid} to  $\theta: \pi\ra\pi'$ to find an element
$h=(d,D)\in \aff(G)$ which conjugates $\pi$ onto $\pi'$.
This gives a weakly equivariant map
$$
(h,\mu(h)): (G,\pi)\ra (G,\pi')
$$
and $h$ gives rise to an affine diffeomorphism $M\ra M'$,
which is homotopic to the original map.
\end{proof}
\bigskip

\begin{cor} 
\begin{enumerate} 
\item {\rm(Bieberbach)} Homotopy equivalent flat manifolds are affinely 
diffeomorphic.
\item {\rm\cite{lr85-1}} Homotopy equivalent infra-nilmanifolds are affinely 
diffeomorphic.
\item Homotopy equivalent infra-solvmanifolds of type {\rm(R)} 
are affinely diffeomorphic.
\end{enumerate}
\end{cor}

Notice, in particular, that the isomorphism $\theta$ maps the lattice
$\gam=G\cap\pi$ onto $\gam'=G\cap\pi'$.
This is also true topologically on the orbifold level.  

\begin{proof}
What is needed for this conclusion is the element $h=(d,D)\in\aff(G)$ which
conjugates $\pi$ onto $\pi'$. We have this element immediately from the
rigidity of the Seifert Construction for lattices of type (S1), (S2) and
(S3). If we wish to avoid the Seifert Construction, then we need to verify
that $H^1(\Psi;G)$ vanishes for each finite subgroup $\Psi\subset\aut(G)$.
This is accomplished in \cite{lr85-1} and \cite{lee95-1}, and well known
for $G$ of type (S1).
\end{proof}

\subsection{Lifting Problem for Homotopy Classes of Self-Homotopy 
Equivalences}  
Let $M$ be an admissible space (see section \ref{topgp}) and $\cale(M)$ be
the $H$-space of homotopy equivalences of $M$ into itself.  
Any $f\in \cale(M)$ induces an isomorphism 
$f_*\col\pi_1(M,x)\ra\pi_1\left(M,f(x)\right)$.  
By choosing a path $\ome$ from $x$ to $f(x)$, we have an automorphism 
$f^\ome_*$ of $\pi_1(M,x)$, defined by
$f^\ome_*\left([\tau]\right)=[\ome^{-1}\cdot(f\circ\tau)\cdot\ome]$.  
A different choice of $\ome$ alters $f^\ome_*$ only by an inner automorphism.  
Therefore, we obtain a homomorphism 
$$
\gamma: \cale(M)\ra\out(\pi),
$$
$\pi=\pi_1(M,x)$.
If $M$ is a $K(\pi,1)$ space, then $\cale_0(M)$ is the kernel of $\gamma$ so
that $\gamma$ factors through $\pi_0\left(\cale(M)\right)
=\cale(M)/\cale_0(M)$, where $\cale_0(M)$ is the self homotopy equivalences
homotopic to the identity.
Moreover, $\gamma$ is onto since every automorphism of $\pi$ can be realized
by a self homotopy equivalence of $M$.

A homomorphism $\vp\col F\ra\out(\pi)\cong\pi_0\left(\cale(M)\right)$ is 
called an {\it abstract kernel}.  An injective abstract kernel is the same 
as a subgroup of homotopy classes of self-homotopy equivalences of $M$.  
A {\it lifting} of $\vp$ as a group of homeomorphisms is a homomorphism 
$\hat\vp\col F\ra\ttop(M)$ which makes
$$
\begin{array}{ccc}
F        & \stackrel{=}\lra &F\\
\bdar{\hat{\vp}}     &                   &\bdar{\vp}\\
\ttop(M) &\lra \cale(M) \lra &\pi_0\left(\cale(M)\right)\\
\end{array}
$$
commutative.

Suppose $F$ acts on $M$.  Let $F^*$ be the group of all liftings of elements of
$F$ to the universal covering $\wt M$.  Then
$$1\ra\pi\ra F^*\ra F\ra 1$$
is exact, and is called the \emph{lifting sequence} of the action $(F,M)$.
Furthermore, for an aspherical $M$, if $(F,M)$ is effective, then
{\it the centralizer of $\pi$ in $F^*$, $C_{F^*}(\pi)$, is torsion-free}
(such an extension was called \emph{admissible} in \cite{lr81-1}); 
and $F^*$ lies in $\ttop(\wt M)$.

Thus we have a necessary condition for the existence of a lifting of an
abstract kernel $\vp\col F\ra\out(\pi)$, as an (effective, resp.) group
action:  the existence of an (admissible, resp.) group extension of 
$\pi$ by $F$ realizing the abstract kernel.  For finite groups, this necessary 
condition is also sufficient for some tractable manifolds.
\bigskip

No examples of closed aspherical manifolds are yet known where the 
existence of a group extension does not also yield a group 
action.
\bigskip

\begin{defn}{\rm
Let $Q$ act properly on an admissible space $W$ (see section \ref{topgp}) 
and $B$ be the quotient
$Q\bs W$. Suppose for each extension $1\ra Q\ra E\ra F\ra 1$ by a finite
group $F$, the action of $Q$ extends to a proper action of $E$ on $W$.
Then we say that the $Q$ action on $W$ is \emph{finitely extendable}.
In particular, then $F$ acts on $B$ preserving the orbit structure.

If $\gam$ is normal in $\pi$, let $\aut(\pi,\gam)$ denote the automorphisms
of $\pi$ that leave $\gam$ invariant.
Since $\inn(\pi)$ leaves $\gam$ invariant, we can put
$\aut(\pi,\gam)/\inn(\pi)= \out(\pi,\gam)$. It is a subgroup of
$\out(\pi)$. See notation in section 1.

Let $M=\pi\bs\wt{M}$ be an orbifold (with the orbifold group $\pi$).
Call a finite abstract kernel $F\ra\out(\pi)$
\emph{geometrically realizable}
if it can be realized as an action of $F$ on $M$.
}\end{defn}

We are interested in realizing a finite abstract kernel $F\ra\out(\pi)$ 
as a group action on a Seifert fiber space $M$ with a typical fiber
$\gam\bs G$. In particular, we want the $F$ action to be fiber-preserving
maps. This means that, on the group level, the extension must leave the
lattice $\gam$ invariant. In other words, we consider only those abstract
kernels which have images in $\out(\pi,\gam)$.

\begin{thm}
\label{thm-5-5-2}
Let $\gam$ be a lattice in the special Lie group $G$ 
(see subsection \ref{special}),  and
$\rho: Q\ra\ttop(W)$ be a proper action on $W$ which is finitely
extendable.
Let $1\ra\gam\ra \pi\ra Q\ra 1$ be an extension 
and $\theta:\pi\ra\topgp$ be any homomorphism.
Then each abstract kernel $\varphi: F\ra\out(\pi,\gam)$  of a finite group 
$F$ can be geometrically realized as a group of Seifert automorphisms
on $M(\theta(\pi))=\theta(\pi)\bs(G\x W)$ if and only if the abstract
kernel $\vp$ admits some extension.
\end{thm}

\begin{proof}
Let  $1\ra\pi\ra E \ra F\ra 1$ be an extension realizing the abstract
kernel $\vp$.
Consider the induced extension
$$
1\lra Q\lra E/\gam \lra F\lra 1.
$$
Since $\rho$ is finitely extendable, 
there exists $\rho':E/\gam\ra\ttop(W)$ extending $\rho:Q\ra\ttop(W)$.
We have the commutative diagram with exact columns and rows:
$$
\begin{array}{ccccccccc}
  &     &              &      &1            &      &1     &      & \\
  &     &              &      &\bda         &      &\bda  &  &\\
1 &\lra &\gam          &\lra  &\pi          &\lra  &Q     &\lra  &1\\
  &     &\bda          &      &\bda         &      &\bda  &      & \\
1 &\lra &\gam          &\lra  &E            &\lra  &E/\gam     &\lra  &1\\
  &     &              &      &\bda         &      &\bda&  &\\
  &     &              &      &F            &\stackrel{=}\lra  &F     &      & \\
  &     &              &      &\bda         &      &\bda  &  &\\
  &     &              &      &1            &      &1     &      & \\
\end{array}
$$
Again by the existence theorem for special lattices, Theorem \ref{const-thm},
there exists $\theta': E\ra\topgp$,
where $\theta'|_{\gam}=\theta|_{\gam}=i: \gam\hra G$, and $\rho'|_Q=\rho$.
Put $\theta'|_{\pi}=\theta'$.
Of course, $\theta'$ may be different from $\theta$, but as $\theta$ and
$\theta'$ agree on $\gam$ and $Q$, we can apply Theorem \ref{const-thm} (2)
to conjugate $\topgp$ by an element of $\mwg\rx\inn(G)$ which carries
$\theta'|_{\pi}$ to $\theta$ so that the new homomorphism $\theta':E\ra\topgp$
is an extension of $\theta:\pi\ra\topgp$.
This yields an action of $F$ on $\theta(\pi)\bs(G\x W)$ as a group of 
Seifert automorphism as desired.
\end{proof}

Since we used the uniqueness of the Seifert Construction, we use
the modified Seifert Construction modelled on $G/K\x W$ instead of $G\x W$
in the semi-simple case.

\begin{cor}
Let $M(\pi)$ be an infra-$G$-manifold with $G$ special, and with
 fundamental group $\pi$.
Suppose $\vp: F\ra\out(\pi)$ is an abstract kernel. Then there exists
a geometric realization of $F$ acting on $M(\pi)$ by affine diffeomorphisms
 if and only if there is an extension $E$ of $\pi$ by $F$ which realizes
 this abstract kernel.
\end{cor}

\begin{proof}
The infra-$G$-manifold $M(\pi)$ is modelled on $G\x\{p\}$
($p$=point)\quad ($G/K\x\{p\}$ for a convenient form of $G$ in the
semi-simple case)\quad and $\ttop_G(G\x \{p\})=\aff(G)$
(resp. $\ol{\aff}(G,K)$), see \cite{lr96-1} for this notation.
Trivially, every $Q\ra\ttop(\{p\})$ extends to $E/\gam\ra\ttop(\{p\})$.
The above theorem then immediately applies, and $F$ acts on $M(\pi)$ by Seifert
automorphisms which are affine diffeomorphisms.
\end{proof}

\begin{rem}{\rm
1. Since we may introduce a metric structure in the Corollary from a left
invariant metric on $G$ (see section \ref{infra-homo}),
$M(\pi)$ has the structure  of a flat, almost flat, Riemannian
infra-solvmanifold or a locally symmetric spaces.
We may also further conjugate $\theta(\pi)$ in $\aff(G)$ so that $F$ now
acts on the conjugated manifold by isometries preserving the flat, etc,
structures.

2. The proper action of $\pi$ in the Theorem is not necessarily free nor
effective.
Thus $M(\pi)$ could very well be a Seifert orbifold. The Corollary
then works for such orbifolds, i.e., infra-$G$-spaces.
In the Euclidean case, $M(\pi)$ would then be a Euclidean ``crystal'' and 
$\pi$ a Euclidean crystallographic group.
In Theorem \ref{thm-5-5-2}, $F$ sends fibers to fibers.

For more about the realizations up to strict equivalences, and finding
examples where $F$ does not lift because there are no extensions realizing
the abstract kernels, the reader is referred to 
\cite{klr83-1}, \cite{lr96-1},
\cite{lr81-1}, \cite{zz79-1}, \cite{lr82-1}, \cite{rs77-1},
\cite{lee82-1}, \cite{lee82-2}, \cite{igod84-2}, \cite{sy79-1} and
\cite{raym79-1}.
}\end{rem}

\subsection{Seifert fiberings with codimension 2 fibers}
\begin{step}{\rm
There is a class of Seifert fiberings which are very close generalizations
of the classical Seifert 3-manifolds. Let the discrete group $Q$ act
effectively and properly on the topological plane $\bbr^2$.
Then $Q$ can be conjugated in $\ttop(\bbr^2)$ to be a  group acting as
isometries on the hyperbolic plane or the Euclidean plane.
In the first case, call $Q$ hyperbolic, and in the latter Euclidean.
If $Q$ preserves orientation, then $Q$ can be conjugated to act as
holomorphic actions on $\bbc$ or the unit disk $D$.
This latter fact, using standard techniques in transformation group
theory, uniformization and topology of 2-manifolds, is not very
difficult to prove.
(Certainly, if $Q$ acts simplicially, the arguments are not hard).
When the orbit space $Q\bs\bbr^2$ is compact, then a theorem,
essentially due to Nielsen, says that any two $Q$ actions are
topologically conjugate.
If the actions are assumed smooth, then they are smoothly conjugate.
In the non-compact case, there are similar results but they can be
very complicated if $Q$ is not finitely presentable. So for our next
example, we restrict to $Q$ for which $Q\bs\bbr^2$ is compact.
Since $Q\bs\bbr^2$ will have the structure of an orbifold, $Q$ will
either be isomorphic to a Euclidean crystallographic group or a
cocompact hyperbolic group.
In the Euclidean case, these are the 17 wall-paper groups. They are
centerless except for $Q=\bbz\x\bbz$ or $\bbz\rx\bbz$ (the fundamental
group of the Klein bottle).
In the hyperbolic case, all the groups are not solvable, and have no
non-trivial normal abelian subgroups.
}\end{step}

\begin{step}
\label{5.6.2}
{\rm
An important class of closed 3-manifolds are those closed 3-manifolds
whose fundamental group contains a normal subgroup $\bbz$.
This class contains all the classical Seifert 3-manifolds with
infinite fundamental group as well as some other non-orientable
manifolds such as closed 3-manifolds admitting a circle action without
fixed points but having 2-dimensional submanifolds whose stabilizer is
$\bbz_2$.
The orientable ones, however, coincide with the classical orientable
Seifert 3-manifolds with infinite fundamental group as defined by
Seifert in \cite{seif33-1}.
Of course, these statements are to be taken modulo the Poincar\'e
Conjecture, for if there were a fake 3-sphere $\Sigma$, then
$M\sharp\Sigma$, where $M$ is a Seifert manifold, would not be a Seifert
manifold but would still have a normal $\bbz$ in its fundamental group.
In any case, a closed 3-manifold whose fundamental group contains a
normal $\bbz$ subgroup, has, modulo the Poincar\'e Conjecture, an
orientable cover admitting an injective $S^1$-action. 
Consequently, they all
have a universal cover that splits into $\bbr\x\bbr^2$ or
$\bbr\x S^2$, where the $\bbr$ descends to the fibers associated with
the normal $\bbz$.
In the $\bbr\x S^2$ case, the manifolds have an orientable cover
isomorphic to $S^1\x S^2$.
Thus, their fundamental groups are just finite extensions of $\bbz$.
We are interested in those that are not finite extensions of $\bbz$.
They are all aspherical since they are covered by $\bbr^3$.
\bigskip

Therefore, for each extension,
$1\ra\bbz\ra\pi\ra Q\ra 1$, where $\rho:Q\ra\ttop(\bbr^2)$ is an
effective proper action on $\bbr^2$, we may construct model Seifert
manifolds. In fact, as $\rho$ is topologically rigid, we can assume
that after conjugating in $\ttop(\bbr^2)$, $\rho$ maps $Q$ into either
$\isom(\bbe^2)=\bbr^2\rx O(2)$ or $\isom(\bbh)=\psl\rx\bbz_2$,
the isometries of the Euclidean plane or the hyperbolic plane.
In this case, as we shall see in section \ref{pslt-geom},
the constructed
Seifert manifolds possess a geometric structure inherited from a
subgroup $\calu$ of $\ttop_{\bbr}(\bbr\x\bbr^2)$ into which
$\pi_1(M)$ is embedded. The group $\pi$ acts freely on
$\bbr\x\bbr^2$ if and only if $\pi$ is torsion free. Otherwise, the
map $\bbr\x\bbr^2 \ra\theta(\pi)\bs(\bbr\x\bbr^2)$ is a non-trivially
branched covering and the underlying topological space of the
3-dimensional orbifold, which is still a 3-manifold with possible
boundary, may not even be aspherical (see the next section).
To obtain all possible orbifolds, we have to enlarge $Q$ to $Q'$, where
$Q$ is of index 2 in $Q'$ and the $\bbz_2$ in $Q'$ does not act effectively,
(it injects into $\aut(\bbz)$). That is,
$\vp\x\rho:Q'\ra\aut(\bbz)\x\ttop(\bbr)\ra\aut(\bbr)\x\ttop(\bbr^2)$
will be injective.
\bigskip

Obviously the existence of the Seifert Construction for every torsion
free extension $1\ra\bbz\ra \pi\ra Q\ra 1$ produces all the known
aspherical 3-manifolds with a normal $\bbz$ in their fundamental group.
Furthermore, rigidity of the construction, since $\bbz$ will almost
always be characteristic, says that any two constructed manifolds with
isomorphic fundamental groups are diffeomorphic.
In fact, except for several exceptions (less than 10), the
isomorphisms must be given by a Seifert automorphism.
(The exceptions can occur for the 3-torus, for example, where the
$\bbz$ of the fiber is not the only possible normal $\bbz$).
That no other manifolds are possible except those given by the Seifert
Construction follows from results of Waldhausen in the Haken case. For
non-Haken manifolds, the result is more recent and is due to Scott,
Casson and others.

The rigidity result is also true in the orbifold case. For example, if
$Q$ is \emph{orientation-preserving} and $Q\ra\aut(\bbz)$ is trivial,
then the $\bbr$ action descends to an $S^1$-action on
$\theta(\pi)\bs(\bbr\x\bbr^2)$.
The group $\bbz$ is the maximal normal abelian subgroup of $\pi$ in
case $Q$ is hyperbolic since $Q$ has no normal abelian subgroups.
In the Euclidean crystallographic case where $Q\neq\bbz^2$, $Q$ is
centerless and so $\bbz$ will be the entire center of $\pi$ and so is
characteristic.
For $Q=\bbz^2$, the $\pi$ are torsion free and are distinguished by
$H_1(\pi;\bbz)=\bbz^2+\bbz_b$.
So, except for 3-torus, \emph{the rigidity theorem for the Seifert
Construction implies two such oriented orbifolds with isomorphic orbifold
groups will be diffeomorphic via an $S^1$-equivariant
orientation-preserving homeomorphism}.
}\end{step}

\begin{thm}[Higher dimensional fibers]
\label{higher}
Let $Q$ be hyperbolic and act properly on $\bbr^2$ with compact
quotient. Let $\gam$ be a special lattice of type {\rm (S1), (S2) or (S3)}
(see subsection \ref{special}),
and
$$
1\lra\gam\lra\pi\lra Q\lra 1
$$
be an extension. Then there exists a homomorphism 
$\theta:\pi\ra\ttop_G(G\x\bbr^2)$
unique up to conjugation. $\pi$ acts freely if and only if $\pi$ is
torsion free and, consequently, $M(\theta(\pi))$ will then be
aspherical. In general, $M(\theta_1(\pi_1))$ is homeomorphic via a Seifert
automorphism to $M(\theta_2(\pi_2))$ if and only if $\pi_1$ is isomorphic
to $\pi_2$.
\end{thm}

\begin{proof}
We know $\theta$ exists by the existence theorem for special lattices.
Since $Q$ has no non-trivial normal solvable subgroups, and $G$ is
solvable, $\pi\cap G=\gam$.
Therefore, $\gam$ is the maximal normal solvable subgroup in $\pi$ and
consequently is characteristic. Since $\rho:Q\ra\ttop(\bbr^2)$ is
rigid, the Seifert Construction is rigid and the theorem follows.
(cf. \cite[section 12]{cr71-3}).
\end{proof}

Nicas and  Stark in \cite{ns85-1}, using controlled surgery, have shown 
that any closed aspherical topological manifold whose fundamental
group is a central extension of $Z^s$ by a 2-dimensional orientation
preserving  cocompact hyperbolic $Q$  is homeomorphic to the corresponding
Seifert manifold constructed in Theorem \ref{higher},
provided the dimension of the manifold is greater than $4$.

\begin{cor}
Let $\vp:F\ra\out(\pi)$ be an abstract kernel, where $F$ is finite and
$\pi$ is as in the above theorem. Then $\vp$ has a geometric
realization of an action of $F$  as a group of Seifert automorphism
on $M(\theta(\pi))$ if and only if
there exists some extension $E$ of $\pi$ which realizes the abstract
kernel.
\end{cor}

\begin{proof}
The extension $E$ induces a short exact sequence
$1\ra Q\ra E/\gam\ra F\ra 1$.
S. Kerckhoff \cite{kerckhoff83-1} has shown, in  his solution to the 
Nielsen Realization Problem,
that a finite extension of a hyperbolic $Q$ with an injective abstract 
kernel is again a hyperbolic group.
If $K$ denotes the kernel of the abstract kernel $F\ra\out(Q)$, then
$(E/\gam)/K$ is a hyperbolic group containing a topological conjugate,
in $\ttop(\bbr^2)$, of $Q$.
Then $E/\gam$, via $E/\gam\ra (E/\gam)/K$, acts properly on $\bbh$ and
(up to conjugation in $\ttop(\bbr^2)$) extends the $Q$ action.
Therefore Theorem \ref{thm-5-5-2} applies.
\end{proof}

\subsection{The  classical Seifert 3-dimensional manifolds}
In a classical paper \cite{seif33-1} in 1933, H.~Seifert described a class
of 3-dimensional manifolds that have turned out to be fundamental for the 
topology of 3-manifolds. We assume the reader has some familiarity with 
these manifolds.

	In the closed orientable (resp. non-orientable) case Seifert's
manifolds whose fundamental groups are  infinite, but not a
finite extension of $\bbz$,  coincide with (resp. are a subset of) those
closed orientable (resp. non-orientable) manifolds, $\cals\calf$,  that can be 
constructed, as in Theorem \ref{const-thm},
from a  torsion free  extension,
$1 \ra Z \ra \pi \ra Q\ra 1$, with $Q$ acting properly, effectively,
and cocompactly on $\bbr^2$. 

	It has recently been shown that if there are no fake homotopy
3-spheres, then the closed 3-manifolds, whose fundamental groups contain a 
normal $\bbz$ subgroup and which are not  finite extensions  of $\bbz$, 
coincide with $\cals\calf$.
See section \ref{5.6.2}.
%---end of section 5 ----------------------------------------------------

\section{Reduction of the Universal Group}
\bigskip

\subsection{Purpose and Requirements}
In a general Seifert Construction, it may happen that $\rho:Q\ra\ttop(W)$
has an image in a subgroup that has geometric or topological significance.
This means that it is likely that the associated Seifert Constructions
inherit some of these properties. 
Let $\calu$ be a subgroup of $\topgp$ containing $\ell(G)$.
Put $\wh{\calu}=\calu\cap(\mwg\rx\inn(G))$ and $\ol{\calu}=\calu/\wh{\calu}$
so that
\begin{equation}
\label{x-1}
{\small
\begin{array}{ccccccccc}
1&\lra &\wh{\calu}     &\lra &\calu  &\lra &\ol{\calu}        &\lra &1\\
 &     &\cap           &     &\cap   &     &\cap              &     & \\
1&\lra &\mwg\rx\inn(G) &\lra &\topgp &\lra &\out(G)\x\ttop(W) &\lra &1\\
\end{array}
}
\end{equation}
commutes. 

If the image of $\theta:\pi\ra\topgp$ lies in $\calu$, we say that the
universal group has been reduced to $\calu$ for $\pi$.
The group $\calu$ can then be used to study extra structures on the Seifert
fiber space $\theta(\pi)\bs(G\x W)$.

We have already seen in section \ref{meaning} that if $W$ is a smooth
manifold, $\gam$ special and $\rho:Q\ra\diff(W)$,
then any injective Seifert Construction can be done smoothly in
$$
\difgp=\cwg\rx(\aut(G)\x\diff(W)),
$$
where $\difgp=\topgp\cap\diff(G\x W)$.

Furthermore, existence, uniqueness and rigidity also hold because
we may use a differential partition of unity in the proof for the vanishing
of the  necessary cohomology groups.

To obtain existence of an injective Seifert Construction in $\calu$, one
has to verify the conditions in Theorem \ref{main-alg} in $\calu$.
If $\calu$ is much smaller than $\topgp$, uniqueness and rigidity are
not likely to hold.
The set of all homomorphisms of $\pi$ into $\calu$ which belong to the same
conjugacy class in $\topgp$ is a deformation space of that conjugacy class.
We shall illustrate these concepts by describing in some detail the
deformation spaces for interesting geometric situations.
\bigskip

\subsection{Injective Holomorphic Seifert Fiberings}
\label{holomorphic}
We assume that $G=\bbc^k$, $W$ is a complex manifold, $\rho:Q\ra\ttop(W)$
has image in the holomorphic homeomorphisms of $W$, and 
$\calh(W,\bbc^k)\subset\mw{\bbc^k}$ are the holomorphic maps.

We also assume that the map $\bbc^k\x W\ra W$ is holomorphically trivial and
so $\calu=\hol_{\bbc^k}(\bbc^k\x W)$ then becomes
$$
\calh(W,\bbc^k)\rx(\gl(k,\bbc^k)\x\hol(W)),
$$
where $\hol(W)$ is the group of holomorphic automorphisms of $W$.
Existence, uniqueness and rigidity do not necessarily hold because we do
not have holomorphic partitions of unity and the groups
$H^i(Q;\calh(W,\bbc^k))$ does not vanish in general.

The reader is referred to \cite{cr71-3} where a general and comprehensive
theory of holomorphic Seifert fiberings whose universal space is a
holomorphic fiber bundle over $W$ with fiber a complex torus or $\bbc^k$ is
given. We shall restrict ourselves here to a special case closely related
to the classical Seifert 3-manifolds.

Let $(\cs,M)$ be an injective, proper, holomorphic $\cs$ action on a
complex 2-manifold $M$ so that the quotient space is compact.
As in the case of an injective $S^1$ action, an injective
proper $\cs$ action lifts to the covering space of $M$ corresponding to the 
image of the evaluation  homomorphism, and yields  a  splitting 
$(\cs,\cs\x W)$ where $W$ is
a simply connected complex 1-manifold, (cf. section \ref{torus-action}).
Therefore, $W$ is $\bbc$, $D$ the open unit disk, or $\bbc P_1$.
We shall restrict ourselves to $W$ being the unit disk $D$. The orbit space
$Q\bs W$ is a closed Riemann surface.
The action of $Q=\pi_1(M)/\bbz$ on $D$ is holomorphic, \pd, but not
necessarily free.
Therefore, $M\lra\cs\bs M=Q\bs W$ is a generalization of a principal
holomorphic $\cs$-bundle over a Riemann surface.

From the exact sequence $1\ra\bbz\ra\bbc\stackrel{\exp}\lra\cs\ra 0$,
we obtain the exact sequence
$$
0\lra \bbz=M(W,\bbz)\lra\calh(W,\bbc)\lra\calh(W,\cs)\lra 0
$$
which gives rise to a long exact sequence of cohomology groups
$$
\cdots\stackrel{\delta^{i-1}}\lra H^i(Q;\bbz)\lra H^i(Q;\calh(W,\bbc))\lra
H^i(Q;\calh(W,\cs))\stackrel{\delta^i}\lra\cdots
$$
The group $Q$ acts on the unit disk $D$ as a cocompact Fuchsian group. That
is, $\rho:Q\ra\hol(D)$, the complex automorphisms of the unit disk.
The action of $Q$ on $\lambda\in\calh(W,\bbc)$ is given by
$$
^\alpha\lambda=\lambda\circ (\rho(\alpha))\inv.
$$

Let us compare the smooth situation with the holomorphic one. We have the
following commutative diagram of exact sequences
\begin{equation}
\label{x-2}
{\small
\begin{array}{ccccccccccccc}
0 &\lra &H^1(Q,\bbz)&\lra &H^1(Q,\hw{\bbc})&\lra &H^1(Q,\hw{\cs})
&\stackrel{\delta}\lra  &H^2(Q,\bbz) &\lra &H^2(Q;\hw{\bbc})\\
 &     &\bdar{=}  &     &\bda   &     &\bda  &     &\bdar{=} & &\bda \\
 &\lra &H^1(Q,\bbz)&\lra &H^1(Q,\cw{\bbc})&\lra &H^1(Q,\cw{\cs})
&\stackrel{\delta}\lra  &H^2(Q,\bbz) &\lra &H^2(Q;\cw{\bbc})\\
\end{array}
}
\end{equation}
For the smooth case, $H^i(Q;\cw{\bbc})=0$, $i>0$,
and as we shall see, $H^2(Q;\calh(W,\bbc))=0$.

For each central extension $0\ra\bbz\ra\pi\ra Q\ra 0$ represented by
$[f]\in H^2(Q;\bbz)$, we have smooth Seifert Constructions
$\theta:\pi\ra\diff_{\bbc}(\bbc\x D)=\diff_{\bbc}(\bbc\x \bbr^2)$.
If we fix $i:\bbz\ra\bbc$ and $\rho:Q\ra\diff(D)$, the construction is
unique up to strict equivalences (subsection \ref{uniqueness}).
We have the smooth Seifert orbifold over $Q\bs W=Q\bs D$ with an induced 
$\cs$ action and therefore an $S^1$ action on $\theta(\pi)\bs(\bbc\x D)$  
($\cong\bbr^1\x \theta(\pi)\bs(\bbr^1\x D)=\bbr^1\x N^3$ since $\cs$ splits
smoothly as $\bbr^1\x S^1$).
The uniqueness says that for any other embedding
$\theta':\pi\ra\diff_{\bbc}(\bbc\x D)$, keeping $i$ and $\rho$ fixed, the
$\cs$ action on $\theta'(\pi)\bs (\bbc\x D)$ is strictly smoothly 
equivalent to that on $\theta(\pi)\bs (\bbc\x D)$.
\bigskip

For the same $\pi$, we have the homomorphism
$H^2(Q;\bbz)\lra H^2(Q;\calh(D,\bbc))$. The second group fortunately
can be identified with the second cohomology of the sheaf of germs of
holomorphic functions over $Q\bs D$. This vanishes since $Q\bs D$ is
(complex) 1-dimensional and the sheaf is coherent (i.e., locally free).
This means that $[f]\in H^2(Q;\bbz)$ maps to $0\in H^2(Q;\hw{\bbc})$. 
But as the groups
are abelian, this becomes exactly the identity (\ref{cond-1}), and we have
$\theta:\pi\ra\hol_{\bbc}(\bbc\x W)$.
Therefore, each $[f]$ has holomorphic realizations for each fixed $i$ and
$\rho$.
Recall from Theorem \ref{main-alg}, the set of all $\theta:\pi\ra
\hol_{\bbc}(\bbc\x D)$ with fixed $i: \bbz\ra\bbc$ and $\rho:Q\ra\calh(D)$, up
to conjugation by elements of $\calh(D,\bbc)$, is in one--one 
correspondence with $H^1(Q;\calh(D,\bbc))$. 
(This complex vector space is the same as 
$H^1(V;h^0_{\bbc})$, the first
cohomology of the sheaf of germs of holomorphic functions where $V$ is
treated as the analytic space $V=Q\bs D$.
That is, for each open $U$ in $V$, we consider $p\inv(U)$ and holomorphic
functions $\lambda: p\inv(U)\ra \bbc$ such that
$\lambda(\rho(\alpha)(w))=\lambda(w)$, $w\in p\inv(U)$ and $p: D\ra Q\bs D$
is the projection.
This defines the sheaf $h^0_{\bbc}$ over $V$).
This group is isomorphic to $\bbc^g$, where $g$ is the genus of $V$.

\begin{thm}[{\cite[\S13]{cr71-3}}]
For each smooth action $(\cs,M)$ corresponding to the unique strict
conjugacy class $\theta(\pi)$, with $\pi$ torsion free,  
there exists a complex $g$-dimensional
family of strictly holomorphically inequivalent $\cs$ actions each strictly
smoothly equivalent to the smooth $(\cs,M)$.
\end{thm}

\begin{proof}
We may interpret $[f']\in H^1(Q;\calh(W,\cs))$ to represent the 
effective holomorphic
$\cs$ action on $\theta(\pi)\bs(\bbc\x D)$ up to strict $\cs$ equivalence.
If $Q$ were torsion free, then $Q$ acts freely and $Q\bs D=V$ is a closed
oriented surface without branch points.
The $\cs$ action is then free and proper yielding a principal holomorphic
$\cs$ bundle, corresponding to a complex line bundle over $V$.
Since $Q$ is not assumed to be torsion free, $f'$ determines the
holomorphic $\cs$ action (cf, \cite[\S5]{cr71-3}).
Given two injective $\cs$ actions $[f']$ and $[f'']$ with the same ``Chern
class'' $\delta[f']=\delta[f'']$,
there exists a $[\lambda]\in H^1(Q;\calh(D,\bbc))$ such that
$[f'']=[f']+[\lambda]$. 
Since $H^1(Q;\calh(D,\bbc))$ is a vector space of complex dimension $g$,
there exists a whole $g$-dimensional family of inequivalent holomorphic $\cs$
actions starting from $[f']$ and ending with $[f'']$.
\end{proof}

Returning to diagram \ref{x-2}, define
$$
\begin{array}{rcl}
\pic(Q\bs D) &\ = &\ H^1(Q;\calh(D,\bbc))/\on{image}(H^1(Q;\bbz))\\
             &\ = &\ {\rm a\ complex\ } g{\rm -torus,\ or\ real\ } 
               2g{\rm -torus}, T^{2g}.
\end{array}
$$
The connected component of the 
group $H^1(Q;\calh(D,\cs))$, in the $Q$ torsion free case, is called
the \emph{Picard group} for the line bundles over $Q\bs D$.
In our case, $H^1(Q;\calh(D,\cs))$
is isomorphic to $T^{2g}\oplus\bbz\oplus{\rm finite\ torsion}$.

We obtain the exact sequence
$$
0\lra\pic(Q\bs D)\lra H^1(Q;\calh(D,\cs))\stackrel{\delta}\lra H^2(Q;\bbz)
\lra 0,
$$
where the middle group is the isomorphism classes of injective holomorphic
$\cs$ actions over $Q\bs D$, $\delta$ sends such an isomorphism class to
its ``Chern class'' and $\pic(Q\bs D)$ represents the deformations.
As before, $H^2(Q;\bbz)\cong \bbz\oplus{\rm Torsion}$.
\bigskip

In the discussion above, we have fixed $i:\bbz\ra\bbc$ and
$\rho:Q\ra\hol(D)$.
If we vary these choices, we don't get anything new in the smooth case
because of smooth  rigidity. That is, $\theta(\pi)$ is conjugate to
$\theta'(\pi)$ in $\diff_{\bbc}(\bbc\x D)$ where conjugation is taken in
the whole group and not just in $\calc(D,\bbc)$ as for strict equivalence.
However, in the holomorphic case, a change in $\rho:Q\ra\hol(D)$ induces a
much larger deformation space than treated above.
We can see this in our next example of reduction of the universal group
where instead of considering complex structures and complex actions, we
replace them by essentially equivalent Riemannian metric structures and
metric preserving $S^1$-actions on $N^3, (M=N\x\bbr^1)$.

\subsection{$\pslt$-geometry}
\label{pslt-geom}
\begin{step}{\rm
When $W$ is homeomorphic to $\bbr^2$ in a Seifert Construction and $\rho(Q)$
is a discrete subgroup of $\ttop(\bbr^2)$, acting properly on $\bbr^2$, then
the group $\rho(Q)$ can be conjugated in $\ttop(\bbr^2)$ to a group which 
acts as
isometries on $\bbr^2$ with the usual Euclidean metric (e.g., $\rho(Q)$ is
crystallographic) or as isometries on $\bbr^2$ with the usual hyperbolic
metric. In the former case, we say $\rho(Q)$ is isomorphic to a (Euclidean)
crystallographic group, and in the latter, $\rho(Q)$ is isomorphic to a
hyperbolic group.
We write $\bbr^2$ with the usual Euclidean metric (resp. hyperbolic metric)
as $\bbe^2$ (resp. $\bbh$).
If $Q\bs \bbr^2$ is compact, then any two embeddings of $Q$ are
conjugate in $\ttop(\bbr^2)$.
However, if $\rho_1,\rho_2:Q\ra\isom(\bbh)$ have compact quotients, then
$\rho_1$ is conjugated in $\isom(\bbh)$ to $\rho_2$ if and only if 
$\rho_2(Q)$ lies in the normalizer of $\rho_1(Q)$ in $\isom(\bbh)$,
(The normalizer is a finite extension of $\rho_1(Q)$).
Thus, if we reduce the universal group $\ttop_G(G\x\bbr^2)$ to $\calu$,
where at least $\ttop(\bbr^2)$ is replaced by $\isom(\bbe^2)$ or
$\isom(\bbh)$, we would expect to find a rich deformation theory for
Seifert Constructions modelled on $G\x\bbr^2$.
We will take $G=\bbr$, and $\calub=\isom(\bbh)\subset\ttop(\bbr^2)$.
}\end{step}

For each central extension $0\ra\bbz\ra\pi\ra Q\ra 0$ with $Q$ cocompact
hyperbolic, the group $\pi$ can be embedded by $\theta:\pi\ra
\ttop_{\bbr}(\bbr\x\bbh)$ so that it is topologically and/or smoothly
rigid.
If $\pi$ is torsion free, $\pi\bs(\bbr\x\bbh)$ is a classical closed
Seifert 3-manifold $N^3$ and with a unique $S^1$ action up to equivalence.
The $S^1$-orbit space or base space is a 2-dimensional orbifold isomorphic
to $Q\bs\bbh$.
The product $\bbr\x\bbh$ carries several geometries so that the
Riemannian metric induced on $\bbr\x\bbh\ra\bbh$ is the
hyperbolic metric.
We will examine first the Riemannian metric on $\pslt$ and then later an
$\bbr$-invariant Lorentz metric on $\bbr\x\bbh$.
\bigskip

Let us denote $\pslt$ by $\tp$.
The space $\tp$ is the universal covering group of $\psl$.
Topologically $\tp$ is homeomorphic to $\bbr\times\bbh$.  
$$
\begin{array}{ll}
  P&=\psl\\
\tp&=\pslt\\
   &\approx\bbr\times\bbh.
\end{array}
$$
Therefore, if we use 
$$
\ttop_{\bbr}(\bbr\times\bbh)=\on{M}(\bbh,\bbr)\rx(\gl(1,\bbr)\times\ttop(\bbh))
$$ 
as our universal group, we will not be able to distinguish the geometries
of $\bbr\times\bbh$, $\bbr\x\bbr^2$ and $\nil$ from that of $\tp$.
\medskip

The Lie group $\psl$ can be viewed as the unit tangent bundle of the
hyperbolic space $\bbh$, and it has a natural Riemannian metric.
This metric pulls back to a Riemannian metric on $\tp$.
It turns out that this metric is right invariant. That is, all right
translations by elements of the group are isometries.
Furthermore, the isometry group is
$$
\isom(\tp)=(\bbr\x_{\bbz}\tp)\rx\bbz_2,
$$
where $\bbr$ is a subgroup of $\tp$ containing the center $\bbz$,
acting as left translations.
These two actions commute with each other, and $\wh{\ell}(z)=\wh{r}(z\inv)$
for $z\in\bbz$, the center of $\tp$.
The finite group $\bbz_2$ is generated by the reflection about the
$y$-axis. In fact, any orientation-reversing isometry of period 2 will do.
While it reverses the orientation of the base space $\bbh$, it also
reverses the orientation of the fiber $\bbr$. Consequently, it preserves
the orientation of $\tp$.

We take the subgroup $\bbr$ described above as our $G$. 
Smoothly, $\tp=\bbr\x\bbh$. Therefore, we have
$$
G=\bbr,\quad W=\bbh.
$$
Since $\bbr\x_{\bbz}\tp$ commutes with the left translation $G=\bbr$, 
and the generator of $\bbz_2$ is an inversion of $\bbr$,
$\isom(\tp)=(\bbr\x_{\bbz}\tp)\rx\bbz_2$ lies inside $\ttop_{\bbr}(\bbr\x\bbh)$.
To make the presentation clearer, we use only the connected component of
$\isom(\tp)$. So, let's take
$$
\calu\ =\ \isom_0(\tp)\ =\ \bbr\x_{\bbz}\tp
$$
so that we have the commuting diagram:
$$
\begin{array}{ccccccccc}
1&\lra &\bbr           &\lra &\calu=\bbr\x_{\bbz}\tp  &\lra &P &\lra &1\\
 &     &\cap           &     &\cap   &     &\cap              &     & \\
1&\lra &\on{M}(\bbh,\bbr) &\lra &\ttop_{\bbr}(\bbr\x\bbh)&\lra 
&\gl(1,\bbr)\x\ttop(\bbh) &\lra &1\\
\end{array}
$$
For this case, $Q$ coincides with the $Q$ examined in section
\ref{holomorphic}, and has a well known presentation
$$
Q=
\langle \ol{x}_1,\cdots,\ol{x}_g,
\ol{y}_1,\cdots,\ol{y}_g,
\ol{w}_1,\cdots,\ol{w}_p \mid
{\ol{w}_j}^{\alpha_j}=1,\quad\prod_{j=1}^p \ol{w}_j
\prod_{i=1}^g [\ol{x}_i,\ol{y}_i]=1\rangle
$$
for $p\ge 0$, $g\ge 0$ and all $\alpha_j \ge 2.$
It is also required that the Euler characteristic of $Q$, defined as
$$
\chi(Q)=(2 - 2g) - \sum_{j=1}^p \left(1 - {1\over{\alpha_j}}\right),
$$
satisfies $\chi(Q)<0$. It is our intention to characterize those
$[\pi]\in H^2(Q;\bbz)$ which embed in $\calu$ and to determine their
deformation spaces.

\begin{step}{\rm
For any subgroup $Q$ of $\isom_0(\bbh)$, 
one can pullback (see subsection \ref{pullback-ext}) the above extension via 
$Q\hra\isom_0(\bbh)$ to get $\wt{Q}$ so that the diagram
$$
\begin{array}{cccccccccccc}
1 &\lra &\bbr     &\lra &\wt{Q} &\lra &Q          &\lra &1\\
  &     &\bdar{=} &     &\bda   &     &\bdar{\cap}&     &\\
1&\lra  &\bbr     &\lra &\isom_0(\tp)&\lra &\isom_0(\bbh)&\lra &1
\end{array}
$$
commutes. Thus, $\bbr$ becomes a trivial $Q$-module.
}
\end{step}

\begin{lem}
\label{q-tilde}
Let $Q$ be a cocompact discrete subgroup of $\isom_0(\bbh)$. Then
\begin{enumerate}
\item $H^2(Q;\bbr)=\bbr$,
\item The class $[\wt{Q}]\in H^2(Q;\bbr)$ is non-zero. 
(That is, the exact sequence does not split).
\end{enumerate}
\end{lem}

\begin{proof}
(1) Let  $1\ra\bbz\ra\tp\ra P\ra 1$ be the universal covering projection, 
and let $1\ra\bbz\ra\wh{Q}\ra Q\ra 1$ be the pullback of this exact
sequence via $Q\hra P$. Then $\bbz\subset \wh{Q}$ sits in
$\bbr\x_{\bbz}\wt{P}=\isom_0{\wt{P}}$ as the center. Denote the inclusion
of $\bbz\hra\bbr\subset \isom_0(\wt{P})$ by $i$. Then
$$
\begin{array}{cccccccccccc}
1 &\lra &\bbz   &\lra &\wh{Q} &\lra &Q  &\lra &1  
&\quad [\wh{Q}]\in &H^2(Q;\bbz)\\
  &     &\bdar{i} &&\bda&&\bdar{=}&&&&\bdar{i_*}\\
1 &\lra &\bbr   &\lra &\wt{Q} &\lra &Q  &\lra &1  
&\quad [\wt{Q}]\in &H^2(Q;\bbr)\\
\end{array}
$$
is commutative so that $i_*[\wh{Q}]=[\wt{Q}]$.
By Selberg's lemma, $Q$ contains a torsion free normal subgroup $Q_0$ of finite
index. Then $H^2(Q_0;\bbr)=H^2(Q_0\bs\bbh;\bbr)\cong\bbr$.
Then, by transfer, $H^2(Q;\bbr)=\bbr$, since $Q/Q_0$ acts on $Q_0\bs\bbh$
preserving orientation.
Since $H^2(Q;\bbz)$ is finitely generated, it follows that
$H^2(Q;\bbz)=\bbz\oplus\rm{Torsion}$ by the Universal Coefficient Theorem,
and the fact that $i_*:H^2(Q;\bbz)\ra H^2(Q;\bbr)$ is given by $\otimes\bbr$.
Thus the elements of infinite order inject and those of finite order are in
the kernel.

(2) We can assume, without loss of generality, $Q\subset P$ is torsion free.
In \cite{rv81-1}, it is shown that $[\wh{Q}]\in
H^2(Q;\bbz)$ is non-zero. In fact, $\wh{Q}$ is the fundamental group of the
unit tangent bundle $\wh{Q}\bs\tp$ of the surface $Q\bs\bbh$. 
In this case, the euler characteristic of $Q$, $2-2g\neq 0$, is the 
characteristic class of the principal $S^1$-bundle $\wh{Q}\bs\tp\ra
Q\bs\bbh$, and is also equal to the negative of the 
cohomology class $[\wh{Q}]\in H^2(Q;\bbz)\cong\bbz$ of the extension $\wh{Q}$.
Therefore, $[\wt{Q}]$ is non-zero in $H^2(Q;\bbr)=H^2(Q;\bbz)\otimes\bbr$. 
We also point out, in this case, $[\wh{Q}]=e(\wh{Q})$, see subsection 
\ref{euler} for definition of $e(\wh{Q})$.
\end{proof}

\begin{step}[euler number]
\label{euler}
{\rm
Let $i:\bbz\hra\bbr$ be the standard inclusion. For each central extension
$0\ra\bbz\ra\pi\ra Q\ra 1$, there is associated a rational invariant called
the \emph{euler number} of $\pi$, and is denoted by $e(\pi)$.
It can be defined in terms of a presentation of $\pi$.
Let 
$$
\begin{array}{lll}
\pi=\langle 
&\tilde{x}_1,\cdots,\tilde{x}_g,\tilde{y}_1,\cdots,\tilde{y}_g,
\tilde{w}_1,\cdots, \tilde{w}_p,\tilde{z} \mid \tilde{z}
{\textrm{\ central,\ }}\\
&{\tilde{w}_j}^{\alpha_j}={\tilde{z}}^{-\beta_j},
\prod_{j=1}^p \tilde{w}_j \prod_{i=1}^g [\tilde{x}_i,\tilde{y}_i]
={\tilde{z}}^b
\rangle.
\end{array}
$$
Then $e(\pi)=-\left(b+\Sigma{{\beta_j}\over{\alpha_j}}\right)$
and $|e(\pi)|$ is an invariant of the isomorphism class of $\pi$.
In \cite[Theorem 4.5]{llr96-1}, it is shown that under the homomorphism
$$
i_*: H^2(Q;\bbz)\lra H^2(Q;\bbr)\cong\bbr
$$
induced by $i$, we have $i_*[\pi]=L\cdot e(\pi)$, where
$L=\on{lcm}[\alpha_1,\cdots,\alpha_p]$. Thus $[\pi]$ has infinite order in
$H^2(Q;\bbz)$ if and only if $e(\pi)\neq 0$.
}\end{step}

We now characterize the cocompact orbifold groups modelled on $\tp$.

\begin{thm}
An abstract group $\pi$ can be embedded into $\isom_0(\tp)$ as a
cocompact discrete subgroup if and only if $\pi$ is a central extension of
$\bbz$ by a discrete cocompact orientation-preserving hyperbolic group $Q$ 
(so that $1\ra\bbz\ra\pi\ra Q\ra 1$ is exact) and $[\pi]\in H^2(Q;\bbz)$
has infinite order.
Further, if this is the case, the subgroup $\bbz$ is the center of $\pi$,
and in any discrete embedding, the image of $\bbz$ is $\pi\cap\bbr$,
where $\bbr\subset\bbr\x_{\bbz}\tp=\isom_0(\tp)$.
\end{thm}

\begin{proof}
Suppose $\pi$ is a cocompact discrete subgroup of $\isom_0(\tp)$. The
subgroup $\bbr$ of $\isom_0(\tp)$ is the radical (maximal connected normal
solvable subgroup) of $\isom_0(\tp)$ and
the quotient $\isom_0(\bbh)=\psl$ has no compact factor.
A theorem of Wang (see \cite[8.27]{ragh72-1}) says that the image $Q$ of
$\pi$ in $\isom_0(\bbh)$ is a lattice so that $Q$ is a 
discrete cocompact orientation-preserving hyperbolic group.
It remains to show that $\pi\cap\bbr$ is non-trivial.
Suppose not. Then $\pi$ is isomorphic to $Q$ and hence, it has 
$\bbr$-cohomological
dimension 2. However, since $\pi$ is cocompact in $\isom_0(\tp)$, its
$\bbr$-cohomological dimension is 3. This contradiction shows that
$\pi\cap\bbr=\bbz$.
Thus $\pi$ is of the form $1\ra\bbz\ra\pi\ra Q\ra 1$.
Clearly, $\bbz$ is the center of $\pi$ since $Q$ is centerless.

We shall now sketch why $[\pi]$ must be of infinite order in $H^2(Q;\bbz)$.
By Selberg's Lemma, there exists a normal subgroup $Q_0$ of $Q$ which is 
torsion free of finite index.
Let $1\ra\bbz\ra\pi_0\ra Q_0\ra 1$ be the pullback of the above exact 
sequence via $Q_0\hra Q$. But if $[\pi]$ has finite order, one can take 
$Q_0$ so that $[\pi_0]$ has order 0 so that $\pi_0=\bbz\x Q_0$.
We claim that this group does not embed discretely into $\isom_0(\tp)$.
Choose a standard presentation for $Q_0$:
$$
\langle a_1,b_1,\cdots,a_g,b_g \mid
\prod_{i=1}^g [a_i,b_i]=1\rangle.
$$
Since $Q_0\subset\psl$, we may think of the $a_i$, $b_i$ as elements of
$\psl$.
In $\isom_0(\tp)$, these elements lift to 
$\{(a_i',t_{a_i})\}$, $\{(b_i',t_{b_i})\}$, were $a_i', b_i'\in\tp$;
$t_{a_i},t_{b_i}\in\bbr$.
These are unique up to the center of $\tp$. Since
$\prod_{i=1}^g [(a_i',t_{a_i}),(b_i',t_{b_i})]=(t^{2g-2},0)$ by
\cite{rv81-1}, it is non-zero. Since $\pi_0=\bbz\x
Q_0\subset\isom_0(\tp)$, this relation $2g-2$ would have to be 0. This gives a
contradiction and so $[\pi]$ must have infinite order.
\bigskip

Conversely, suppose $1\ra\bbz\ra\pi\ra Q\ra 1$ is exact, where $Q$ is a
cocompact discrete subgroup of $\isom_0(\bbh)=\psl$; and $[\pi]\in
H^2(Q;\bbz)$ has infinite order. 
Then with the natural inclusion $i:\bbz\hra\bbr$ and the induced homomorphism
$i_*: H^2(Q;\bbz)\ra H^2(Q;\bbr)$, $i_*[\pi]$ is the pushout 
$[\bbr\pi]\in H^2(Q;\bbr)$ (see subsection \ref{pushout-ext}), and is non-zero. 
By Lemma \ref{q-tilde}, $[\wt{Q}]\in H^2(Q;\bbr)$ is also
non-zero. Therefore, there exists $\epsilon\in\bbr$ so that
$(\epsilon\circ i)_*[\pi]=\epsilon_*[\bbr\pi]=[\wt{Q}]$. 
This implies that there exists a
homomorphism of $\pi$ into $\wt{Q}$ with the diagram
$$
\begin{array}{ccccccccccc}
1 &\lra &\bbz   &\lra &\pi    &\lra  &Q      &\lra &1\\
  &     &\bdar{i} &&\bda&&\bdar{=}&&\\
1 &\lra &\bbr   &\lra &\bbr\pi &\lra &Q  &\lra &1\\
  &     &\bdar{\epsilon} &&\bda&&\bdar{=}&&\\
1 &\lra &\bbr   &\lra &\wt{Q} &\lra &Q  &\lra &1\\
\end{array}
$$
commutative and with injective vertical maps. Since $\bbz$ and $Q$ acts \pd
ly with compact quotient on $\bbr$ and $\bbh$ respectively, $\pi$ is
cocompact and discrete. This completes the proof.
\end{proof}

\begin{cor}
Let $\rho:Q\ra\psl$ be a discrete cocompact subgroup.
For an extension $1\ra \bbz\ra\pi\ra Q\ra 1$, there exists an injective 
homomorphism $\theta:\pi\ra\isom(\tp)$ so that the diagram
\begin{equation}
\label{psl-pi}
\begin{array}{ccccccccccc}
1 &\lra &\bbz   &\lra &\pi    &\lra  &Q      &\lra &1\\
  &     &\bdar{\epsilon} &&\bdar{\theta}&&\bdar{\rho}&&\\
1 &\lra &\bbr   &\lra &\isom_0(\tp) &\lra &\psl  &\lra &1\\
\end{array}
\end{equation}
commutes if and only if $[\pi]\in H^2(Q;\bbz)$ has infinite order.
\end{cor}

\begin{cor}
{\rm (Structure)}  Let $M$ be a closed orbifold modelled on
$(\isom_0(\tp),\tp)$-geometry. 
Then $M$ is an orientable closed Seifert orbifold over a hyperbolic base
with $e(M)\neq 0$.
\par\noindent
{\rm (Realization)}  Let $M$ be a compact orientable Seifert orbifold
over a hyperbolic base. Then $M$ admits an $(\isom_0(\tp),\tp)$-geometry 
if and only if $e(M)\neq 0$.
\end{cor}

\begin{step}
\label{h1-qr}
{\rm
Consider the commuting diagram (\ref{psl-pi}). With fixed $\epsilon$ and $\rho$,
how many $\theta$'s are there to make the diagram commutative?
Such maps are classified by $H^1(Q;\bbr)$, see Theorem \ref{main-alg}.
Since the action of $Q$ on $\bbr$ is trivial,  
$$
H^1(Q;\bbr)=\bbr^{2g},
$$
where $2g$ is the first Betti number of the group $Q$.
Compare this with
$$
H^i(Q;\on{M}(\bbh,\bbr))=0,\qquad (i\geq 1)
$$
so that, for any extension $1\ra\bbz\ra \pi\ra Q\ra 1$, a homomorphism
$\pi\ra\ttop_{\bbr}(\bbr\x\bbh)$ exists and is unique, up to conjugation by
elements of $\on{M}(\bbh,\bbr)$.
}\end{step}

\subsection{Lorentz Structures and $\wt{\psl}$-Geometry}  
The 
spaces admitting $\wt{\PSL_2\bbr}$-geometry have another interesting geometric
structure.  
Here is a more explicit description of our problem.
Let us denote $\pslt$ by $\pinf$.
The space $\pinf$ is the universal covering group of $\pone=\psl$.
Topologically $\pinf$ is homeomorphic to $\bbr\times\bbh$.  
$$
\begin{array}{ll}
  P_1&=\psl\\
\pinf&=\pslt.
\end{array}
$$
Consider the {\it indefinite} metric of signature $++--$ on $\bbr^4$.
The unit sphere of this space is
\[
\begin{array}{ll}
S^{1,2} &=\left\{(x,y)\mid x,y\in\bbr^2\,,\quad |x|^2-|y|^2=1\right\}\\
        &\approx {\textrm{O}}(2,2)/{\textrm{O}}(1,2)\,.
\end{array}
\]
The linear map of $\bbr^4$ defined by the matrix
$$
\left[\
\begin{array}{rrrrr}
1\ & 0\  &\ 1 &\ \ 0\\
0\ & 1\  &\ 0 &\ \ 1\\
0\ & -1\  &\ 0 &\ \ 1\\
1\ & 0\  &\ -1 &\ \ 0
\end{array}\
\right]
$$
transforms $S^{1,2}$ to

\[
\begin{array}{ll}
P_2&=\left\{(x,y,z,u)\in \bbr^4\mid xu-yz=1\right\}\\
   &=\on{SL}_2\bbr\,.
\end{array}
\]
Thus $P_2$ has a complete {\it Lorentz metric} of signature $+,-,-$ and {\it
constant sectional curvature} $=1$.

A space is called a {\it Lorentz orbifold} if it is the quotient of $\pinf$
by a discrete group of Lorentz isometries acting {\it properly
discontinuously.}  One can show that such a group contains normal subgroups of
finite index which act freely.  A Lorentz orbifold for which the discrete group
acts freely is called a {\it Lorentz space-form}.  The Lorentz structure on a
space-form is non-singular and it has $\pinf$ as its (metric) universal
covering.
Then $1\ra\bbz\ra \pinf\ra P_1\ra 1$ is actually a central extension 
with $\bbz$ being the entire center of $\pinf$.
It turns out that the identity component of $\isom(\pinf)$ is
$(\pinf\times \pinf)/\bbz$ where $\bbz$ is the diagonal central subgroup
corresponding to the center of each of the $\pinf$-factors.  The action
$\pinf\times \pinf$, as isometries, on $\pinf$ is given by
$$(\alp,\beta)\cdot x=\alp\cdot x\cdot\beta^{-1}\,.$$
Moreover, $\isom_0(\pinf)=\pinf\x_{\bbz} \pinf$ has index 4
in $\isom (\pinf)$.  These Lorentz space-forms are analogous to the
complete spherical space-forms in the Riemannian case.  

Let us describe some
obvious ones which turn out to be {\it homogeneous} in the sense that
$\isom(M)$ acts {\it transitively} on $M$.  Take $\pi\subset \pinf\times
e\subset \pinf\x_{\bbz} \pinf$ as a discrete subgroup.  Then
surely the centralizer of $\pi$ in $\isom(\pinf)$,
$C_{\isom(\pinf)}(\pi)$, contains $e\times \pinf$ in
$(\pinf\x_{\bbz}\pinf)$.  In fact, it is exactly $e\times
\pinf$ (unless $\pi\approx\bbz$ and sits in the center of $\pinf\times e$
in $\isom(\pinf)$).  Such groups $\pi$ are classified in \cite{rv81-1} and are
certain Seifert manifolds over a hyperbolic base, for there is an obvious $S^1$
action on $\pi\bs \pinf$ induced from $e\times
\pinf\subset\isom_0(\pinf)$.  A surprising fact is that {\it all
homogeneous Lorentz orbifolds are} actually {\it homogeneous Lorentz space
forms and coincide with those just described above},
\cite[\S 10]{kr85-1}. 

If $\pi\subset\isom(\pinf)$ so that $M=\pi\bs \pinf$ is compact then it
is shown in \cite[\S7]{kr85-1} that $M$ {\it is homeomorphic to an orientable
Seifert orbifold over a hyperbolic base.}  However, the connections between the
Seifert structure and the Lorentz structure is unclear.  This is due to the
fact that $\isom(\pinf)$ does not act properly on $\pinf$.  By selecting
a maximal subgroup of $\isom(\pinf)$ which acts properly on
$\pinf$, these two disparate structures can be related.

The subgroup 
$$
J(\pinf)=(\pinf\x_{\bbz}\bbr)\rx\bbz_2\subset\isom(\pinf),
$$
where $\pinf\x_{\bbz}\bbr=(\pinf\x\bbr)/\bbz$, where $\bbz$
is the central diagonal subgroup of $\pinf\x\bbr$, and
$\bbz_2$ reverses the orientation of time $(=\bbr)$ and space
$(=\bbh=\Bbb \pinf/\bbr)$ at the same time is the same Lie group as the
$\isom(\pinf)$ in the Riemannian case.   A {\it Lorentz orbifold}
(resp; {\it space-form}) $M=\pi\bs \pinf$, where $\pi\subset
J(\pinf)\subset\isom(\pinf)$ is called {\it standard}.  Therefore the
standard Lorentz orbifolds and space-forms coincide with the orbifolds and
space-forms of the $\wt{\PSL(2,\bbr)}$-geometry.  
A Lorentz space-form is homogeneous if the full group of Lorentz isometries
acts transitively.
It is known that homogeneous Lorentz space-forms $M$
admit nonstandard complete Lorentz structures if $H^1(M;{\bbr})\neq 0$.
See \cite{gold85-1}.

The {\it Seifert fibering}, in the theorem below, on $M$ {\it descends} 
from the $\bbr$-action by the second $\bbr$-factor in $J(\pinf)$ on $\pinf$.  
The following are due to Kulkarni and Raymond and are stated here, 
for simplicity, in the {\it closed cases}.
\medskip

\begin{thm}[{\rm \cite[(8.5)]{kr85-1}}]
\label{lorentz-theorem}
{\rm (Structure)}  A compact standard
Lorentz space-form (resp. Lorentz orbifold) is an orientable Seifert manifold
$M$ (resp. Seifert orbifold) over a hyperbolic base $B$ with $e(M)\neq 0$.

\noindent
{\rm (Realization)}.  Let $M$ be a compact orientable Seifert manifold (resp.
orientable Seifert orbifold) over a hyperbolic base with $e(M)\neq 0$.  Then,
$M$ admits a structure of a standard Lorentz space-form (resp. Lorentz
orbifold).
\end{thm}

We remark that the orbifold part breaks into two separate cases.  If all fibers
are $\approx S^1$, then {\it all} closed orientable Seifert {\it manifolds}
with $e(M)\neq 0$ appear as Lorentz {\it orbifolds} (and conversely).  This is
similar to having the topological sphere appear as a 2-dimensional hyperbolic
orbifold.  If some fibers are arcs then the Lorentz orbifolds are homeomorphic
to connected sums of lens spaces (including $S^3$ and $S^2\times S^1$).  
We should also mention that the statements in Theorem \ref{lorentz-theorem}
are for the full $J(\pinf)$ and not the connected component of the
identity, and correspond to the cocompact discrete subgroups of
$\isom(\tp)$ in the Riemannian case instead of $\isom_0(\tp)$ as described
in the preceding section.
See \cite[\S8,9]{kr85-1} for details and treatment of cases other than the 
compact ones.

\subsection{Deformation Spaces for $\pslt$-geometry}  
\begin{defn}{\rm
Let $\calr(\pi;\calu)$ be the space of all injective homomorphisms
$\theta:\pi\ra \calu$ such that
$\theta(\pi)$ is cocompact acting properly on $P_\infty$ and, 
is discrete in $\calu$.
We topologize $\calr(\pi;\calu)$ as a subset of $\calu^{\pi}$.  
%as a subset of $\prod_{\alp\in\pi}\calu_\alp$.  
In general, if $\calu$ is a Lie group, then $\calr(\pi;\calu)$ will be a real 
analytic space.
}\end{defn}

The space $\calr(\pi;\calu)$ is called \emph{the space of discrete
representations of $\pi$ into $\calu$} or \emph{the Weil
space of $(\pi;\calu)$}.
When there is no confusion likely, we denote $\calr(\pi;\calu)$ simply by
$\calr(\pi)$.
\bigskip

Recall that $\mu$ denotes conjugation.
The inner-automorphisms group $\inn(\calu)$ acts on $\calr(\pi)$ from 
the left by $$\mu(u)\cdot\tth=\mu(u)\circ\tth$$
for $u\in\calu$ and $\tth\in\calr(\pi)$.  
Denote the orbit space of this action by
$$
\calt(\pi)=\inn(\calu)\bs\calr(\pi)\,.
$$
It is called the {\it Teichm\"uller space} of $\pi$ (or of $M$).

$\aut(\pi)$ acts on $\calr(\pi)$ from the right by
$$
\tth\cdot f=\tth\circ f
$$
for $\tth\in\calr(\pi)$ and $f\in\aut(\pi)$.  Denote the orbit space of this 
action by
$$
\cals(\pi)=\calr(\pi)/\aut(\pi)\,.
$$
$\cals(\pi)$ is the {\it space of discrete subgroups of $\calu$ each isomorphic to
$\pi$}, or the \emph{Chabauty space}.  

Since the two actions of $\inn(\calu)$ and $\aut(\pi)$ commute 
with each other, $\aut(\pi)$ acts on $\calt(\pi)$, and $\inn(\calu)$ acts on 
$\cals(\pi)$.
$\aut(\pi)$ has an obvious kernel $\inn(\pi)$.  Consequently, we get an action of
$\out(\pi)$ on $\calt(\pi)$.  We denote the orbit space by
$$
\calm(\pi)=\calt(\pi)/\out(\pi)\,.$$
It is called the {\it moduli} (or {\it Riemann}) {\it space} of $\pi$.  It is
also obtained as the orbit space
$$\calm(\pi)=\inn(\calu)\bs \cals(\pi)\,.$$
Summarizing in the form of a commutative diagram of orbit mappings, we have
\begin{equation}
\begin{array}{rlcl}
(\inn(\calu),&\calr(\pi),\aut(\pi)) &\quad\stackrel{\inn(\calu)\bs}\lra\quad
                        &(\calt(\pi),\out(\pi))\\
             &\quad\bda{/\aut(\pi)}&&\quad\bda{/\out(\pi)}\\
(\inn(\calu),&\cals(\pi)) &\stackrel{\inn(\calu)\bs}\lra &\calm(\pi)
\end{array}
\end{equation}

We shall now describe the deformation spaces for closed $M$ which have a 
geometric structure modelled on $(\calu,\tp)=(\isom_0(\tp),\tp)$ 
(or equivalently the {\it standard} Lorentz structures).  
We shall see that these deformation spaces all have Seifert fiberings 
over well-studied deformation spaces of 
discrete cocompact orientation-preserving hyperbolic groups.  
\bigskip

Let $\pi$ be a cocompact discrete subgroup of $\calu=\isom_0(\pinf)$.
We have the central extension
$$
1\ra\bbz\ra \pi\ra Q\ra 1\,,\qquad [\pi]\in H^2(Q;\bbz)
$$
with $Q\subset P$, and having infinite order in $H^2(Q;\bbz)$.
\bigskip

Since our group $\calu$ embeds into $\ttop_{\bbr}(\bbr\times \bbh)$ we have 
{\it topological rigidity} in the {\it strong sense} that if $\tth_1$ and 
$\tth_2$ are two embeddings of cocompact $\pi$ into 
$\ttop_{\bbr}(\bbr\times \bbh)$, then they are conjugate in 
$\ttop_{\bbr}(\bbr\times \bbh)$.
In almost all cases they will not be conjugate in $\calu$.  The
elements of $\calm(\pi)$ then represent the different $\calu$-structures 
on $M$, and we may expect large deformation spaces.
\bigskip

In order to understand $\calr(\pi)$, we need to study $\aut(\pi)$.
We would like to describe $\aut(\pi)$ in terms of $\aut(Q)$.
Since $\bbz$ is characteristic, any automorphism of $\pi$ induces an
automorphism of $Q$.

\begin{defn} Let $\aut(Q(\pi))$ be the image of $\aut(\pi)\ra\aut(Q)$.
That is, 
$$
\aut(Q(\pi))=\{\tthb\in\aut(Q): \exists\tth\in\aut(\pi)\ 
\rm{inducing}\ \tthb\},
$$
the group of automorphisms of $Q$ which can be 
lifted to an automorphism of $\pi$.
\end{defn}
\bigskip

Because $\bbz\subset\bbr$ has the unique isomorphism extension property,
one can form a pushout (see subsection \ref{pushout-ext}) to get $\bbr\pi$ fitting 
the commuting diagram:
$$
\begin{array}{ccccccccccc}
1 &\lra &\bbz  &\lra &\pi    &\lra  &Q &\lra &1\\
  &     &\bda &&\bda&&\bda&&\\
1 &\lra &\bbr  &\lra &\bbr\pi    &\lra  &Q &\lra &1\\
\end{array}
$$

\begin{lem}
\label{2-4-1}
There is a commutative diagram with exact rows and injective vertical maps:
\begin{equation}
\begin{array}{ccccccccccc}
1 &\lra &\hhom(Q,\bbz)   &\lra &\aut(\pi)    &\lra  &\aut(Q(\pi)) &\lra &1\\
  &     &\bda &&\bda&&\bda&&\\
1 &\lra &\hhom(Q,\bbr)   &\lra &\aut(\bbr\pi)    &\lra  &\aut(Q) &\lra &1\\
\end{array}
\end{equation}
\end{lem}

\begin{proof}
The crucial fact for the proof is $H^2(Q;\bbr)=\bbr$.
Since $[\pi]\in H^2(Q;\bbz)$ has infinite order,
$[\bbr\pi]\in H^2(Q;\bbr)$ is non-zero.
Since $\bbr$ is characteristic in $\bbr\pi$, any $f\in\aut(\bbr\pi)$
induces an automorphism $\fb\in\aut(Q)$.
Suppose $\fb=\rm{id}$. Let $\fh:\bbr\ra\bbr$ be the restriction of
$f$. Then $\fh_*[\bbr\pi]=\fb^*[\bbr\pi]=[\bbr\pi]$ since $\fb=\on{id}$. 
Since $[\bbr\pi]$ is non-zero, $\fh=\on{id}$.
Therefore, $f$ is of the form
$$
f(\alpha)=\lambda(\alpha)\cdot\alpha
$$
for some map $\lambda:\pi\ra\bbr$. One easily sees that $\lambda$ factors
through $Q$ and it satisfies the cocycle condition
$$
\lambda(\ol{\alpha}\ol{\beta})=\lambda(\ol{\alpha})+
\ol{\alpha}\lambda(\ol{\beta})
$$ 
for all $\ol{\alpha}$, $\ol{\beta}\in Q$ so that $\lambda\in Z^1(Q,\bbr)$.
However, since $\bbr$ is central in $\bbr\pi$, $\bbr$ is a trivial
$Q$-module so that $Z^1(Q,\bbr)=\hhom(Q,\bbr)$.
Conversely, any such a $\lambda\in \hhom(Q,\bbr)$ yields an automorphism 
$f\in\aut(\bbr\pi)$.
Moreover, $\hhom(Q,\bbr)\cap\aut(\pi)=\hhom(Q,\bbz)$.

Let $\ol{g}\in\aut(Q)$. We would like to find $g\in\aut(\bbr\pi)$ which
induces $\ol{g}$ on $Q$. The automorphism $\ol{g}$ induces an automorphism
$\ol{g}^*: H^2(Q;\bbr)\ra H^2(Q;\bbr)$.
Since $H^2(Q;\bbr)=\bbr$ and $[\bbr\pi]\neq 0$, there is a real number
$\epsilon$ for which $\ol{g}^*[\bbr\pi]=\epsilon[\bbr\pi]$. 
This number $\epsilon$ is non-zero and can be viewed as an automorphism of
$\bbr\pi$ such that
$$
\begin{array}{ccccccccccc}
1 &\lra &\bbr   &\lra &\bbr\pi  &\lra  &Q &\lra &1\\
  &     &\bdar{\epsilon} &&\bda&&\bdar{\ol{g}}&&\\
1 &\lra &\bbr   &\lra &\bbr\pi  &\lra  &Q &\lra &1\\
\end{array}
$$
is commutative. Thus we have shown that $\aut(\bbr\pi)\ra\aut(Q)$ is
surjective. 
\end{proof}
\bigskip

The deformation spaces of $\pi$ will be studied via those of $Q$. To this
end, it is necessary to define the following:

\begin{defn}{\rm
Let $\calr(Q;\calub)$ be the space of all injective homomorphisms
$\tthb(Q)$ such that $\tthb(Q)$ is cocompact and discrete in $\calub$.

There is a left action of $\inn(\calub)$ on $\calr(Q;\calub)$, and also a
right action of $\aut(Q(\pi))$ on $\calr(Q;\calub)$. We the define
$$
\begin{array}{rcl}
\calt(Q;\calub) &\ =\ &\inn(\calub)\bs\calr(Q;\calub)\\
\cals(Q(\pi);\calub) &\ =\ &\calr(Q;\calub)/\aut(Q(\pi))\\
\calm(Q(\pi);\calub) &\ =\ &\calt(Q;\calub)/\out(Q(\pi)).\\
\end{array}
$$
These are the \emph{Teichm\"uller space},  the \emph{restricted Chabauty
space}, and the \emph{restricted moduli space} of $Q$, respectively.
}
\end{defn}

It is known that $\aut(Q(\pi))$ is a subgroup of $\aut(Q)$ of 
finite index, see \cite[Appendix]{klr86-1}.
This implies that
$\cals\left(Q(\pi)\right)=\calr(Q)/\aut(Q(\pi))$ is a finite regular covering 
of $\cals(Q)$.
Also $\calm(Q(\pi))=\calt(Q)/\out(Q(\pi))$, where 
$\out(Q(\pi))=\aut(Q(\pi))/\inn(Q)$.
\bigskip

\begin{lem}
\label{z-1-free}
The action of $\aut(\pi)$ on $\calr(\pi;\calu)$ extends to an action of 
$\aut(\bbr\pi)$ on $\calr(\pi;\calu)$.
The subgroup $Z^1(Q;\bbr)=H^1(Q;\bbr)=\hhom(Q,\bbr)$ acts on 
$\calr(\pi;\calu)$ freely and properly.
Moreover, $\calr(\pi;\calu)/\hhom(Q,\bbr)=\calr(Q;\calub)$.
\end{lem}

\begin{proof}
Clearly, an element $\theta\in\calr(\pi;\calu)$ determines a homomorphism 
$\tht:\bbr\pi\ra\calu$ uniquely.
The action of $\aut(\bbr\pi)$ on $\calr(\pi;\calu)$ is defined as follows:
For $\theta\in\calr(\pi;\calu)$ and $f\in\aut(\bbr\pi)$,
$$
\theta\cdot f=\tht\circ f|_{\pi}.
$$
The image of $\pi$ under $\tht\circ f$ is a discrete cocompact subgroup of
$\calu$ so that $\tht\circ f|_{\pi}\in\calr(\pi;\calu)$.

Suppose $\theta,\theta'\in\calr(\pi;\calu)$ induce the same representation
$\tthb=\tthb'\in\calr(Q;\calub)$. 
Let $\tht,\tht':\bbr\pi\ra\calu$ be the homomorphisms induced from
$\theta,\theta'$ as described above. Since  $\tthb=\tthb'$, the embeddings 
$\tht$ and $\tht'$ are related by
$\tht'(\alpha)=\lambda(\alpha)\tht(\alpha)$ for some map
$\lambda:\bbr\pi\ra\bbr\subset\calu$. Since $\tht$ and $\tht'$ must be
equal on $\bbr$, the map $\lambda$ factors through $Q$, and hence
$$
\lambda: Q\lra\bbr.
$$
The map $\lambda$ satisfies the cocycle condition so that $\lambda\in
\hhom(Q,\bbr)$. Conversely, let $f\in\aut(\bbr\pi)$ inducing the identity on
$Q$. Then by reversing the order of arguments given above, one sees that
$\theta$ and $\theta\circ f$ represent the same element of
$\calr(Q;\calub)$.

Clearly, unless $f$ is the identity, $\theta$ and $\theta\circ f$ will be
different, which shows that the action of $\hhom(Q,\bbr)$ on
$\calr(\pi;\calu)$ is free and proper since the orbit space is Hausdorff.
\end{proof}

\begin{thm}[{\rm \cite[Theorem 2.5]{klr86-1}}]
\label{deform-thm}
Let $\pi$ be a compact orbifold group with
$(\isom_0(\tp),\tp)$--geometry.
Let $g$ be the genus of the base orbifold. Then,\par
\begin{tabular}{rlll}
$\calr(\pi)$ &$=$&$\calr(Q)\x\bbr^{2g}$ &trivial principal
$\bbr^{2g}$-bundle over $\calr(Q)$\\
$\calt(\pi)$ &$=$&$\calt(Q)\x\bbr^{2g}$ &trivial principal
$\bbr^{2g}$-bundle over $\calt(Q)$\\
$\cals(\pi)$ &&&\kern-0.9in 
$T^{2g}$-bundle over $\cals(Q(\pi))$\\
$\calm(\pi)$ &&&\kern-0.9in 
Seifert fiber space over $\calm(Q(\pi))$ with typical fiber $T^{2g}$.\\
\end{tabular}
\par\noindent
Furthermore, $\cals(Q(\pi))$ is a finite sheeted covering of $\cals(Q)$, and
$\calm(Q(\pi))$ is a finite sheeted branched covering of $\calm(Q)$.
\end{thm}

\begin{proof}
Let $\bbz$ be the center of $\pi$.
Then $Q=\pi/\bbz$ is the base orbifold group. 
Hence in Lemma \ref{2-4-1}, $\hhom(Q,\bbr)=H^1(Q;\bbr)=\bbr^{2g}$, and
$$
1\lra\bbr^{2g}\lra\aut(\bbr\pi)\lra\aut(Q)\lra 1
$$
is exact.
By Lemma \ref{z-1-free}, the group $\bbr^{2g}$ acts on $\calr(\pi)$ freely
and properly so that the orbit map becomes a principal bundle
$$
\bbr^{2g}\ra\calr(\pi)\ra \calr(Q).
$$
Since $\bbr^{2g}$ is contractible, its classifying space is a point and
consequently $\calr(\pi)$ splits as
$(\calr(\pi),\bbr^{2g})=(\calr(Q)\x\bbr^{2g},\bbr^{2g})$ equivariantly,
where $\bbr^{2g}$ acts on the second factor as translations.

Since $\bbr$ is the center of $\isom_0(\tp)$, we have 
$\mu(\isom_0(\tp))=\mu(P)$ and 
$$
\calt(\pi)=\mu(\isom_0(\tp))\bs\calr(\pi)=\mu(P)\bs(\calr(Q)\x\bbr^{2g}).
$$
Now $\mu(P)$ acts on $\calr(Q)$ freely and properly with quotient $\calt(Q)$
\cite{ms75-1} which has the homotopy type of the set of two points.
Therefore $\calt(\pi)$ is a product $\calt(Q)\x\bbr^{2g}$.
Moreover, it is known that $\calt(Q)$ is diffeomorphic to two copies of
$\bbr^{6g-6+2p}$, where $p$ is the number of 
non-free orbit types of the $Q$ action on $\bbh$. (This is the number of 
distinct conjugacy classes of maximal finite subgroups of $Q$).

For the space of subgroups
$$
\cals(\pi)=\calr(\pi)/\aut(\pi)\cong(\calr(Q)\x\bbr^{2g})/\aut(\pi),
$$
note that $\aut(\pi)\cap\bbr^{2g}=\bbz^{2g}$ and the quotient
$\aut(\pi)/\bbz^{2g}$, which we called $\aut(Q(\pi))$, is a subgroup of
$\aut(Q)$. By first dividing out by $\bbz^{2g}$, we get
$\cals(\pi)=(\calr(Q)\x T^{2g})/\aut(Q(\pi))$.
Since $\aut(Q(\pi))$ acts on $\calr(Q)$ freely, we have a genuine fibration
$$
T^{2g}\lra\cals(\pi)\lra\cals(Q(\pi)).
$$
The action of $\aut(\pi)$ on $\bbr^{2g}\x\calr(Q)=\calr(\pi)$ is weakly
$\bbr^{2g}$-equivariant, because $\bbz^{2g}=\hhom(Q;\bbz)$ is normal in
$\aut(\pi)$. In other words, 
$$
\aut(\pi)\ \hra\ \ttop_{\bbr^{2g}}(\bbr^{2g}\x\calr(Q))
=\on{M}(\calr(Q),\bbr^{2g})\rx(\gl(2g,\bbr)\x\ttop(\calr(Q)))
$$
so that
$$
\begin{array}{ccccccccc}
1&\lra&\bbz^{2g}&\lra &\aut(\pi) &\lra&\aut(Q(\pi))&\lra &1\\
 &    &\bda     &     &\bda      &    &\bda        &      \\
1&\lra&\on{M}(\calr(Q),\bbr^{2g})
&\lra &\ttop_{\bbr^{2g}}(\bbr^{2g}\x\calr(Q))
&\lra &\gl(2g,\bbr)\x\ttop(\calr(Q))&\lra &1
\end{array}
$$
is commutative.
Thus, the structure group is a subgroup of the affine group of the torus
$T^{2g}\circ\gl(2g,\bbz)$.
Let $\theta\in\calr(\pi)$. Then for $\alpha\in\pi$,
$\theta\cdot\mu(\alpha)=\mu(\theta(\alpha))\circ\theta$. Therefore, on
$\calt(\pi)$,
$\inn(\pi)=\inn(\bbr\pi)=Q$ acts trivially as does $\inn(Q)\cong Q$ on
$\calt(Q)$. 
Consequently, we have \pd\ actions of $\out(\bbr\pi)$ and $\out(Q(\pi))$ on
$\calt(\pi)$ and $\calt(Q)$.

The space of moduli $\calm(\pi)=\calt(\pi)/\out(\pi)$ requires more care.
Recall that $\calr(\pi)=\bbr^{2g}\x\calt(Q)$ with $\bbr^{2g}$ action by
translations on the first factor.
Since
$$
1\lra\bbr^{2g}\lra\out(\bbr\pi)\lra\out(Q)\lra 1
$$
is exact, we also have the commutative diagram:
\begin{equation}
\begin{array}{ccc}
(\calt(\pi)=\bbr^{2g}\x\calt(Q),\out(\bbr\pi))
&\quad\stackrel{/\bbr^{2g}}\lra
&(\calt(Q),\out(Q))\\
\bdar{/\out(\pi)}&&\bdar{/\out(Q(\pi))}\\
\calm(\pi)=\calt(\pi)/\out(\pi) 
&\quad\stackrel{q}\lra &\calm(Q(\pi))
\end{array}
\end{equation}
The actions and maps arise from the embedding
\begin{equation}
\begin{array}{ccccccccc}
1 &\lra &\bbz^{2g} &\lra &\out(\pi) &\lra &\out(Q(\pi)) &\lra &1\\
  &     &\bda      &     &\bda      &     &\bda         &     & \\
1 &\lra &\bbr^{2g} &\lra &\out(\bbr\pi) &\lra &\out(Q) &\lra &1  
\end{array}
\end{equation}
obtained from Lemma \ref{2-4-1} by dividing out the ineffective $Q$.
Note $\inn(\bbr\pi)\cap\bbr^{2g}=1$.
Now as $\out(\pi)$ normalizes $\bbr^{2g}$ and $\bbz^{2g}$ sits in
$\bbr^{2g}$ as a lattice, the mapping $q$ is a Seifert fibering with
typical fiber the torus $\bbr^{2g}/\bbz^{2g}$.
In general, the fibering will not be locally trivial. 
In fact, if $F=(\out(Q(\pi)))_{[\tthb]}$ for some $[\tthb]\in\calt(Q)$,
then the induced extension
$$
1\lra\bbz^{2g}\lra E\lra F\lra 1
$$ 
acts affinely on $\bbr^{2g}\x[\tthb]$ sitting over $[\tthb]$.
The orbit over $[\tthb]$ under $\out(Q(\pi))$ determines a 2-dimensional
hyperbolic orbifold up to isometry.
Over this hyperbolic orbifold is the set $E\bs\bbr^{2g}=F\bs T^{2g}$ of
metric Seifert orbifolds in $\calm(\pi)$ with base this
hyperbolic orbifold.
\end{proof}

\subsection{Polynomial structure on virtually poly-$\bbz$ manifolds}
We have seen in Theorem \ref{6-1-2} that any torsion-free virtually 
poly-$\bbz$ group $\pi$ 
is the fundamental group of a closed $K(\pi,1)$-manifold.  For example, 
see (4.1).  On the other hand, Milnor \cite{miln77-1} has shown that such 
a group $\pi$ can be the fundamental group of a complete affinely flat 
(not necessarily compact) manifold.  
The question arises as to whether one can find a {\it compact}
complete affinely flat manifold with such fundamental group.  It is true when
$\pi$ is 3-step nilpotent \cite{sche74-1}.  

However, counter-examples produced by Benoist \cite{beno92-1} and Burde and
Grunewald \cite{bg93-1} show
that certain compact nilmanifolds do not admit a complete affinely flat 
structure.
Consequently, the question has a negative answer even in the nilpotent case.

Note that one can look at such a complete affinely flat structure as 
a ``polynomial structure of degree 1''.
In a situation where ``polynomial structures of degree 1'' fails to exist, 
the next best structure will be ``polynomial structure of higher degree''.
The main reference for this section is \cite{dil96-1}.
\medskip

A \tf filtration for a torsion-free finitely generated nilpotent group 
$\gam$ is a central series of the form:
\[\Gamma_\ast:\;\;\Gamma_0=1\subset \Gamma_1 \subset \Gamma_2 \subset
\cdots \subset \Gamma_{c-1} \subset \Gamma_c \subset \Gamma_{c+1}=\Gamma\]
for which
\[\Gamma_{i+1}/\Gamma_{i}\cong \Z{k}\mbox{
for $1\leq i \leq c$ and some $k\in \N{}_0$}\]
Moreover, each $\Gamma_i$ can be chosen to be a characteristic subgroup of
$\Gamma$.
We use $K$ to denote the Hirsh number (or rank) of $\Gamma$.
Often, we will also use $K_i=k_i+k_{i+1}+\cdots+k_c$.
It follows that $K=K_1$.
\bigskip

\noindent
\begin{notation} {\rm
(1) $P(\R{K},\R{k})\subset\on{M}(\R{K},\R{k})$ denotes the 
vector space of polynomial mappings $p:\R{K}\rightarrow \R{k}$.  
So $p$ is given by $k$ polynomials in $K$ variables.\par
\noindent
(2) $P(\R{K})\subset\ttop(\R{K})$ will be used to indicate
the set of all polynomial diffeomorphisms of
$\R{K}$, with an inverse which is also a polynomial mapping. 
This is a group where the multiplication is given by composition.
}\end{notation}

It is not hard to verify that
$P(\R{K},\R{k})$ is a ${\rm Aut}(\Z{k}) \times P(\R{K})$-module
and that the resulting semi-direct product group
\[P(\R{K},\R{k})\rx(\aut(\Z{k}) \times P(\R{K}))\subseteq P(\R{K+k})\]
by defining
$\forall\; (p,g,h)\in P(\R{K},\R{k})\rx(\aut(\Z{k}) \times P(\R{K})),
\;\forall \;(x,y)\in \R{k+K}:$
\[ ^{(p,g,h)}(x,y)=(gx-p(h(y)),h(y)).\]
Restrict $\on{M}(\R{K},\R{k})$ to $P(\R{K},\R{k})$
and $\ttop(\R{K})$ to $P(\R{K})$, and we will speak of canonical type
polynomial representations.

From now on, if we speak of a polynomial representation
$\rho : \gam \rightarrow P(\R{K})$ which is of canonical type, 
we mean of canonical type with respect to some \tf central series.

\begin{thm} \label{niluniek}
{\rm (Existence)}
Let $\gam$ be a torsion-free, finitely generated nilpotent
group of rank $K$. For any \tf central series $\gam_\ast$,
there exists a canonical type polynomial representation
$\rho: \gam\ra P(\R{K})$.
\par\noindent
{\rm (Uniqueness)}
If $\rho_1,\rho_2$ are two
such polynomial representations
of $\gam$ into $P(\R{K})$ with respect to the same $\gam_\ast$,
there exists a polynomial map $p\in P(\R{K})$
such that $\rho_2=p^{-1}\circ \rho_1 \circ p$.
Consequently, the manifolds
${\mbox{$\rho_1(\gam)$}}\!\backslash \!\R{K}$ and
${\mbox{$\rho_2(\gam)$}}\!\backslash \!\R{K}$
are ``polynomially diffeomorphic''.
\end{thm}

\begin{proof} 
The fundamental fact that we shall use is the cohomology vanishing
$$
H^i(\gam/\gam_1,P(\R{K_2},\R{k_1}))=0
$$
for $i>0$. See \cite{dil96-1} for a proof.
\par\noindent
{\rm(Existence)}
We will proceed by induction on the nilpotency class $c$ of $\gam$. If $\gam$ is
abelian, then the existence is well known. Now, suppose that $\gam$ is of
class $c>1$ and the existence is guaranteed for lower nilpotency classes.
Using the induction hypotheses, the group $\gam/\gam_1$ can be furnished with
a canonical type polynomial representation
\[\bar{\rho}: \gam/\gam_1 \rightarrow P(\R{K_2}).\]
We obtain an embedding $i:\gam_1\cong \Z{k_1}\rightarrow P(\R{K_2},\R{k_1})$ if
we define $i(z):\R{K_2}\rightarrow \R{k_1}:x\mapsto z$. We are looking for
a map $\rho$ making the following diagram commutative:
\[
\begin{array}{ccccccccc}
1 & \lra & \gam_1 & \lra & \gam         & \lra & \gam/\gam_1 &\lra & 1 \\
  &      &\bdar{i}&    &\bdar{\rho} &      &\bdar{\varphi\x\bar{\rho}} && \\
1 & \lra & P(\R{K_2},\R{k_1}) &\lra & P \rx(A \times P) &\lra 
&\aut(\Z{k_1})\x P(\R{K_2}) &\lra & 1 \\
 & & & & \cap & & & & \\
 & & & & P(\R{K}) & & & &
\end{array} \]
where
$P\rx (A\times P)$ denotes $P(\R{K_2},\R{k_1})\rx(\aut(\Z{k_1}) \times
P(\R{K_2}))$ and
$\varphi: \gam/\gam_1\rightarrow \aut(\Z{k_1})=\aut(\gam_1)$ denotes the
morphism induced by the extension $1\ra \gam_1 \ra \gam \ra \gam/\gam_1 \ra 1 $.\\
The existence of such a map $\rho$ is now guaranteed by the surjectiveness
of $\delta$ in the long exact cohomology sequence
\par\bigskip
\leftline{$
\cdots\lra 0=H^1(\gam/\gam_1,P(\R{K_2},\R{k_1})) \lra  
H^1(\gam/\gam_1,P(\R{K_2},\R{k_1})/\Z{k_1})$}
\rightline{$
\stackrel{\delta}{\lra} H^2(\gam/\gam_1,\Z{k_1}) \lra 
H^2(\gam/\gam_1,P(\R{K_2},\R{k_1}))=0\lra \cdots$}
\bigskip

\noindent
{\rm (Uniqueness)}
Again we proceed by induction on the nilpotency class $c$ of $\gam$. For
$c=1$ the result is again well known. Indeed, two canonical type representations
of a virtually abelian group are even known to be affinely conjugated.

So we suppose that $\gam$ is of class $c>1$ and that the theorem holds for
smaller nilpotency classes.  The representations $\rho_1,\rho_2$ induce
two canonical type polynomial representations
\[\bar{\rho_1},\bar{\rho_2}:\gam/\gam_1\rightarrow \R{K_2}.\]
By the induction hypothesis, there exists a polynomial map
$\bar{q} : \R{K_2}\rightarrow \R{K_2}$ such that
\[\bar{\rho}_2=\bar q^{-1} \circ \bar \rho_1 \circ \bar q .\]
Lift this $\bar q$ to a polynomial map $q$ of $\R{K_1}$ as follows:
\[ \forall x\in \R{k_1}, \forall y \in \R{K_2}:\;q:\R{K_1}\rightarrow \R{K_1}:
(x,y)\mapsto q(x,y)=(x,\bar q(y)). \]

When restricted to $\gam_1$, $\rho_1$ and $\rho_2$ are mapping elements onto
pure translations of $\R{k_1}$.  We know that there exist an affine
mapping $\tilde{A}:\R{k_1}\rightarrow \R{k_1}$ for which
\[\rho_{2|_{\mbox{\small{\R{k_1}}}}}=
\tilde A^{-1} \circ \rho_{1|_{\mbox{\small{\R{k_1}}}}} \circ \tilde A.\]
Extend $\tilde{A}$ to an affine mapping of $\R{K}$ by defining
\[\forall x\in \R{k_1}, \forall y \in \R{K_2}:\; A:\R{K_1}\rightarrow \R{K_1}:
(x,y)\mapsto A(x,y)=(\tilde A(x),y) \]

Let us denote $\psi=A^{-1}\circ q^{-1}\circ
\rho_1\circ  q\circ  A$. Then we see that $\psi$ and
$\rho_2$ are two canonical type polynomial representations of $\gam$, which
coincide with each other on $\gam_1$ and which induce the same representation
of $\gam/\gam_1$. This means that $\psi$ and $\rho_2$ can be seen as the result
of a Seifert construction with respect to the same data (cf. the commutative
diagram above).

Now we use the injectiveness of $\delta$, which implies that
$\psi$ and $\rho_2$ are conjugated to each other by an element
$r \in P(\R{K_2},\R{k_1})$ (seen as an element of $P(\R{K})$!). So we may
conclude that
\[\rho_2= r^{-1}\circ \psi\circ r
        = r^{-1}\circ A^{-1}\circ q^{-1}\circ \rho_1\circ q\circ A\circ r
        = p^{-1}\circ  \rho_1\circ  p\]
if we take $p=q\circ A \circ r$.
\end{proof}
%---end of section 6 ----------------------------------------------------

\ifx\undefined\bysame
\newcommand{\bysame}{\leavevmode\hbox to3em{\hrulefill}\thinspace}
\fi

%---end of references----------------------------------------------------

\bigskip
\noindent
Copywritten K.B.Lee and Frank Raymond 1999.\par

\vfill\eject
\end{document}